\documentclass[12pt]{article}
\usepackage{fontspec}
\usepackage{amsmath}
\usepackage{amssymb}
\usepackage{amsthm}
\usepackage{enumerate}
\usepackage{graphicx}
\usepackage{float}
\usepackage{indentfirst}
\usepackage[left=1in,right=1in,top=1in,bottom=1in]{geometry}
\usepackage{newpxtext}
\usepackage{pifont}
\usepackage{lipsum}
\usepackage{xunicode, xltxtra, url,parskip}
\usepackage{hyperref}
\usepackage{cleveref}
\usepackage{upgreek}
\usepackage{bbm}
\newcommand{\R}{\mathbb{R}}

\newcommand{\N}{\mathbb{N}}

\newcommand{\CC}{\mathbb{C}}

\newcommand{\n}{\nabla}
\newcommand{\D}{\Delta}

\newcommand{\DD}{\mathbb{D}}
\newcommand{\dd}{\mathrm{d}}
\newcommand{\HH}{\mathbb{H}}
\newcommand{\A}{\mathbb{A}}

\theoremstyle{definition}
\newtheorem{thm}{Theorem}[section]

\newtheorem{defn}{Definition}[section]
\newtheorem{lem}[thm]{Lemma}
\newtheorem{prop}[thm]{Proposition}
\newtheorem{cor}[thm]{Corollary}

\newtheorem{rmk}[thm]{Remark}
\newtheorem{ex}{Example}[section]

\numberwithin{equation}{section}

\newtheorem{introtheorem}{Theorem}

\title{Renormalized volume, Generalized Loewner energy, and Grunsky determinant of multiply connected domains}
\author{Yan Luo\thanks{Department of Mathematics, ETH Zürich, Switzerland. Email: yan.luo@math.ethz.ch}}
\begin{document}
	\maketitle
	\begin{abstract}
		 We define and study three M\"obius invariant functionals for multiply connected planar domains with finitely many boundary components: the Loewner energy, the renormalized volume, and the Grunsky determinant. Using quasiconformal variations, we prove that when the boundary curves are sufficiently regular, the three functionals differ only by a modulus-dependent constant. We also prove that the Loewner energy is finite if and only if every boundary curve is a Weil--Petersson quasicircle. Moreover, we prove that the three functionals all attain their unique extremal value at the circle domains. In particular, the Loewner energy is zero if and only if it is evaluated at a circle domain. Finally, we use the minimization of Loewner energy to give a new proof of the finitely connected version of Koebe uniformization theorem.
	\end{abstract}
    \tableofcontents
	\section{Introduction}
    Let $\Omega\subset\CC$ be a multiply connected domain with $n$ disjoint complements $D_1,\dots, D_n$ each of which is a Jordan domain. 
    Koebe's uniformization theorem \cite{Koe1910} implies that $\Omega$ is conformally equivalent to a \emph{circle domain} $\Omega_c$, namely, a domain whose boundary curves are all round circles. Moreover, $\Omega_c$ is unique up to M\"obius transformations of $\widehat {\mathbb C} = \mathbb C \cup\{\infty\}$. 
The goal of this work is to introduce three natural M\"obius-invariant functionals to measure how different $\Omega$ is from a circle domain and to describe relations between these functionals. 
    
    In the simply connected case, describing a Jordan domain is equivalent to describing its boundary Jordan curve. There are several M\"obius invariant functionals associated with Jordan curves, including Loewner energy \cite{Wang19,RW21}, universal Liouville action \cite{TT06}, M\"obius energy \cite{FHW1994,Bis25}, Fredholm determinants \cite{Sch57,TT06}, renormalized volume in hyperbolic 3-space \cite{BBVW25,Mar26}, free energy of Coulomb gases \cite{Joh22,JV26,wiegmann_zabrodin_22}, etc. Strikingly, except for the M\"obius energy, all of these functionals turn out to be the same up to a multiplicative constant, yet the classes of curves for which each of these functionals is finite all coincide. This class of curves is the Weil--Petersson quasicircles, see Section~\ref{subsec:WP}. Moreover, these functionals all attain their minimum when the curve is a circle.  

    \medskip
    
    
    In this paper we investigate the generalization of three of the functionals mentioned above to multiply connected domains. Note that while all Jordan domains are conformally equivalent, that is, their moduli space is trivial, this is not the case for multiply connected domains and there is \emph{a priori} no reason why the functionals should coincide. 
    

    \textbf{The Loewner energy of a multiply connected domain}\medskip\\
The Loewner energy of a Jordan curve is defined as the Dirichlet energy of a Loewner driving function associated with the curve \cite{RW21}. It originated from the study of the random fractal Schramm--Loewner evolution curves (SLE$_\kappa$). For instance, it is the large deviation rate function as $\kappa \to 0$ of SLE, e.g., \cite{PW24,Guskov_LDP}, and the Onsager--Machlup functional of SLE loop measures \cite{CW24} for fixed $\kappa \le 4$. Wang showed in \cite{Wang19} that the Loewner energy is a multiple of the universal Liouville action introduced in \cite{TT06}.
    
    In \cite{LM25}, we defined the \emph{generalized Loewner potential} $\mathcal{H}_{2,g}$ for an annular domain $A$ with smooth boundary via the zeta-regularized determinant of Laplacians. In this paper, we further generalize this notion to a nonnegative Loewner energy for finitely connected planar domains, with Weil--Petersson boundary regularity.
    \begin{figure}
        \centering
        \includegraphics{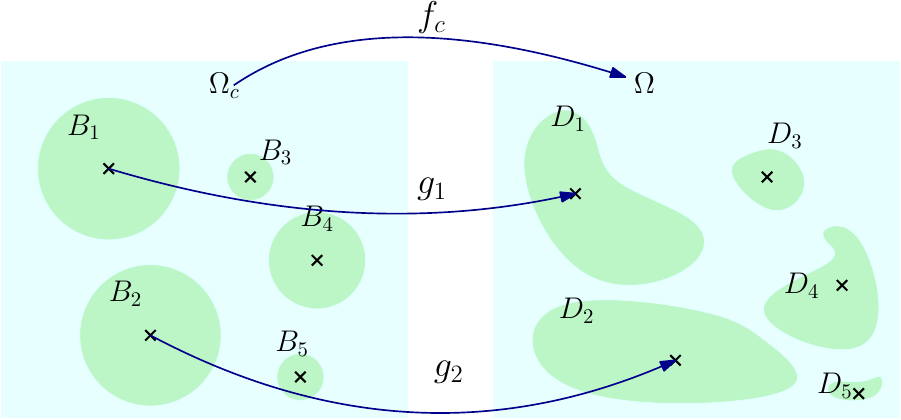}
        \caption{The conformal maps between domains.}
        \label{fig:Multiply}
    \end{figure}


   Let us now discuss the definition of the Loewner energy of a general multiply connected domain. 
    We start by introducing the notations. Let $\Omega\subset \widehat{\CC}$, such that $\infty\in\Omega$, and the complement $\widehat{\CC}\setminus\overline{\Omega}=D_1\cup\dots\cup D_n$ is composed of disjoint bounded Jordan domains in $\widehat{\CC}$. By Koebe's uniformization theorem, there is a circle domain $\Omega_c$, and conformal maps $f_c:\Omega_c\to \Omega$, $g_j:B_j\to D_j$, such that $f_c(z)=z+O(z^{-1})$ near $\infty$. The circle domain $\Omega_c$ conformally equivalent to $\Omega$ is unique up to M\"obius transformations, and we fix it throughout this paper. Let $\widehat{\CC}\setminus\overline{\Omega_c}=B_1\cup\dots\cup B_n$, where $B_j=B_{r_j}(z_j)$ are round disks in $\CC$. We define the \emph{generalized Loewner energy of }$\Omega$ via the following pre-Schwarzian formula similar to the universal Liouville action \cite{TT06} in the single curve case (see \Cref{prop3.2} for details):
    \begin{equation}
        I^L(\Omega)=\frac{1}{\pi}\int_{\Omega_c}\left\vert\frac{f_c''}{f_c'}\right\vert^2 \dd^2z+\frac{1}{\pi}\sum_{j=1}^n\int_{B_j}\left\vert\frac{g_j''}{g_j'}\right\vert^2 \dd^2z+4\sum_{j=1}^n(m_j(g_j)-m_j(f_c)),
    \end{equation}
    where the $m_j$ term is the integral on $\gamma_j$ of the signed geodesic curvature and the log-derivative of the conformal map, and $\dd^2z=\dd x\dd y$ is the Lebesgue measure. Moreover, we have the following properties analogous to the simply connected case \cite{TT06,Wang19,SW24}.
    \begin{introtheorem}\label{thm:A}
        For an $n$-connected domain $\Omega$ introduced as above: 
        \begin{itemize}
            \item[(a)]The Loewner energy $I^L(\Omega)$ is finite if and only if every boundary component of $\partial\Omega$ is a Weil--Petersson quasicircle;
            \item[(b)]The Loewner energy is  $\mathrm{PSL}_2(\CC)$ invariant, nonnegative, with minimum $0$ if and only if $\Omega$ is a circle domain;
            \item[(c)]When the boundary $\partial\Omega$ is smooth, the Loewner energy can be expressed in terms of the zeta-regularized determinants of Laplacians:
            \begin{align}
                I^L(\Omega)&=12(\mathcal{H}_{n,g}(\Omega)-\mathcal{H}_{n,g}(\Omega_c))\notag\\
                &=12\log\frac{\det_\zeta(-\D_{g,\Omega_c})}{\det_\zeta(-\D_{g,\Omega})}\prod_{j=1}^n\frac{\det_\zeta(-\D_{g,B_j})}{\det_\zeta(-\D_{g,D_j})},
            \end{align}
            where $g$ is any metric on $\widehat{\CC}$ conformal to the spherical metric and $\mathcal{H}_{n,g}$ is the Loewner potential defined in \Cref{def:lp};
            \item[(d)]When all the boundary components of $\partial\Omega$ are Weil--Petersson, for any quasiconformal variation $w_t:\CC\to\CC$ driven by \emph{harmonic Beltrami differentials} $\dot{\mu}=\oplus_{j=1}^n\dot{\mu}_j$ supported on $D=\cup_{j=1}^nD_j$, we have:
            \begin{equation}
            \delta I^L(\Omega)=\left.\frac{\dd}{\dd t}\right|_{t=0}I^L(w_t(\Omega))=-\frac{4}{\pi}\mathrm{Re}\sum_{j=1}^n\int_{D_j}\dot{\mu}_j S[g_j^{-1}]\dd^2z,
        \end{equation}
        where $S[g_j^{-1}]$ is the Schwarzian derivative of $g_j^{-1}$.
        \end{itemize}
    \end{introtheorem}
    The detailed proof and discussion can be found in \Cref{sec3}. 

    Let us make a few remarks about the minimization and nonnegativity. 
    Unlike the simply connected case, where $I^L(\Omega) = I^L(\partial \Omega)$ can be equivalently expressed as the Dirichlet energy of the Loewner driving function of the Jordan curve $\partial \Omega$, there is no Loewner evolution adapted to the multiply connected case. So we do not know a priori the nonnegativity of the generalized Loewner energy. In \Cref{sec3.2}, using the cascade relation \Cref{eq43}, we showed that Koebe's construction gives us a non-increasing sequence of domains with smooth boundary for the Loewner energy, which converges to the circle domain $\Omega_c$. By a lower-semicontinuity argument, the circle domain minimizes the energy, which further implies $ I^L \ge 0$. 
    The fact that the energy is non-increasing along the circularization sequence also means intuitively that ``the Loewner energy gets smaller when the domain gets rounder.''
    
    \medskip
    
    \textbf{The Grunsky determinant for multiply connected domains}\medskip\\
    The \emph{Grunsky operator} is a classical tool in geometric function theory. For a Jordan curve $\gamma\subset \CC$ not passing through $\infty$, let $g:\DD^*\to D^*$ be the conformal map with $g(\infty)=\infty,g'(\infty)>0$. Then the Grunsky operator $B$ acting on $\ell^2(\N^*)$ is defined as:
    \begin{equation}
        \log\frac{g(z)-g(w)}{z-w}=\log g'(\infty)-\sum_{k=1}^\infty\sum_{l=1}^\infty a_{kl}z^{-k}w^{-l},\quad z,w\in \DD^*,\quad B :=(\sqrt{kl}a_{kl})_{k,l\in\N^*}.
    \end{equation}
    The curve $\gamma$ is Weil--Petersson if and only if $B$ is Hilbert--Schmidt, and various relations to the Loewner energy of $\gamma$ were obtained  via the Fredholm determinant of $\mathbbm{1}-BB^*$, see \cite{TT06,Joh22,JV26,FVW26}. This determinant was originally considered by Schiffer \cite{Sch57}, in the context of the analysis of Fredholm eigenvalues of the Neumann--Poincaré operator of $\gamma$. Schiffer showed that they, up to sign, coincide with the reciprocal of the singular values of the Grunsky operator of $\gamma$. For a multiply connected domain $\Omega$, Schiffer \cite{Sch59} also considered Fredholm eigenvalues for multiply connected domains, and a corresponding Grunsky operator. 
    This can be easily characterized using conformal maps $\vec{g}=(g_1,\dots,g_n)$, when $g_j:\DD\to D_j$. Let $B_{\vec{g}}$ be an $n\times n$ block operator matrix with each entry being a linear operator on $\ell^2(\N^*)$. The diagonal entries are the $B_1$ component of the corresponding generalized Grunsky matrix (see \cite[Section 2.2.1]{TT06}) of $\gamma_j=\partial D_j$, while the off-diagonal entries have bounded kernels. See \Cref{sec6.1} for details. We also mention that there is a similar construction of multiply connected Grunsky operator in \cite{RSS20}. We prove the following:
    \begin{introtheorem}\label{thm:C}
        For an $n$-connected domain $\Omega$ such that every boundary curve is $C^{3,\alpha}$ with $\alpha>0$, then:
        \begin{itemize}
            \item[(a)]The block Grunsky operator is Hilbert--Schmidt. Consequently, the Grunsky determinant is well-defined, and can be expressed as the Fredholm determinant of $B_{\vec{g}}$:
            \begin{equation}
                \log D_F(\partial\Omega)=\log\det(\mathbbm{1}-B_{\vec{g}}B^*_{\vec{g}}).
            \end{equation}
            \item[(b)]For any quasiconformal variation $w_t:\CC\to\CC$ driven by Beltrami differentials $\dot{\mu}=\oplus_{j=1}^n\dot{\mu}_j$ \emph{compactly supported} on $D=\cup_{j=1}^nD_j$, we have:
            \begin{equation}
                \delta\log D_F(\partial\Omega)=\left.\frac{\dd}{\dd t}\right|_{t=0}\log D_F(w_t(\partial\Omega))=\frac{1}{3\pi}\mathrm{Re}\sum_{j=1}^n\int_{D_j}\dot{\mu}_jS[g_j^{-1}]\dd^2z.
            \end{equation}
            \item[(c)]Consequently,
            \begin{equation}
                \log D_F(\Omega)-\log D_F(\Omega_c)=-\frac{1}{12}I^L(\Omega).
            \end{equation}
        \end{itemize}
    \end{introtheorem}
    We briefly sketch the proof of \Cref{thm:C}. In \cite{Sch59}, Schiffer introduced and derived the point (or Schiffer-) variation of the Grunsky determinant, which deforms the curves with the Cauchy transform of a Dirac point mass. Heuristically, by "integrating" these variations over the support of a Beltrami differential, we may get the quasiconformal variation of the Fredholm determinant. In \Cref{sec6}, we make this argument rigorous by translating the variation vector fields onto Banach spaces, and use the Bochner integral to realize the integration of variations.
    
    We remark on the $C^{3,\alpha}$ regularity assumption used here. In \cite{Sch59}, Schiffer's variation of the Grunsky determinant was studied for domains with only $C^3$ boundary regularity. In our case, however, we used Kellogg's theorem in the proof of \Cref{thm:C}, this is the reason why we impose a slightly stronger H\"older regularity.

    \medskip

    \textbf{The renormalized volume of multiply connected domains}\medskip\\
    The third functional we consider involves the renormalized volume in the hyperbolic $3$-space $\mathbb H^3$. In fact, the group $\mathrm{PSL}_2(\CC)$ of M\"obius transformations of $\widehat{\CC} = \partial_\infty \mathbb H^3$ extends exactly to the group of orientation-preserving isometries of $\mathbb H^3$. Therefore, it is natural to look for a geometric quantity in $\mathbb H^3$ associated to $\Omega$, which is then automatically $\mathrm{PSL}_2(\CC)$-invariant.

    The renormalized volume of a convex cocompact hyperbolic 3-manifold, introduced in \cite{KS08,TT_Liouville}, expresses the Liouville action associated with a conformal metric on the conformal boundary holographically. Taking the conformal metric to be the hyperbolic metric on the boundary, the renormalized volume becomes a functional on the Teichmüller space of the conformal boundary. 
    In \cite{BBVW25}, the authors extended the definition of renormalized volume to the universal Teichm\"uller space (now identified with quasicircles in $\widehat \CC$), and proved that for sufficiently regular curves, it is equal to the universal Liouville action defined in \cite{TT06} up to a multiplicative constant (see also the very recent \cite{Mar26} extending this result to all Weil--Petersson curves). This implies deep relations between the renormalized volume and geometric function theory.\medskip\\
    In this paper, we generalize the construction in \cite{BBVW25} to multiply connected domains assuming the same regularity. In \cite{Ep2024}, Epstein associated each conformal metric on a planar domain with a surface in $\HH^3$. By taking the complete hyperbolic metric $g_U=\rho^2_U|\dd z|^2$ on each component $U$, we get the \emph{Epstein--Poincaré} map $\mathrm{Ep}_U:U\to \HH^3$, where $U\in \{\Omega,D_1,\dots,D_n\}$. The images $\Sigma_\Omega,\Sigma_{D_1},\dots,\Sigma_{D_n}$ of the Epstein maps are called the \emph{Epstein--Poincaré surfaces}. They extend continuously to the boundary of each domain on $\widehat{\CC}=\partial_\infty\HH^3$. We showed that the $n+1$ surfaces here are pairwise disjoint. Thus, they bound a region $M_\Omega\subset\HH^3$. This region is not compact, and the surfaces may have self-intersections (see \Cref{ex:4.2}). We define the (signed) volume of $M_\Omega$ by pulling back the volume form via a diffeomorphism extending the Epstein--Poincar\'e maps. 
    Then the \emph{renormalized volume} of $\Omega$ is defined as:
    \begin{equation}
            V_R(\Omega):=V(\Omega)-\frac{1}{2}\int_{\Sigma_\Omega\cup(\cup_{j=1}^n\Sigma_{D_j})}H\dd A.
        \end{equation}
        Here $V(\Omega)$ is the signed volume of $M_\Omega$, $H=\frac{1}{2}\mathrm{tr}(B)$ is the mean curvature, where $B$ is the corresponding shape operator, and $\dd A$ is the area measure induced from the hyperbolic metric in $\HH^3$. We prove the following results, parallel to the single-curve results in \cite{BBVW25}.
    \begin{introtheorem}\label{thm:B}
        For an $n$-connected domain $\Omega$ such that every boundary curve is  $C^{5,\alpha}$ with $\alpha>0$, then:
        \begin{itemize}
            \item[(a)]The Epstein surfaces $\Sigma_\Omega,\Sigma_{D_1},\dots,\Sigma_{D_n}$ are pairwise disjoint;
            \item[(b)]The renormalized volume $V_R(\Omega)$ is well-defined, finite, and $\mathrm{PSL}_2(\CC)$-invariant. Moreover,
            \begin{equation}
                V_R(\Omega)=V(\Omega)-\frac{1}{2}\int_{\Omega}\frac{|\vartheta_\Omega|^2}{\rho^2_\Omega}\dd^2z-\frac{1}{2}\sum_{j=1}^n\int_{D_j}\frac{|S[g_j^{-1}]|^2}{\rho^2_{D_j}}\dd^2z.
            \end{equation}
            Here $S[g_j^{-1}]$ is the Schwarzian derivative of $g_j^{-1}$, and $\vartheta_\Omega=S[\pi^{-1}]$, where $\pi:\HH\to \Omega$ is a universal covering map. \item[(c)]\label{B.(c)}Let $t\mapsto \mu_t$ be a path in $C^1(I,\mathrm{Bel}^{4,\alpha}(\Omega))$ (see \Cref{eq99} for definition), where $I$ is an interval containing $0$, and $w_t$ is the corresponding normalized solution to Beltrami equation.  Assume that $\mu_t=t(\oplus_{j=1}^n\dot{\mu}_j+\dot{\mu}_\Omega)+o(t)$, which means that $\dot{\mu}_j,\dot{\mu}_\Omega$ are supported on $D_j,\Omega$ respectively. Then:
            \begin{equation}
            \delta V_R(\Omega)=\left.\frac{\dd}{\dd t}\right|_{t=0}V_R(w_t(\Omega))=-\mathrm{Re}\left\{\sum_{j=1}^n\int_{D_j}\dot{\mu}_j\cdot S[g_j^{-1}] \dd^2z+\int_{\Omega}\dot{\mu}_\Omega\cdot\vartheta_{\Omega}\dd^2z\right\}.
        \end{equation}
            \item[(d)]There exists a path of Beltrami differentials $t\mapsto\mu_t\in C^1([0,1],\mathrm{Bel}^{4,\alpha}(\Omega))$, such that the corresponding normalized quasiconformal map $w_t$ is conformal on $\Omega_c$, and $w_0(\Omega)=\Omega_c$, $w_1(\Omega)=\Omega$. Consequently,
            \begin{equation}
                V_R(\Omega)-V_R(\Omega_c)=\frac{\pi}{4}I^L(\Omega).
            \end{equation}
        \end{itemize}
    \end{introtheorem}
    We discuss the detailed proofs in \Cref{sec4} and \Cref{sec5}. We make several brief remarks about this theorem here. First, the technical conditions for the quasiconformal variation in (c) are due to the change of modulus. We allow the Beltrami differential to be supported on the multiply connected domain $\Omega$, so generally the family $\{\Omega_t\}_{t\in I}$ does not belong to the same conformal class. The regularity assumption on Beltrami differentials lets us control the Epstein maps uniformly in time $t$ on $\Omega$ and obtain the final variation formula. Second, we refer the reader to another renormalized volume constructed for multiply connected circle domains with varying conformal metric \cite{BV23}, where they also proved relations to determinants of Laplacians as we have in (d). 
    
    Finally, the assumption on the $C^{5,\alpha}$ regularity of $\partial \Omega$ is the same as in \cite{BBVW25}; we believe an extension to multiply connected domains with the rougher Weil--Petersson boundary is plausible but presents technical difficulties that is not the main focus of this work. 
    For the single-curve case, it was established very recently in \cite{Mar26}, and we think it would be an interesting question for the multiply connected case.
    
    \medskip

    \textbf{The uniformization of planar domains}\medskip\\
   We also give a proof of Koebe's uniformization theorem for an $n$-connected domain $\Omega$, based on our minimization argument of the Loewner potential (\Cref{thm:lpnonnegative}). 
    
    Let $\mathcal{F}(\Omega)$ be all the conformal maps on $\Omega$, such that for $f\in\mathcal{F}(\Omega)$, $f(z)=z+O(z^{-1})$ near $\infty$. For $\Omega_f=f(\Omega)$, a priori, the Loewner potential is well-defined if $\partial \Omega_f$ is smooth. However, we define $\mathcal{H}_{n,g}(\Omega_f)$ to be $\infty$ if $\Omega_f$ is not smooth. We show that for all $f\in\mathcal{F}(\Omega)$, there is a lower bound on the logarithmic capacity of the boundary curves of $\Omega_f$. Then we prove that this lower bound on capacity implies a spectral gap for the Laplace--Beltrami operators, similar to the idea of \cite{Cou87}. After that, using the analysis of zeta-regularized determinant of Laplacians, we may prove a lower bound on $\mathcal{H}_{n,g}(\Omega_f)$ for all $f\in\mathcal{F}(\Omega)$. This lower bound helps us to control the single-curve Loewner energy of every boundary curve in a convergence sequence, which leads to a circle domain as the limit. To sum up, we obtain Koebe's uniformization:
    \begin{introtheorem}
        For any $n$-connected domain $\Omega$, the Loewner potential $\mathcal{H}_{n,g}$ is bounded from below among all domains $\Omega_f=f(\Omega)$ for $f\in\mathcal{F}(\Omega)$. Moreover, there exists $f_0\in\mathcal{F}(\Omega)$, such that $f_0$ maps $\Omega$ to a circle domain conformally.
    \end{introtheorem}
   See \Cref{sec7.1} for detailed discussions. 
\begin{rmk}
   In particular, (b) and (c) of Theorem~\ref{thm:A} then imply that $f_0(\Omega)$ minimizes the Loewner potential. 
\end{rmk}
    
    \medskip
    
 Finally, we discuss possible generalizations of our domain functionals to countably connected domains in \Cref{sec7.2}. Let $\Omega$ be countably connected, and $\widehat{\CC}\setminus\overline{\Omega}=\cup_{j=1}^\infty D_j$ is composed of Jordan domains with $C^{5,\alpha}$ boundary for some $\alpha>0$. We define the renormalized volume of $\Omega$ as the supremum of the Loewner energy of all the domains such that their complement is a finite subset of $\{D_j\}_{j\in\N^*}$. We prove in \Cref{prop7.9} that the critical point of this infinitely connected renormalized volume is a circle domain. We believe it would be interesting to compare our approach to the He--Schramm uniformization theorem \cite{HeSchramm1993} for countably connected domains. 
 
    \bigskip
    
    \textbf{Acknowledgements:} I would like to thank Yilin Wang for many helpful discussions, and her help in the preparation of the manuscript. I also thank Fredrik Viklund for his suggestion to look into Schiffer's work on Fredholm determinants and his comments on the manuscript, and  Baojun Wu for helpful discussions. Part of the work was done at l'Institut des Hautes Études Scientifiques (IHÉS), and I would like to thank them for the support. The author acknowledges support from the Swiss State Secretariat for Education, Research and Innovation (SERI) MB25.00004.
    
    \textbf{AI disclosure:} In the final stage of manuscript preparation, ChatGPT was used for proofreading and literature search.

    \section{Preliminaries} \label{sec:prelim}
    We introduce a few notations which will be used in the whole paper. Let $\widehat{\CC}=\CC\cup\{\infty\}$ be the Riemann sphere. Consider $\Omega\subset\hat{\CC}$ to be an $n$-connected domain containing $\infty$, and $f_c$ is a conformal map from a circle domain $\Omega_c$ to $\Omega$, such that $f_c(z)=z+O(z^{-1})$ near $\infty$. They have complementary domain as $\widehat{\CC}\setminus\overline{\Omega}=D=\cup_{j=1}^n D_j$, and $\widehat{\CC}\setminus\overline{\Omega_c}=B=\cup_{j=1}^nB_j$. We denote $\partial\Omega_c=\cup_{j=1}^nC_j$ and $\partial\Omega=\cup_{j=1}^n\gamma_j$ as disjoint union of round circles or Jordan curves respectively. For a conformal map $f$, it has Schwarzian and pre-Schwarzian derivatives as:
    \begin{equation}
        S[f]=\left(\frac{f''}{f'}\right)'-\frac{1}{2}\left(\frac{f''}{f'}\right)^2,\quad \mathcal{A}[f]=\frac{f''}{f'}.
    \end{equation}
    If $h$ is another conformal map such that $f\circ h$ is defined, then:
    \begin{equation}\label{eq13}
        S[f\circ h]=(S[f]\circ h)\cdot(h')^2+S[h],\quad \mathcal{A}[f\circ h]=(\mathcal{A}[f]\circ h)\cdot h'+\mathcal{A}[h].
    \end{equation}
    There are universal covering maps $\pi_c:\DD\to \Omega_c$ and $\pi:\DD\to \Omega$, such that $f_c\circ \pi_c=\pi$. Locally, let $\vartheta_{\Omega_c}=S[\pi^{-1}_c]$ and $\vartheta_\Omega=S[\pi^{-1}]$. Since the Schwarzian is invariant for $\mathrm{PSL}_2(\CC)$ post-composition, $\vartheta_\Omega$ and $\vartheta_{\Omega_c}$ are well-defined holomorphic quadratic differentials on $\Omega_c$ and $\Omega$ respectively. The domains admit complete hyperbolic metric via the push forward of covering maps of the Poincaré metric on $\DD$. We denote as $g_\Omega=\rho_\Omega^2|\dd z|^2$.
    \subsection{Weil--Petersson Teichmüller space} 
    \label{subsec:WP}
    We present some essential definitions and results on universal Teichmüller theory. Let $U\subset\widehat{\CC}$ be a domain. A Beltrami differential on $U$ is a complex-valued measurable function in the space
    \begin{equation}
        L^\infty_1(U)=\{\mu\in L^\infty(U):\Vert\mu\Vert_\infty<1\}.
    \end{equation}
    By the measurable Riemann mapping theorem \cite{Ah06}, for any $\mu\in L^\infty_1(U)$, there is a unique (up to pre-composition of conformal maps) quasiconformal map $\varphi^{\mu}$ on $\Omega$ satisfying the Beltrami equation:
    \begin{equation}
        \partial_{\bar{z}}\varphi^\mu=\mu\cdot\partial_z\varphi^\mu.
    \end{equation}
    When $\mu\in L^\infty_1(\CC)$, such that $\mu$ is compactly supported, we will use the following normalization very often. Let $\varphi^\mu:\CC\to \CC$ be the solution to Beltrami equation of $\mu$. Since $\mu=0$ near $\infty$, $\varphi^\mu$ is conformal near $\infty$. Thus we may require that $\varphi^\mu(z)=z+O(z^{-1})$ near $\infty$, and such a solution is unique. Moreover, by \cite[Theorem 4.7.4]{Hub06}:
    \begin{prop}\label{prop2.1}
        If $\mu\in L^\infty_1(\CC)$ is compactly supported, then the map:
        \begin{equation}
            \epsilon\mapsto \varphi^{\epsilon\mu},\quad \epsilon\in\CC,\quad \Vert\epsilon\mu\Vert_{L^\infty}<1
        \end{equation}
        is analytic. When $\epsilon$ is small:
        \begin{equation}
            \varphi^{\epsilon\mu}(z)=z+\epsilon V_\mu(z)+o(\epsilon),\quad V_\mu(z)=\frac{1}{\pi}\int_{\CC}\frac{\mu(w)}{z-w}\dd^2w,\quad z\in\CC.
        \end{equation}
    \end{prop}
    \begin{defn}\label{def2.1} Let $\mu\in L^\infty_1(\DD^*)$, there are two ways to extend it to $\hat{\CC}$:
        \begin{itemize}
        \item \textbf{Extension A}: Let $\mu|_\DD=0$, then there is a unique quasiconformal map $w^\mu:\hat{\CC}\to\hat{\CC}$ solving the Beltrami equation, such that $w^\mu(0)=0, (w^\mu)'(0)=1,(w^\mu)''(0)=0$.
        \item \textbf{Extension B}: Let $\mu|_\DD(z)=\overline{\mu(1/\bar{z})}(z/\bar{z})^2$, then there is a unique quasiconformal map $w_\mu:\hat{\CC}\to\hat{\CC}$ solving the Beltrami equation while fixing $1,-1,-i$. Moreover, by the reflection symmetry of $\mu$, $w_\mu\vert_{\mathbb{S}^1}$ is a homeomorphism of $\mathbb{S}^1$.
    \end{itemize}
    \end{defn}
    The universal Teichmüller space $T(1)$ can be defined using the solution to Beltrami equation of the extended differentials.
    \begin{prop}(See \cite[Theorem 3.1.2]{Lehto1987}) Define equivalence relation on $L^\infty_1(\DD^*)$ via:\\
    (a) $\mu\sim\nu$ if $w^\mu|_\DD=w^\nu|_\DD$; (b) $\mu\sim\nu$ if $w_\mu|_{\mathbb{S}^1}=w_\nu|_{\mathbb{S}^1}$.\\
    Then the two quivalence relations are the same. The universal Teichm\"uller space is defined as $T(1)=L^\infty_1(\DD^*)/\sim$. We denote $\Phi:L^\infty_1(\DD^*)\to T(1)$ as the quotient map. 
    \end{prop}
    On a domain $U$ with hyperbolic metric $g_U=\rho_U^2|\dd z|^2$, for $1\le p<+\infty$, we use the following notations of Bergman spaces following \cite{TT06}:
    \begin{align}
        A_\infty(U)&=\left\{\phi\text{ holomorphic on }U:\Vert\phi\Vert_{A_\infty(U)}=\sup_{z\in U}\frac{\vert\phi(z)\vert}{\rho_{U}^2(z)}<+\infty\right\};\\
        A_p(U)&=\left\{\phi\text{ holomorphic on }U:\Vert\phi\Vert_{A_p(U)}=\left(\int_{U}\vert\phi\vert^p\rho^{2-2p}_{U}\dd^2z\right)^{1/p}<+\infty\right\};\\
        A_p^1(U)&=\left\{\phi\text{ holomorphic on }U:\Vert\phi\Vert_{A_p^1(U)}=\left(\int_{U}\vert\phi\vert^p\dd^2z\right)^{1/p}<+\infty\right\}.
    \end{align}
    The universal Teichm\"uller space has a natural group structure given as follows. For $[\mu],[\nu]\in T(1)$, $[\lambda]=[\mu]*[\nu]$ is given by the relation $w_\lambda|_{\mathbb{S}^1}=w_\mu|_{\mathbb{S}^1}\circ w_\nu|_{\mathbb{S}^1}$. The group structure helps to characterize tangent spaces and metrics of $T(1)$. At $[0]\in T(1)$, one can compute that:
    \begin{equation}
        T_{[0]}T(1)\cong L^\infty(\DD^*)/\mathrm{Ker}\{\dd\Phi_0\}\triangleq\Omega^{-1,1}(\DD^*)=\{\dot{\nu}\in L^\infty(\DD^*):\dot{\nu}(z)=\overline{\phi(z)}\rho^{-2}_{\DD^*}(z),\phi\in A_\infty(\DD^*)\}.
    \end{equation}
    For any $[\mu]\in T(1)$, let $R_{[\mu]}$ be the right translation by $[\mu]$, one has $T_{[\mu]}T(1)\cong \dd R_{[\mu]}(T_{[0]}T(1))$. Consider the following subspace of $T_{[0]}T(1)$, which is called the space of harmonic Beltrami differentials on $\DD^*$ (see \cite[Lemma 1.2.1]{TT06}):
    \begin{equation}
        H^{-1,1}(\DD^*)=\{\dot{\nu}\in L^\infty(\DD^*):\dot{\nu}(z)=\overline{\phi(z)}\rho^{-2}_{\DD^*}(z),\phi\in A_2(\DD)\},
    \end{equation}
    By right translations, it gives a distribution $\mathfrak{D}=\{\dd R_{[\mu]}(H^{-1,1}(\DD^*))\}_{[\mu]\in T(1)}$ on the tangent bundle of $T(1)$. Note that for any Jordan domain $D$, we may define $H^{-1,1}(D)$ and $A_2(D)$ using the hyperbolic metric $g_D=\rho_D^2dz^2$.
    \begin{defn} (See \cite[Theorem 1.2.13]{TT06})
        The Weil--Petersson Teichmüller space $T_0(1)$ is the identity component of $T(1)$ in the Hilbert manifold topology. It is also the integral manifold of the distribution $\mathfrak{D}$, so for $[\mu]\in T_0(1)$, $T_{[\mu]}T_0(1)=\mathfrak{D}_{[\mu]}$.
    \end{defn}
    There are various equivalent characterizations of the Weil--Petersson Teichm\"uller space (see \cite{Bis25}). The Weil--Petersson Teichm\"uller space $T_0(1)$ has a unique Kähler structure which is defined through the Weil--Petersson metric:
    \begin{equation}
        \langle\dot{\mu},\dot{\nu}\rangle_{\mathrm{WP}}=\int_{\DD^*}\dot{\mu}\cdot\overline{\dot{\nu}}\rho^2_{\DD^*}\dd^2z,\quad\dot{\mu},\dot{\nu}\in H^{-1,1}(\DD^*)\cong T_{[0]}T_0(1).
    \end{equation}
    Moreover, let $h:\DD^*\to D$ to be a conformal map, then we may define the pull-back differential $h^*\dot{\mu}=(\dot{\mu}\circ h)\frac{\overline{h'}}{h'}$, for $\dot{\mu}\in H^{-1,1}(D)$. This allows us to define the Weil--Petersson metric on any Jordan domain $D$ as:
    \begin{equation}
        \langle\dot{\mu},\dot{\nu}\rangle_{\mathrm{WP},D}=\langle h^*\dot{\mu},h^*\dot{\nu}\rangle_{\mathrm{WP},\DD^*}=\int_{D}\dot{\mu}\cdot\overline{\dot{\nu}}\rho_D^2\dd^2z,\quad \dot{\mu},\dot{\nu}\in H^{-1,1}(D).
    \end{equation}
    The universal Liouville action defined in \cite{TT06} is a K\"ahler potential of the Weil--Petersson metric. An explicit construction is as follows. For $[\mu]\in T_0(1)$, let $D=w^\mu(\DD)$, and $f=w^\mu|_\DD$. Let $g:\DD^*\to D^*$ such that $g(\infty)=\infty$. Then for $\gamma=\partial D$:
    \begin{equation}
        \mathbf{S}([\mu])=\tilde{\mathbf{S}}(\gamma)=\int_\DD\left|\frac{f''}{f'}\right|^2\dd^2z+\int_{\DD^*}\left|\frac{g''}{g'}\right|^2\dd^2z+4\pi\log\left|\frac{f'(0)}{g'(\infty)}\right|.
    \end{equation}
    By \cite[Theorem 8.1]{Wang19}, $\mathbf{S}([\mu])=\pi I^L(\gamma)$, with the notations above. 
    \subsection{Regularity of Jordan curves}\label{sec2.2}
    Let $\gamma\in\widehat{\CC}$ be a Jordan curve. In this paper, we will mainly consider three different regularity notions of $\gamma$. All of them can be characterized by the conformal maps $f:\DD\to D$, where $D$ is a component of $\widehat{\CC}\setminus \gamma$.\medskip\\
    \textbf{The $C^{k,\alpha}$ Jordan curves:}\medskip\\
    Let $\gamma$ be rectifiable, we say it is of $C^{k,\alpha}$ regularity for $k\in \N$ and $0<\alpha<1$, if under the arc length parameterization, $\dd^k\gamma/\dd s^k\in C^{0,\alpha}$. The classical Kellogg--Warschawski theorem states that, $\gamma$ is $C^{k,\alpha}$ if and only if $f'$ does not vanish on $\gamma$, and $f\in C^{k,\alpha}(\overline{\DD})$. This can be further generalized to multiply connected case:
    \begin{prop}\label{prop:kellogg}(Multiply connected Kellogg's theorem, \cite[Theorem 1.1]{HM09})
         Let $f:\Omega_c\to\Omega$ be conformal. For  $k\in\mathbb{N}$, $0<\alpha<1$,  every boundary curve of $\Omega$ is of  $C^{k,\alpha}$ regularity if and only if $f'\ne 0$ on $\partial\Omega_c$ and  $f\in C^{k,\alpha}(\overline{\Omega_c})$.
    \end{prop}
    \textbf{The Weil--Petersson quasicircles:}\medskip\\
    A Jordan curve $\gamma$ is called a quasicircle if there is a quasiconformal map $w:\CC\to\CC$, such that $\gamma=w(\mathbb{S}^1)$.
     A quasicircle $\gamma\in\widehat{\CC}$ is Weil--Petersson if there is $[\mu]\in T_0(1)$, such that up to post-composing a M\"obius map, $\gamma=w^\mu(\mathbb{S}^1)$. The Beltrami differential $\mu$ can be found as follows. Let $f:\DD\to D$ be conformal as defined in the beginning. Since $\gamma$ is a quasicircle, $f$ can be extended to a quasiconformal map $w:\CC\to\CC$. The Beltrami differential $\mu$ of $w$ lies in $L^\infty_1(\DD^*)$. Now take the quasiconformal map $w^\mu$ in \Cref{def2.1}, by uniqueness of solution of Beltrami equation, there is $A\in\mathrm{PSL}_2(\CC)$, such that $w^\mu=A\circ w$.
     \begin{prop}\label{prop2.3}\cite[Theorem 2.1.12]{TT06} Let $\gamma$ be a quasicircle $[\mu]\in T(1)$, and take the Beltrami differential $\mu\in L^\infty_1(\DD^*)$ determined by $\gamma$ as above. Then $\gamma$ is Weil--Petersson if and only if  $S[w^\mu|_\DD]\in A_2(\DD)$.
    \end{prop}
    \textbf{The asymptotically conformal Jordan curves:}\medskip\\
    Let $q \dd z^2$ be a holomorphic quadratic differential on $\Omega$, and $g=\rho^2 |\dd z|^2$ be a conformal metric on $\Omega$. Let $|q|_g(z)=|q(z)|\rho^{-2}(z)$, and use the notation $|q|_\Omega$ when $g=g_\Omega$ is taken to be the hyperbolic metric. Then $\gamma$ is said to be \emph{asymptotically conformal} one of the following holds:
    \begin{equation}
        \lim_{|z|\to1^-}|S[f]|_\DD(z)=0.
    \end{equation}
    See equivalent definitions in \cite[Theorem 11.1]{Pom92}. By \cite[Corollary II.1.4]{TT06}, a Weil--Petersson quasicircle is asymptotically conformal.
    \subsection{Regularity of boundary of multiply connected domains via quadratic differentials}
        For simplification of notations, we write $\vartheta_c$ for $\vartheta_{\Omega_c}$, and $\vartheta$ for $\vartheta_\Omega$. The main technical lemma we are going to use is:
    \begin{lem}\label{lem2.5}\cite[Theorem 3.2]{RSS17}
        Let $q$ be a holomorphic quadratic differential on $\Omega$. For $1\le p\le \infty$, $q\in A_p(\Omega)$ 
        if and only if one of the following holds:
        \begin{itemize}
            \item[(a)] For every $z_0\in\partial\Omega$, there is a chart $\varphi:U\to V\subset\DD$ around $z_0$, such that $(\varphi^{-1})^*q=(q\circ\varphi^{-1})\cdot((\varphi^{-1})')^2\in A_p(V)$ with respect to the hyperbolic metric on $\DD$;
            \item[(b)] For every boundary curve $\gamma_j$ of  $\Omega$, there is a collar chart $\gamma_j\subset A_j\subset\Omega$ with a conformal map $\varphi_j:A_j\to\mathbb{A}(r,1)=\{z\in\CC:r<|z|<1\}$ for some $0<r<1$, such that  $\varphi_j(\gamma_j)=\mathbb{S}^1$, and $(\varphi^{-1})^*q\in A_p(\A(r,1))$ with respect to the hyperbolic metric on $\DD$.
        \end{itemize}
    \end{lem}
    \begin{lem}\label{lem:quadraA2}With the notations above:
         $\vartheta_{c}\in A_2(\Omega_c)$.
    \end{lem}
    \begin{proof}
        By \Cref{lem2.5}, it is equivalent to show that at every point $z\in\partial \Omega_c$, there is a chart $z\in U\subset\Omega_c$, such that $\vartheta_c\in A_2(U)$. To do this, consider $z\in C_j$, and choose $U$ simply connected and small enough, such that $\pi_c^{-1}\vert_U$ is a well-defined conformal map onto a domain $V\subset \DD$. In particular, $\pi_c^{-1}(\partial U\cap C_j)\subset\partial \DD$. Then by Schwarz reflection theorem, $\pi_c^{-1}$ can be conformally extended to a domain across the boundary circle $C_j$. Therefore, $\vartheta_c=S[\pi_c^{-1}]$ extends to be a holomorphic function across $C_j$, and is in particular bounded near $\partial U\cap C_j$. This implies that:
        \begin{equation}\label{eq24}
            \int_U\frac{|S[\pi_c^{-1}]|^2}{\rho_{\Omega_c}^2}\dd^2z=\int_V|S[\pi_c]|^2\frac{(1-|z|^2)^2}{4}\dd^2z\lesssim\int_V\frac{(1-|z|^2)^2}{4}\dd^2z<+\infty.
        \end{equation}
    \end{proof}
    \begin{prop}\label{prop:quadraWP}
        Suppose that every boundary curve of $\Omega$ is a quasicircle. Then we have:
        \begin{itemize}
            \item[(a)]Every boundary curve of $\Omega$ is Weil--Petersson if and only if $\vartheta\in A_2(\Omega)$;
            \item[(b)]Every boundary curve of $\Omega$ is asymptotically conformal if and only if $\lim_{z\to\partial\Omega}|\vartheta|_\Omega(z)=0$.
        \end{itemize}
    \end{prop}
    \begin{proof}(a)
        Let $f_c:\Omega_c\to\Omega$ conformal, such that $f_c\circ \pi_c=\pi$. By \Cref{eq13}, we have that:
        \begin{equation}
            (\vartheta\circ f_c)\cdot(f_c')^2+S[f_c]=\vartheta_c.
        \end{equation}
        Note that $f_c^*\rho_\Omega=\rho_{\Omega_c}$, thus:
        \begin{equation}
            \Vert S[f_c]\Vert_{A_2(\Omega_c)}-\Vert \vartheta_c\Vert_{A_2(\Omega_c)}\le \Vert \vartheta\Vert_{A_2(\Omega)}\le \Vert S[f_c]\Vert_{A_2(\Omega_c)}+\Vert \vartheta_c\Vert_{A_2(\Omega_c)}.
        \end{equation}
        Thus by \Cref{lem:quadraA2}, $\vartheta\in A_2(\Omega)$ if and only if $S[f_c]\in A_2(\Omega_c)$. Now consider one boundary curve $\gamma_j$. By M\"obius maps we may assume that $\gamma_j=f_c(C_j)$, $C_j=\partial\DD$, and $\Omega_c\subset\DD$. For $\gamma_j$, by discussion in \Cref{sec2.2}, we may further assume that there is $\mu\in L^\infty_1(\DD^*)$, such that $\gamma=w^\mu(\mathbb{S}^1)$, and $\Omega\subset w^\mu(\DD)$. Take a small collar $\A_j$ at $C_j$ inside $\Omega_c$, and let $A_j'=(w^\mu)^{-1}\circ f_c(\A_j)$, which is another collar of $C_j$ in $\Omega_c$. Moreover,
        \begin{equation}
            h_j=(w^\mu)^{-1}\circ f_c\vert_{\A_j}:\A_j\to A_j'
        \end{equation}
        is conformal, mapping $\mathbb{S}^1$ to $\mathbb{S}^1$. Thus by Schwarz reflection theorem, it extends conformally to a bigger annular domain containing $\mathbb{S}^1$ in the interior. This implies that $S[h_j]\in A_2(\A_j)$. By \Cref{eq13}
        we have that:
        \begin{equation}
            \Vert S[f_c]\Vert_{A_2(\A_j)}-\Vert S[h_j]\Vert_{A_2(\A_j)}\le \left(\int_{A_j'}\frac{|S[w^\mu]|^2(z)}{((h_j)_*\rho_{\Omega_c})^2(z)}dz^2\right)^{1/2}\le \Vert S[f_c]\Vert_{A_2(\A_j)}+\Vert S[h_j]\Vert_{A_2(\A_j)}.
        \end{equation}
        Thus $S[f_c]\in A_2(\A_j)$ if and only if the integral in the middle is finite. By \cite[Lemma 3.5]{RSS17}, we may regard $h_j:\A_j\to A'$ as a collar chart map of $\Omega_c$ into $\DD$, and then there is $0<K\le 1$, such that:
        \begin{equation}\label{eq29}
            \frac{1}{K}\rho_\DD(z)\le (h_j)_*\rho_{\Omega_c}(z)\le K\rho_\DD(z),\quad z\in A_j'.
        \end{equation} Thus $S[f_c]\in A_2(\A_j)$ if and only if $S[w^\mu]\in A_2(A_j')$. Since our choice of $1\le j\le n$ is arbitrary, by \Cref{prop2.3} and \Cref{lem2.5}, this is equivalent to every boundary curve is Weil--Petersson. 
    \end{proof}
    (b) Consider the same settings as in part (a). 
        By \Cref{lem:quadraA2}, we know that on $\A_j$ $\vartheta_c$ is bounded  thus for $z\in \A_j$,
        \begin{equation}
            \vert\vartheta_c|_{\Omega_c}(z)=\frac{|S[\pi_c](\pi_c^{-1}(z))|}{\rho_\DD^2(\pi_c^{-1}(z))}\lesssim\frac{(1-|\pi_c^{-1}(z)|^2)^2}{4}\to 0
        \end{equation}
        as $z\to \partial\Omega_c$. Then by \Cref{eq13}, on $\A_j$:
        \begin{equation}
            |S[f_c]|_{\Omega_c}(z)-|\vartheta_c|_{\Omega_c}(z)\le |\vartheta|_\Omega(f_c(z))\le |S[f_c]|_{\Omega_c}(z)+|\vartheta_c|_{\Omega_c}(z).
        \end{equation}
        So we have:
        \begin{equation}
            \lim_{z\to C_j}\vert S[f_c]|_{\Omega_c}(z)=0\Leftrightarrow\lim_{z\to\gamma_j}\vert\vartheta|_{\Omega}(z)=0.
        \end{equation}
        Then as in part (a), take $h_j:\A_j\to A_j'$ conformal, by the Schwarz reflection $S[h_j]$ is bounded near $C_j$. Then by \Cref{eq13} and \Cref{eq29} we have that:
        \begin{equation}
            \lim_{z\to C_j}|S[f_c]|_{\Omega_c}(z)=0\Leftrightarrow\lim_{z\to\partial\DD} |S[w^\mu|_\DD]|_\DD(z)=0,
        \end{equation}
        which is equivalent to say that $\gamma_j$ is asymptotically conformal.
        \subsection{Zeta-regularized determinant of Laplacian and Loewner potential}
    Let $(M,g)$ be a compact Riemann surface, such that $M$ is closed or with smooth boundary. The Dirichlet Laplace--Beltrami operator $-\D_g$ has discrete spectrum $0\le\lambda_0\le\lambda_1\le\dots$. Let $\mathrm{tr}(e^{t\D_g})$ be the trace of the heat kernel of $g$, the Mellin transform of it gives:
    \begin{defn}
        The zeta function of $-\D_g$ is given by:
        \begin{equation}\label{eq34}
            \zeta_{-\D_g}(s)=\sum_{\lambda_j>0}\lambda_j^{-s}=\frac{1}{\Gamma(s)}\int_0^\infty t^{s-1}(\mathrm{tr}(e^{t\D_g})-m)\dd t,
        \end{equation}
        where $m\in\{0,1\}$ is the multiplicity of the zero eigenvalue, i.e., $m=0$ for surface with boundary, and $m=1$  for closed surface.
        The zeta function is well-defined through the formula above on $\{\mathrm{Re}(s)>1\}$, and admits meromorphic extension to $\CC$. 
    \end{defn}
    Using the short time expansions of the heat trace, one can show that $\zeta_{-\D_g}$ extends analytically to a neighbourhood of $0\in \CC$ (see e.g.\cite{BGV04}). This ensures that the following definition is well-defined:
    \begin{defn}
        The logarithmic zeta-regularized determinant of $-\D_g$ is defined as $-\zeta'_{-\D_g}(0)$. We use notation $\log{\det}'_\zeta(-\D_g)$ for closed $M$, while $\log{\det}_\zeta(-\D_g)$ for $M$ with boundary. 
    \end{defn}
    The zeta-regularized determinant of Laplacian has nice behaviour under conformal change of the metric.
    \begin{thm}\label{thm:PA}(The Polyakov--Alvarez anomaly formula, \cite{OPS1988}). Let $e^{2\phi}g$ be a Weyl-scaling of $g$, where $\phi\in C^\infty(M)$. Then:
    \begin{itemize}
        \item[(1)] If $M$ is closed:
        \begin{equation}
            \log\frac{\det'_\zeta(-\D_{e^{2\phi}g})}{\det'_\zeta(-\D_g)}=\log\frac{\mathrm{Area}_{e^{2\phi}g}(M)}{\mathrm{Area}_g(M)}-\frac{1}{12\pi}\int_M|\n\phi|^2\dd A_g-\frac{1}{6\pi}\int_MK_g\phi \dd A_g.
        \end{equation}
        \item[(2)] If $M$ has smooth boundary: 
        \begin{align}\label{eq36}
            \log\frac{\det_\zeta(-\D_{e^{2\phi}g})}{\det_\zeta(-\D_g)}&=-\frac{1}{12\pi}\int_M|\n\phi|^2\dd A_g-\frac{1}{6\pi}\int_MK_g\phi \dd A_g\notag\\
            &-\frac{1}{6\pi}\int_{\partial M}k_g\phi \dd s_g-\frac{1}{4\pi}\int_{\partial M}\partial_{\vec{n}}\phi \dd s_g.
        \end{align}
        Here $K_g$ is the Gauss curvature of $M$, and $k_g$ is the geodesic curvature of $\partial M$. On $\partial M$, the outer normal derivative is $\partial_{\vec{n}}$. 
    \end{itemize}
        
    \end{thm}
    The Loewner potential was first defined in \cite{Wang19} for a single smooth Jordan curve in $\widehat{\CC}$, and generalized to pair of disjoint Jordan curves in \cite{LM25}. Here we define the generalized Loewner potential to a multiply connected domain $\Omega$ with smooth boundary in a similar way:
    \begin{defn}\label{def:lp}
        For $\Omega$ with smooth boundary, and $g$ is a metric on $\widehat{\CC}$ which is conformal to the spherical metric $g_{\widehat{\CC}}$. Then we define:
        \begin{equation}\label{eq38}
            \mathcal{H}_{n,g}(\Omega)=\log\frac{\det'_\zeta(-\D_{g,\widehat{\CC}})}{\det_\zeta(-\D_{g,\Omega})\prod_{j=1}^n\det_\zeta(-\D_{g,D_j})\mathrm{Area}_g(\widehat{\CC})}.
        \end{equation}
        In the single Jordan curve case, use the notation $\mathcal{H}_g(\gamma)$, and replace the logarithmic determinant terms in the denominator of \Cref{eq38} by $\det_\zeta(-\D_{g,D})\det_\zeta(-\D_{g,D^*})$. By \cite[Theorem 7.3]{Wang19}, $I^L(\gamma)=12(\mathcal{H}_g(\gamma)-\mathcal{H}_g(\mathbb{S}^1))$.
    \end{defn}
    
    Via the Polyakov--Alvarez formula (\Cref{eq36}), it easily follows that $\mathcal{H}_{n,e^{2\phi}g}=\mathcal{H}_{n,g}$ for any smooth $\phi$. In particular, $\mathcal{H}_{n,g}$ is M\"obius invariant. 
    \begin{prop}\label{prop2.8}
    Let $g_j:B_j(z_j)\to D_j$ be conformal maps, where $B_j(z_j)$ are complements of $\Omega_c$ in $\widehat{\CC}$. Let $f_c:\Omega_c\to\Omega$ be a conformal map. Define the interior-boundary term of  $U\in\{\Omega_c,B_1,\dots,B_n\}$, $h\in\{f_c,g_1,\dots,g_n\}$ as:
    \begin{equation}\label{prop3.2}
        J_{U}(h)=\frac{1}{\pi}\int_{U}\vert\mathcal{A}[h]\vert^2\dd^2z+\frac{2}{\pi}\int_{\partial U}k\cdot\log|h'|\dd s,
    \end{equation}
    where $k$ is the (signed) geodesic curvature of $\partial U$. Then we have:
    \begin{align}\label{eq2.31}
        \mathcal{H}_{n,g}(\Omega)-\mathcal{H}_{n,g}(\Omega_c)&=\frac{1}{12}\sum_{j=1}^nJ_{B_j}(g_k)+\frac{1}{12}J_{\Omega_c}(f_c)\notag\\
        &=\frac{1}{12\pi}\int_{\Omega_c}\left\vert\frac{f_c''}{f_c'}\right\vert^2 \dd^2z+\frac{1}{12\pi}\sum_{j=1}^n\int_{B_j}\left\vert\frac{g_j''}{g_j'}\right\vert^2 \dd^2z+\frac{1}{3}\sum_{j=1}^n(m_j(g_j)-m_j(f_c)),
    \end{align}
    where \begin{equation}
        m_j(h)=\frac{1}{2\pi}\lim_{s\to r_j}\int_0^{2\pi}\log|h'(z_j+se^{i\theta})|\dd\theta.
    \end{equation}
    In particular, when $h=g_j$, $m_j(g_j)=\log|g_j'(z_j)|$.
\end{prop}
\begin{proof}
    Since $\mathcal{H}_{n,g}$ is invariant under Weyl scalings of $g$, we may take a conformal metric $g=e^{2\phi}|\dd z|^2$, such that $\phi=2/(1+|z|^2)$ near $\infty$, and $\phi\equiv 0$ on a big ball $B_R(0)$ containing $\partial\Omega$ and $\partial\Omega_c$ . Then we have that:
    \begin{align}
        \mathcal{H}_{n,g}(\Omega)-\mathcal{H}_{n,g}(\Omega_c)&=\log\frac{\det_\zeta(-\D_{g,\Omega_c})}{\det_\zeta(-\D_{g,\Omega})}+\log\prod_{j=1}^n\frac{\det_\zeta(-\D_{g,B_j})}{\det_\zeta(-\D_{g,D_j})}\notag\\
        &=\log\frac{\det_\zeta(-\D_{g,\Omega_c})}{\det_\zeta(-\D_{f_c^*g,\Omega_c})}+\log\prod_{j=1}^{n}\frac{\det_\zeta(-\D_{\dd z^2,B_j})}{\det_\zeta(-\D_{g_j^*\dd z^2,B_j})}.
    \end{align}
    Then by the Polyakov-Alvarez formula again, we may express the terms above as $L^2$ integrals and boundary integrals. Using Stokes' formula and cancellation observations we get the formula. See \cite[Theorem 3.3]{LM25} for detailed proof.
\end{proof}
If we pick out one boundary curve $\gamma_j$, then we may get the following cascade relation of Loewner potential:
\begin{equation}\label{eq42}
    \mathcal{H}_{n,g}(\Omega)=\mathcal{H}_g(\gamma_j)+\log\frac{\det_\zeta(-\D_{g,\widehat{\CC}\setminus D_j})}{\det_\zeta(-\D_{g,\Omega})\prod_{i\ne j}\det_\zeta(-\D_{g,D_i})}=\mathcal{H}(\gamma_j)+\mathcal{R}_j(\Omega).
\end{equation}
\begin{lem}\label{lem:remainconformal}
    If $\varphi$ is conformal on $\widehat{\CC}\setminus D_j$, then $\mathcal{R}_{g,j}(\varphi(\Omega))=\mathcal{R}_{g,j}(\Omega)$.
\end{lem}
This lemma is an easy consequence of \Cref{thm:PA}. Indeed, note that $\det_\zeta(-\D_{g,\varphi(\Omega)})=\det_\zeta(-\D_{\varphi^*g,\Omega})$. All the corresponding terms in the Polyakov--Alvarez formula will cancel out. 
    \section{Loewner energy of multiply connected domains}\label{sec3}
    
    \subsection{The minimizing problem of generalized Loewner potential}\label{sec3.2}
In this section, we show that the generalized Loewner potential is minimized at the circle domain among all the multiply connected domains with the same modulus. We begin with the following construction:\medskip\\
\textbf{Koebe's construction} (See e.g. \cite[Chapter 17.6]{Hen86}): Let $\Omega$ be an $n$-connected domain, with $\infty\in\Omega$. Fix an order of the exterior domains $D_1,\dots, D_n$. Define a family of domains $\Omega^{(k)}$ inductively as follows. Let $\Omega=\Omega^{(0)}$. Then at step $k\in\N$, take a conformal map $\varphi^{(k)}$ mapping $\widehat{\CC}\setminus D_{j_k}^{(k)}$ to the exterior of a round disk, such that $\varphi^{(k)}(z)=z+O(z^{-1})$ near $\infty$. Here $1+k\equiv j_k(\mathrm{mod} n)$. Then define $\Omega^{(k+1)}=\varphi^{(k)}(\Omega^{(k)})$. By Koebe's uniformization theorem, $(\Omega^{(k)},\infty)\to(\Omega_c,\infty)$ in Carathéodory sense, where $\Omega_c$ is a circle domain. If we also use the notation $\widehat{\CC}\setminus\overline{\Omega_c}=\cup_{j=1}^nB_j$, where $B_j=B_{r_j}(z_j)$, and then $(D_j^{(k)},z_j)\to(B_j,z_j)$.
\begin{thm}\label{thm:lpnonnegative}
    Let $\Omega$ be an $n$-connected domain with smooth boundary, then we have:
    \begin{equation}\label{eq43}
        \mathcal{H}_{n,g}(\Omega)-\mathcal{H}_{n,g}(\Omega_c)\ge\max_{1\le j\le n}\frac{1}{12}I^L(\gamma_j)\ge0
    \end{equation}
    Moreover, $\mathcal{H}_{n,g}(\Omega)=\mathcal{H}_{n,g}(\Omega_c)$ if and only if $\Omega$ itself is a circle domain.
\end{thm}
\begin{proof}
Consider the Koebe's construction introduced as above. At each step, since the curve $\gamma_{j_k}^{(k)}$ is mapped into a round circle, and $\varphi^{(k)}$ is conformal on $\widehat{\CC}\setminus D_{j_k}^{(k)}$ we have by \Cref{lem:remainconformal}:
\begin{align}
    \mathcal{H}_{n,g}(\Omega^{(k+1)})-\mathcal{H}_{n,g}(\Omega^{(k)})&=\mathcal{H}_g(\mathbb{S}^1)+\mathcal{R}_{g,j_k}(\varphi^{(k)}(\Omega^{(k)}))-(\mathcal{H}_g(\gamma_{j_k}^{(k)})+\mathcal{R}_{g,j_k}(\Omega^{(k)}))\notag\\
    &=-\frac{1}{12}I^L(\gamma_{j_k}^{(k)}).
\end{align}
This implies that for any $N\in\N$:
\begin{equation}\label{eq45}
    \mathcal{H}_{n,g}(\Omega)-\mathcal{H}_{n,g}(\Omega^{(N+1)})=\frac{1}{12}\sum_{k=0}^NI^L(\gamma_{j_k}^{(k)})\ge 0.
\end{equation}
 Since $(\Omega^{(k)},\infty)\to(\Omega_c,\infty)$, $(D_j^{(k)},z_j)\to(B_j,z_j)$ in Carathéodory sense, by \cite[Theorem 2.3]{You18}, if we let $\hat{f}_k:\Omega_c\to\Omega^{(k)}$, $\hat{g}_{j,k}:B_i\to D_j^{(k)}$ conformally, then $\hat{g}_{k,j}\to\mathbbm{1}_{B_j}$, $\hat{f}_k\to\mathbbm{1}_{\Omega_c}$ locally uniformly. Using this we can show the following lower-semicontinuity of $\mathcal{H}_{n,g}$. Indeed, consider any compact set $K$ in $\Omega_c$ or $B_j$, we have: 
\begin{equation}
    0=\int_K\lim_{k\to\infty}\left|\frac{\hat{f}_k''}{\hat{f}_k'}\right|^2\dd^2z=\lim_{k\to\infty}\int_K\left|\frac{\hat{f}_k''}{\hat{f}_k'}\right|^2\dd^2z\le\liminf_{k\to\infty}\int_{\Omega_c}\left|\frac{\hat{f}_k''}{\hat{f}_k'}\right|^2\dd^2z.
\end{equation}
By taking an exhaustion of $\Omega_c$, we have that:
\begin{equation}
    0\le\liminf_{k\to\infty}\mathcal{H}_{n,g}(\Omega^{(k)})-\mathcal{H}_{n,g}(\Omega_c).
\end{equation}
This together with \Cref{eq45} implies that:
\begin{equation}
    \mathcal{H}_{n,g}(\Omega_c)\le\liminf_{k\to\infty}\mathcal{H}_{n,g}(\Omega^{(k)})\le\mathcal{H}_{n,g}(\Omega)
\end{equation}
To get \Cref{eq43}, we take the first step of Koebe's  construction for $\gamma_j$ with any $1\le j\le n$, then by \Cref{eq45}:
\begin{equation}
    \mathcal{H}_{n,g}(\Omega)\ge\liminf_{k\to\infty}\mathcal{H}_{n,g}(\Omega^{(k)})+\frac{1}{12}I^L(\gamma_j)\ge\mathcal{H}_{n,g}(\Omega_c).
\end{equation}
Note that $I^L(\gamma_j)=0$ if and only if $\gamma_j$ is a round circle, then we have that $\mathcal{H}_{n,g}(\Omega)=\mathcal{H}_{n,g}(\Omega_c)$ if and only if $\Omega$ is a circle domain.
\end{proof}
\subsection{Definition of Loewner energy and Weil--Petersson finiteness}
\begin{lem}\label{lem:Jinv}
    For any biholomorphism $\varphi$ of $\Omega_c$, we have $J_{\Omega_c}(f_c\circ\varphi)=J_{\Omega_c}(f_c)$.
\end{lem}
\begin{proof}
    Since $\Omega_c$ is a circle domain, $\varphi\in\mathrm{PSL}_2(\CC)$. To simplify the notation, let:
    \begin{equation}
        u=\log|f_c'|,\quad v=\log|\varphi'|.
    \end{equation}
    We compute directly via Stokes formula:
    \begin{align}
        J_{\Omega_c}(f_c\circ \varphi)&=\frac{1}{\pi}\int_{\Omega_c}(|\n(u\circ\varphi)|^2+|\n v|^2+2\langle\n(u\circ\varphi),\n v\rangle)\dd^2z+\frac{2}{\pi}\int_{\partial\Omega_c}k(u\circ\varphi+v)\dd s\notag\\
        &=J_{\Omega_c}(f_c)+\frac{2}{\pi}\int_{\partial\Omega_c}((u\circ\varphi)(\partial_{\vec{n}}v+k)-ku)\dd s+\frac{1}{\pi}\int_{\partial\Omega_c}v(\partial_{\vec{n}}v+2k)\dd s\notag\\
        &=J_{\Omega_c}(f_c)+\frac{1}{\pi}\int_{\partial\Omega_c}(v\circ\varphi^{-1}-v)k\dd s=J_{\Omega_c}(f_c).
    \end{align}
    Here we've used the formula $\partial_{\vec{n}}v=(k\circ\varphi)|\varphi'|-k$, a change of variable, and $v\circ\varphi^{-1}=\log|\varphi'\circ\varphi^{-1}|=-v$. 
\end{proof}
\begin{defn}(Loewner energy of multiply-connected domain).
    For any $\Omega$ conformally equivalent to $\Omega_c$ with smooth boundary, we may define:
    \begin{equation}\label{eq52}
        I^L(\Omega)=12(\mathcal{H}_{n,g}(\Omega)-\mathcal{H}_{n,g}(\Omega_c))\ge 0.
    \end{equation}
    The definition can be extended to the case when $\Omega$ only has Weil--Petersson boundary by \Cref{prop:WPfinite} below:
    \begin{equation}
        I^L(\Omega)=\frac{1}{\pi}\int_{\Omega_c}\left\vert\frac{f_c''}{f_c'}\right\vert^2 \dd^2z+\frac{1}{\pi}\sum_{j=1}^n\int_{B_j}\left\vert\frac{g_j''}{g_j'}\right\vert^2 \dd^2z+4\sum_{j=1}^n(m_j(g_j)-m_j(f_c)).
    \end{equation}
\end{defn}
By \Cref{lem:Jinv}, the definition above does not depend on the choice of conformal maps $f_c,g_1,\dots,g_n$, so it is well-defined.
\begin{prop}\label{prop:WPfinite}
    The Loewner energy $I^L(\Omega)<+\infty$ if and only if $\Omega$ has Weil--Petersson boundary.
\end{prop}
\begin{proof}
    Consider one boundary circle $C_j$, and pick any annulus $A_j=\{z:r_j<s_1<|z-z_j|<s_2\}\subset\Omega_c$ near $C_j$. Then by the holomorphicity of $\log f_c'$ we have that:
    \begin{equation}
        \frac{1}{2\pi}\int_0^{2\pi}\log|f_c'(z_j+s_1e^{i\theta})|\dd\theta=\frac{1}{2\pi}\int_0^{2\pi}\log|f_c'(z_j+s_2e^{i\theta})|\dd\theta.
    \end{equation}
    This implies that $m_j(f_c)$ is well-defined and finite. So we only need to consider the pre-Schwarzian integrals. For $g_j$ this follows from \cite[Theorem 2.1.12]{TT06} that $\gamma_j$ is Weil--Petersson if and only if $\mathcal{A}[g_j]\in L^2(B_j)$. For $f_c$ the proof is essentially the same as \Cref{prop:quadraWP}, just by replacing the Schwarzian with pre-Schwarzian.
\end{proof}
We consider the following approximation construction. Consider any sequence $\{\lambda_k\}_{k\in\N^*}\subset(0,1)$, such that $\lambda_k\to 1$. Let:
\begin{equation}
    B_{k,j}=B_{\lambda_kr_j}(z_j),\quad D_{k,j}=g_j(B_{k,j}),\quad\Omega_k=\widehat{\CC}\setminus(\cup_{j=1}^n\overline{D_{k,j}}).
\end{equation}
Then consider the uniformizing map $f_{k,c}:\Omega_{k,c}\to\Omega_k$, such that $f_{k,c}(z)=z+O(z^{-1})$ near $\infty$. Here $\Omega_{k,c}$ is a circle domain conformally equivalent to $\Omega_k$, and we mention that in general, $\Omega_{k,c}\ne\widehat{\CC}\setminus(\cup_{j=1}^n\overline{B_{k,j}})$
\begin{prop}\label{prop:3.7}
    As $k\to\infty$, $(D_{k,j},g_j(z_j))\to (D_j,g_j(z_j))$, $(\Omega_{k},\infty)\to(\Omega,\infty)$, and $(\Omega_{k,c},\infty)\to(\Omega_c,\infty)$ in Carathéodory sense. Moreover, $f_{k,c}\to f_c$, $f_{k,c}^{-1}\to f_c^{-1}$ locally uniformly. 
\end{prop}
\begin{proof}
    Since $\Omega\subset\Omega_k$ for every $k$, we know that every compact set in $\Omega$ is contained in every $\Omega_k$. On the other hand, if $\infty\in U$ open, such that $U\subset\Omega_k$ for infinitely many $k$. If $U$ is not contained in $\Omega$, then there is one exterior component $D_j\cap U\ne\varnothing$. Pick $z_0\in D_j\cap U$, and we have $d=\mathrm{dist}(g_j^{-1}(z_0),C_j)>0$. But then there is $N\in \N^*$, such that when $k>N$, $r_j-d<\lambda_kr_j$, which implies that $z_0\in D_{k,j}$, and leads to a contradiction. Therefore, $(\Omega_k,\infty)$ converges to $(\Omega,\infty)$. The convergence of $(D_{k,j},g_j(z_j))\to(D_j,g_j(z_j))$ is similar. By \cite[Theorem 2.3, Lemma 3.1]{You18}, we have the convergence of conformal maps as well as the circle domains.
\end{proof}
To proceed, we let $\widehat{\CC}\setminus\overline{\Omega_{k,c}}=\cup_{j=1}^nB_{r_{k,j}}(z_{k,j})$, and $l_{k,j}(z)=r_{k,j}z+z_{k,j}$, $l_j(z)=r_jz+z_j$ for $1\le j\le n$.
\begin{lem}\label{lem3.5}
    There exists $s_0>1$, such that for every $1\le j\le n$:
    \begin{equation}
        \lim_{k\to\infty}r_{k,j}\mathcal{A}[f_{k,c}]\circ l_{k,j}= r_j\mathcal{A}[f_c]\circ l_j\quad\text{ in }L^2(\A(1,s_0)).
    \end{equation}
\end{lem}
\begin{proof}
    Let $h_{k,j}:\DD^*\to D_{k,j}^*$, $h_j:\DD^*\to D_j^*$ be conformal, such that $h_{k,j}(\infty)=h_j(\infty)=\infty$, $h_{k,j}'(\infty)>0$, $h_j'(\infty)>0$. Then since $\partial D^*_{k,j}=g_j(\lambda_kr_j\mathbb{S}^1+z_j)$ forms a family of equipotential curves, by \cite[Appendix A]{TT06}, we have $\mathcal{A}[h_{k,j}]\to\mathcal{A}[h_j]$ in $L^2(\DD^*)$. By the convergence of circle domains, there exists $s>1$, such that for large $k$ and every $1\le j\le n$, the collars $l_{k,j}(\A(1,s))$ are disjoint. Then define  $p_{k,j}=h_{k,j}^{-1}\circ f_{k,c}\circ l_{k,j}$ and $p_j=h_j^{-1}\circ f_c\circ l_j$, which are conformal on $\A(1,s)$, and map $\mathbb{S}^1\to\mathbb{S}^1$. Thus by Schwarz reflection, $p_{k,j}$ extends conformally to $\mathbb{A}(1/s,s)$, and $\{p_{k,j}\}_{k\in\N^*}$ is uniformly bounded. Therefore, by Montel theorem, $\{p_{k,j}\}_{k\in\N^*}$ has a locally uniformly convergent subsequence (which we still denote as $\{p_{k,j}\}_{k\in\N^*}$). This implies that if we take $1<s_0<s$, then $(\mathcal{A}[h_{k,j}]\circ h_{k,j})h_{k,j}'$ converges to $(\mathcal{A}[h_j]\circ p_j)p_j'$ in $L^2(\A(1,s_0))$. Thus by chain rule of pre-Schwarzian, we have:
    \begin{align}
        r_{k,j}\mathcal{A}[f_{k,c}]\circ l_{k,j}
        =(\mathcal{A}[h_{k,j}]\circ p_{k,j})p_{k,j}'+\mathcal{A}[p_{k,j}],
    \end{align}
    which converges to $(\mathcal{A}[h_j]\circ p_j)p_j'+\mathcal{A}[p_j]=r_j\mathcal{A}[f_c]\circ l_j$ in $L^2(\A(1,s_0))$.
\end{proof}
Using the lemma above, we prove the following equipotential approximation of multiply-connected Loewner energy, as an analogue to \cite[Corollary 2.8]{VW2020}.
\begin{prop}\label{prop:leapprox}
    With the notations above, we have that $I^L(\Omega_k)\to I^L(\Omega)$ as $k\to\infty$.
\end{prop}
\begin{proof}
    Let $A_{k,j}=l_{k,j}(\A(1,s_0))$ and $A_j=l_j(\A(1,s_0))$, then by \Cref{lem3.5}, $\Vert \mathcal{A}[f_{k,c}]\Vert_{L^2(A_{k,j})}\to\Vert\mathcal{A}[f_c]\Vert_{L^2(A_j)}$. 
    Outside the collars, we divide the domain into two parts. For $M>0$ large, we define the remaining set $R_k(M)=(\Omega_{k,c}\cap M\DD)\setminus(\cup_{j=1}^nA_{k,j})$, $R(M)=(\Omega_c\cap M\DD)\setminus(\cup_{j=1}^n A_j)$. By the Carathéodory convergence, there is a  compact set $K_M\subset\Omega_c$, such that for any $k$, $R_k(M)\subset K_M$, and $\mathrm{Area}(R_k(M)\triangle R(M))\to 0$ as $k\to\infty$. On $K_M$ we have that $\mathcal{A}[f_{k,c}]$ converges to $\mathcal{A}[f_c]$ uniformly. Therefore, we have, as $k\to\infty$:
    \begin{align}
        &\left|\int_{R_k(M)}\left|\mathcal{A}[f_{k,c}]\right|^2\dd^2z-\int_{R(M)}\left|\mathcal{A}[f_c]\right|^2\dd^2z\right|\notag\\
        &\le\mathrm{Area}(K_M)\sup_{K_M}|\mathcal{A}[f_{k,c}]-\mathcal{A}[f_c]|^2+\sup_{ K_M}|\mathcal{A}[f_c]|^2|\mathrm{Area}(R_{k}(M)\triangle R(M))\to 0.
    \end{align}
    On $M\DD^*$, we let $f_{k,c}(z)=z+\sum_{p=1}^\infty a_{k,p}z^{-p}$, then by area theorem (see e.g. \cite[Theorem 2.1.2]{Lehto1987}), we have $|a_{k,p}|\le M^2$ for every $k$. Thus $|\mathcal{A}[f_{k,c}]|\le CM^2|z|^{-3}$ on $M\DD^*$, where $C$ is independent of $k,M$. Thus:
    \begin{equation}
        \sup_{k\in\N^*}\int_{M\DD^*}|\mathcal{A}[f_{k,c}]|^2\dd^2z\le CM^{-2}\to 0\quad\text{ as }M\to +\infty.
    \end{equation}
    The three parts together show that $\Vert \mathcal{A}[f_{k,c}]\Vert_{L^2(\Omega_{k,c})}\to\Vert\mathcal{A}[f_c]\Vert_{L^2(\Omega_c)}$. For the disk components, the corresponding $L^2$ convergence of conformal maps follows from \cite[Corollary 2.8]{VW2020}, and the convergence of $m_i$ is a direct consequence of the locally uniform convergence. 
\end{proof}
\begin{thm}\label{thm:lenonnegative}The Loewner energy $I^L$ is a well-defined, M\"obius invariant functional for all the finitely connected domains with Weil--Petersson boundary. Moreover, $I^L(\Omega)\ge 0$, with equality if and only if $\Omega$ is a circle domain.
\end{thm}
\begin{proof}
    We take the equipotential approximation ensemble $(\Omega_k,D_{k,1},\dots,D_{k,n})$ as stated above. Since $\Omega_k$ has smooth boundary, by the M\"obius invariance of Loewner potential, we have that for any $A\in\mathrm{PSL}_2(\CC)$, $I^L(A\cdot\Omega_k)=I^L(\Omega_k)$. Note that $\{A\cdot\Omega_k\}_{k\in\N^*}$ forms the equipotential approximation of $A\cdot\Omega$, we have by \Cref{prop:leapprox}:
    \begin{equation}
        I^L(A\cdot\Omega)=\lim_{k\to\infty}I^L(A\cdot \Omega_k)=\lim_{k\to\infty}I^L(\Omega_k)=I^L(\Omega),
    \end{equation}
    which proved the M\"obius invariance for $I^L$ under Weil--Petersson boundary regularity. For the nonnegativity, note that by \Cref{thm:lpnonnegative}, $I^L(\Omega_k)=12(\mathcal{H}_{n,g}(\Omega_k)-\mathcal{H}_{n,g}(\Omega_{k,c}))\ge 0$, which implies that $I^L(\Omega)=\lim_{k\to\infty}I^L(\Omega_k)\ge0$. When $\Omega=\Omega_c$ is a circle domain, every $\Omega_k$ is also a circle domain by construction, thus $I^L(\Omega_k)=0$ and further implies that $I^L(\Omega)=0$. On the other hand, if $I^L(\Omega)=0$, by \Cref{eq43} we have that for each $1\le j\le n$:
    \begin{equation}
        0\le I^L(\gamma_j)=\lim_{k\to\infty}I^L(\gamma_{k,j})\le \lim_{k\to\infty}I^L(\Omega_k)=I^L(\Omega)=0.
    \end{equation}
    This implies that $\gamma_j$ is a round circle for each $j$, and hence $\Omega$ is a circle domain.
\end{proof}
\subsection{Quasiconformal variation of Loewner energy}
   We extend the notion of remainder term $\mathcal{R}_{g,j}(\Omega)$ for a curve $\gamma_j$ defined in \Cref{eq43} to the general Weil--Petersson boundary case:
    \begin{lem}\label{lem:3.10}
        Let
        \begin{equation}
            \mathcal{R}_{g,j}(\Omega)=\frac{1}{12}\left(I^L(\Omega)-I^L(\gamma_j)\right)+\mathcal{H}_{n,g}(\Omega_c)-\mathcal{H}_g(\mathbb{S}^1).
        \end{equation}
        Then $\mathcal{R}_{g,j}(\Omega)$ coincides with \Cref{eq43} when $\Omega$ has smooth boundary. Moreover, for a map $\varphi$ conformal on $\widehat{\CC}\setminus D_j$, we have $\mathcal{R}_{g,j}(\varphi(\Omega))=\mathcal{R}_{g,j}(\Omega)$.
    \end{lem}
    \begin{proof}
        The coincidence at smooth boundary follows from \Cref{eq42}, \Cref{eq52}. In general Weil--Petersson boundary case, denote $D_j^*=\widehat{\CC}\setminus D_j$, with conformal map $h_j:B_j^*\to D_j^*$, $h_j(\infty)=\infty$, we have:
        \begin{equation}
            I^L(\Omega)-I^L(\gamma_j)=\sum_{i\ne j}J_{B_i}(g_i)+J_{\Omega_c}(f_c)-J_{B_j^*}(h_j).
        \end{equation}
        We compute the corresponding term in $I^L(\varphi(\Omega))-I^L(\varphi(\gamma_j))$. Indeed, let $u_i=\log|g_i'|$, $v=\log|\varphi'|$, then:
        \begin{align}
            \pi J_{B_i}(\varphi\circ g_i)&=\int_{B_i}(|\n(v\circ g_i)|^2+|\n u_i|^2+2\langle\n(v\circ g_i),\n u_i\rangle)\dd^2z+2\int_{\partial B_i}k_{\partial B_i}(v\circ g_i+u_i)\dd s\notag\\
            &=\pi J_{B_i}(g_i)+\int_{D_i}|\n v|^2\dd^2z+2\int_{\partial B_i}(v\circ g_i)(k_{\partial D_i}\circ g_i)|g_i'|\dd s\notag\\
            &=\pi J_{B_i}(g_i)+\int_{D_i}|\n v|^2+2\int_{\partial D_i}vk_{\partial D_i}\dd s.
        \end{align}
        Note that this computation also works for the $f_c$ and $h_j$ terms, so we have:
        \begin{align}
        &12\pi(\mathcal{R}_{g,j}(\varphi(\Omega))-\mathcal{R}_{g,j}(\Omega))\notag\\
        &=
            \pi\sum_{i\ne j}\left(J_{B_i}(\varphi\circ g_i)-J_{B_i}(g_i)\right)+\pi\left(J_{\Omega_c}(\varphi\circ f_c)-J_{\Omega_c}(f_c)\right)-\pi\left(J_{B_j^*}(\varphi\circ h_j)-J_{B_j^*}(h_j)\right)\notag\\
            &=\left(\int_{\cup_{i\ne j}D_i}|\n v|^2+\int_{\Omega}|\n v|^2-\int_{D_j^*}|\n v|^2\right)+2\left(\int_{\cup_{i\ne j}\partial D_i}v k_{\partial D_i}+\int_{\partial\Omega}v k_{\partial\Omega}-\int_{\partial D_j^*}v k_{\partial D_j^*}\right)\notag\\
            &=0.
        \end{align}
        The last equality follows from the cancellation of signed curvature terms in terms of the orientation.
    \end{proof}
    \begin{thm}\label{thm:wpvarle}
        For a quasiconformal variation $w_t:\CC\to\CC$ on $\Omega$ with Beltrami differential $\mu_t=t\dot{\mu}+o(t)$, such that $\dot{\mu}=\oplus_{j=1}^n\dot{\mu}_j$, $\dot{\mu}_j\in H^{-1,1}(D_j)$, we have:
        \begin{equation}
            \delta I^L(\Omega)=\left.\frac{d}{dt}\right|_{t=0}I^L(w_t(\Omega))=-\frac{4}{\pi}\mathrm{Re}\sum_{j=1}^n\int_{D_j}\dot{\mu}_jS[g_j^{-1}]\dd^2z.
        \end{equation}
    \end{thm}
    \begin{proof}
        First assume that $\dot{\mu}=\dot{\mu}_j$ for some $1\le j\le n$, i.e., $\mu_t$ is only supported on $D_j$. In this case, $w_t$ is conformal on $D_j^*=\CC\setminus D_j$, by \Cref{lem:3.10}, \cite[Theorem 2.3.5]{TT06}, we have:
        \begin{equation}
            \left.\frac{\dd}{\dd t}\right|_{t=0}I^L(w_t(\Omega))=\left.\frac{\dd}{\dd t}\right|_{t=0}I^L(w_t(\gamma_j))=-\frac{4}{\pi}\mathrm{Re}\int_{D_j}\dot{\mu}_jS[g_j^{-1}]\dd^2z.
        \end{equation}
        Now we consider the general case, when $\dot{\mu}$ is supported on every $D_j$. We consider small real parameters $s=(s_1,\dots, s_n)$. For $1\le k\le n$, let $\mu^{(k)}_s=\sum_{j\le k}s_j\dot{\mu}_j+o(s)$, and suppose $w^{(k)}_s:\CC\to \CC$ quasiconformally, with Beltrami differential $\mu^{(k)}_s$, and normalization  $w_s^{(k)}(z)=z+O(z^{-1})$ at $\infty$. Then we have: $w_s^{(k)}\circ (w_s^{(k-1)})^{-1}$ is conformal except on $D_k$, and its Beltrami differential is of the form:
        \begin{equation}
            s_k\dot{\nu}_{k-1,k}(z)=s_k\left(\dot{\mu}_k\cdot\frac{(w^{(k-1)}_s)'}{\overline{(w^{(k-1)}_s)'}}\right)\circ(w^{(k-1)}_s)^{-1}(z),\quad z\in w_s^{(k-1)}(D_k).
        \end{equation} Thus we have:
        \begin{align}
            I^L(w^{(k)}_s(\Omega))-I^L(w_s^{(k-1)}(\Omega))&=-s_k\frac{4}{\pi}\mathrm{Re}\int_{w^{(k-1)}_s(D_k)}\dot{\nu}_{k-1,k}S[(w_s^{(k-1)}\circ g_k)^{-1}]\dd^2z+o(s_k)\notag\\
            &=-s_k\frac{4}{\pi}\mathrm{Re}\int_{D_k}\dot{\mu}_k(S[g_k^{-1}]-S[w_s^{(k-1)}])\dd^2z+o(s_k).
        \end{align}
        Sum the equalities from $k=1$ to $n$, and by variation of Schwarzian (see e.g. \cite[Theorem 6.12.5]{Hub06}), we have:
        \begin{align}
            I^L(w_s^{(n)}(\Omega))-I^L(\Omega)=-\frac{4}{\pi}\sum_{k=1}^ns_k\mathrm{Re}\int_{D_k}\dot{\mu}_kS[g_k^{-1}]\dd^2z+o(s).
        \end{align}
        Thus $I^L$ is real differentiable with respect to each direction of $\dot{\mu}$. By letting $s_1=\dots=s_n=t$, we get the theorem proved.
    \end{proof}
    \section{The Grunsky determinant of multiply connected domains}\label{sec6}
\subsection{The Grunsky operator and Grunsky determinant}\label{sec6.1}
In this section we assume that $\Omega$ is $n$-connected with $C^{3,\alpha}$ boundary, such that $\infty\in\Omega$. We start with the Green's function $G_D$, where $G_D(z,w)=G_{D_j}(z,w)$ if $z,w\in D_j$, and $G_D(z,w)=0$ if $z,w$ are not in the same simply connected domain. Then consider the following kernel as in \cite{Sch59}:
\begin{equation}
    L_D(z,w)=-\frac{2}{\pi}\frac{\partial^2}{\partial z\partial w}G_D(z,w)=\frac{1}{\pi}\frac{1}{(z-w)^2}-l_D(z,w).
\end{equation}
Moreover, one has:
\begin{lem}\label{lem:kernel}
    Let $g_j:\DD\to D_j$ be conformal, $z,w\in \DD$, then:
    \begin{equation}
        l_D(g_i(z),g_j(w))g_i'(z)g_j'(w)=\frac{1}{\pi}\left(\frac{g_i'(z)g_j'(w)}{(g_i(z)-g_j(w))^2}-\frac{1}{(z-w)^2}\delta_{ij}\right).
    \end{equation}
\end{lem}
In particular,  $l_D(g_j(z),g_j(z))(g_j'(z))^2=\frac{1}{6\pi}S[g_j](z)$. We fix a choice of Riemann maps $\vec{g}=(g_1,\dots,g_n)$ of $D$ from now on. By Kellogg's theorem (\Cref{prop:kellogg}) each $g_j\in C^{3,\alpha}(\overline{\DD})$. So $\vec{g}$ lies in the space:
\begin{equation}
    \mathcal{\mathcal{O}}^{3,n}(\DD)=\{\vec{f}=(f_1,\dots,f_n), f_i\text{ holomorphic on }\DD:\Vert\vec{f}\Vert_{\mathcal{O}^{3,n}}=\sup_{1\le j\le n}\Vert f_j\Vert_{C^3(\overline{\DD})}<+\infty\}.
\end{equation}
By \cite{Sch59}, define the following integral kernel $K^{(ij)}_{\vec{g}}:\DD\times\DD\to \CC$ as:
\begin{equation}
    K^{(ij)}_{\vec{g}}(z,w)=l_D(g_i(z),g_j(w))g'_i(z)g_j'(w).
\end{equation}
If we define the Hilbert space $H_\DD=\oplus_{j=1}^nA_2^1(\DD)$. Then we may proceed the definition of \emph{block Grunsky operator}:
\begin{align}\label{eq:4.5}
    B_{\vec{g}}: \overline{H_\DD}&\to H_\DD,\notag\\
    u&\mapsto\sum_{j=1}^n\int_\DD K^{(ij)}_{\vec{g}}(z,w)\overline{u^{(j)}(w)}\dd^2w,
\end{align}
where $g_i(z)\in D_i$, and $u^{(j)}$ is the $j$-th component of $u\in H_\DD$. If restricts $B_{\vec{g}}$ to the block from the $j-$th component of $H_\DD$ to the $j-$th component of $H_\DD$, it is equal to the $B_4$ Grunsky operator of $D_j$.
\begin{lem}\label{lem6.2}
    The block Grunsky operator $B_{\vec{g}}$ lies in the $2-$Schatten class $\mathcal{S}_2(\overline{H_\DD},H_\DD)$, i.e.,
    \begin{equation}\label{eq:HS}
        \Vert B_{\vec{g}}\Vert^2_{\mathcal{S}_2}=\mathrm{tr}(B_{\vec{g}} B_{\vec{g}}^*)=\sum_{i,j=1}^n\int_\DD\int_\DD|K_{\vec{g}}^{(ij)}(z,w)|^2\dd^2z\dd^2w<+\infty.
    \end{equation}
\end{lem}
\begin{proof}
    The equality follows from expressing each $K^{(ij)}_{\vec{g}}$ according to the basis $\{e_n(z)=\sqrt{\frac{\pi}{n}}z^{n-1}\}_{n\in\N^*}$, and then integrate on $\DD\times\DD$. For the finiteness, since $\partial\Omega$ is $C^{3,\alpha}$, it is also Weil--Petersson. By \cite[Corollary 2.2.9]{TT06}, each diagonal component of the integral above is finite. When $k\ne l$,  we have $|K^{(ij)}_{\vec{g}}(z,w)|\le \frac{1}{\pi}(\mathrm{dist}(D_i,D_j))^{-2}\Vert g_i\Vert_{C^1}\Vert g_j\Vert_{C^1}$, which is also finite.
\end{proof}
\begin{rmk}\label{rmk4.3}
    We mention that for the operator norm of $B_{\vec{g}}$, using \cite[Theorem 4.7]{RSS20}, we may prove that $\Vert B_{\vec{g}}\Vert_{\mathrm{op}}<1$ if every boundary curve is a quasicircle. This is parallel to the simply connected Grunsky operator result. 
\end{rmk}
By \Cref{lem6.2}, $B_{\vec{g}}$ has discrete spectrum. In Schiffer's convention \cite{Sch59}, the inverse of nontrivial eigenvalues of $B_{\vec{g}}$ as $1<\lambda_1\le\lambda_2\le \dots$are called the \emph{Fredholm eigenvalues} of $\Omega$.
Since $B_{\vec{g}}\in\mathcal{S}_2$, the following  (logarithmic) \emph{Grunsky determinant} of $\Omega$ as the Fredholm determinant of $B_{\vec{g}}$ is well-defined:
\begin{equation}
    \log D_{F}(\partial\Omega)=\log\det(\mathbbm{1}-B_{\vec{g}}B_{\vec{g}}^*)=\log\prod_{j=1}^\infty\left(1-\frac{1}{\lambda_j^2}\right).
\end{equation}
\begin{rmk}
    The definition of $D_F(\partial\Omega)$ is independent of the choice of $\vec{g}$. Indeed, if $\vec{\tilde{g}}$ is a different choice of conformal maps, then there are M\"obius maps $\varphi_j\in \mathrm{PSU}(1,1)$, such that $\tilde{g}_j=g_j\circ\varphi_j$. We may define the following unitary operator:
    \begin{align}
        U:H_\DD&\to H_\DD,\notag\\
        u&\mapsto \oplus_{j=1}^n(u^{(j)}\circ \varphi_j)\cdot\varphi_j'.
    \end{align}
    By definition, we have $B_{\vec{\tilde{g}}}=UB_{\vec{g}}\overline{U}^*$, thus:
    \begin{equation}
        \det(\mathbbm{1}-B_{\vec{\tilde{g}}}B_{\vec{\tilde{g}}}^*)=\det(\mathbbm{1}-UB_{\vec{g}}B_{\vec{g}}^*U^*)=\det(\mathbbm{1}-B_{\vec{g}}B_{\vec{g}}^*).
    \end{equation}
\end{rmk}
Schiffer derived the following point variation of the Grunsky determinant of $\Omega$:
\begin{thm}\label{thm:schvar}\cite[Section 7]{Sch59}. Let $z_0\in D_j$, where $D_j$ is one component of $\widehat{\CC}\setminus\overline{\Omega}.$ Then for $|t|<\mathrm{dist}^2(z_0,\partial\Omega)$, the function $\varphi_t(z)=z+t/(z-z_0)$
is conformal on a neighbourhood of $\Omega$. Then it drives a modulus preserving variation:
\begin{equation}
    \left.\frac{\dd}{\dd t}\right|_{t=0}\log D_F(\varphi_t(\partial\Omega))=-2\pi\mathrm{Re}\{l_D(z_0,z_0)\}=\frac{1}{3}\mathrm{Re}\{S[g^{-1}_j](z_0)\},
\end{equation}
\end{thm}
We study a more general type of variation of the Grunsky determinant. To do this, let $\Psi(\vec{g})=\log\det(\mathbbm{1}-B_{\vec{g}}B_{\vec{g}}^*)$. Since $\vec{g}\in\mathcal{O}^{3,n}(\DD)$, $\Vert g_k'\Vert_{C^0}$ is bounded from below for every $1\le k\le n$. Thus take a small neighbourhood $\vec{g}\in W\subset\mathcal{O}^{3,n}(\DD)$, such that every $\vec{f}\in W$ is conformal, and extend to be $C^3$ on the circles. Then $\Psi:W\to \mathcal{S}_2$ is well-defined.
\begin{prop}\label{prop:differential}
    Let $\vec{h}\in T_{\vec{g}}\mathcal{O}^{3,n}(\DD)\cong\mathcal{O}^{3,n}(\DD)$ small, such that $\vec{g}+t\vec{h}\in W$ for $|t|<1$. Then there exists a constant $C\ge0$ only depends on $\vec{g}$, such that: 
    \begin{equation}\label{eq:differentialoriginal}
        |(\dd\Psi)_{\vec{g}}(\vec{h})|=\left|\left.\frac{\dd}{\dd t}\right|_{t=0}\log\det\left(\mathbbm{1}-B_{\vec{g}+t\vec{h}}B_{\vec{g}+t\vec{h}}^*\right)\right|\le C\cdot\Vert \vec{h}\Vert_{C^3}.
    \end{equation}
\end{prop}
\begin{proof}
    To simplify the notation we set $B_t=B_{\vec{g}+t\vec{h}}$. Then $B_t=B_0+t\dot{B}+o(t)$. Let $R=(\mathbbm{1}-B_0B_0^*)^{-1}$, by direct computation:
    \begin{align}
        \Psi(\vec{g}+t\vec{h})&=\log\det(\mathbbm{1}-B_0B_0^*)+\log\det(\mathbbm{1}-tR(B_0\dot{B}^*+\dot{B}B_0^*))+o(t)\notag\\
        &=\Psi(\vec{g})-t\cdot\mathrm{tr}(R(B_0\dot{B}^*+\dot{B}B_0^*))+o(t).
    \end{align}
    Thus let $\Vert B\Vert_{\mathrm{op}}$ be the operator norm of $B$, we have:
    \begin{align}
        |\dd\Psi_{\vec{g}}(\vec{h})|=|\mathrm{tr}(R(B_0\dot{B}^*+\dot{B}B_0^*))|=|2\mathrm{Re}\{\mathrm{tr}(R\dot{B}B_0^*)\}|\le\frac{2\Vert B_0\Vert_{\mathcal{S}_2}}{1-\Vert B_0\Vert_{\mathrm{op}}^2}\Vert\dot{B}\Vert_{\mathcal{S}_2}.
    \end{align}
    We may compute $\dot{B}$ using the kernels $K^{(kl)}$:
    \begin{align}
        \Vert\dot{B}\Vert^2_{\mathcal{S}_2}=\left\Vert\left.\frac{\dd}{\dd t}\right|_{t=0}B_{\vec{g}+t\vec{h}}\right\Vert^2_{\mathcal{S}_2}&=\sum_{k,l=1}^n\int_\DD\int_\DD\left\vert\left.\frac{\dd}{\dd t}\right\vert_{t=0}K^{(kl)}_{\vec{g}+t\vec{h}}\right\vert^2(z,w)\dd^2z\dd^2w.
    \end{align}
    To compute this, we consider the difference quotient as:
    \begin{equation}
        D^Qg_{kl}(z,w)=
    \begin{cases}
        \frac{g_k(z)-g_l(w)}{z-w},\quad &k=l,\\
        g_k(z)-g_l(w),\quad &k\ne l.
    \end{cases}
    \end{equation}
    Then there is a constant $C>0$, such that:
    \begin{equation}
        \Vert D^Qg_{kl}\Vert_{C^0}\le C\Vert\vec{g}\Vert_{C^1}, \Vert(D^Qg_{kl})_z\Vert_{C^0}\le C\Vert \vec{g}\Vert_{C^2},\Vert (D^Qg_{kl})_{zw}\Vert_{C^0}\le C\Vert \vec{g}\Vert_{C^3}.
    \end{equation}
    In the following we fix $1\le k,l\le n$, and to simplify the notation we temporarily drop the subscript $k,l$. Then by \Cref{lem:kernel}, we may compute that:
    \begin{equation}
        K^{(kl)}_{\vec{g}}=\frac{1}{\pi}\left(\frac{(D^Qg)_{zw}}{D^Qg}-\frac{(D^Qg)_z(D^Qg)_w}{(D^Qg)^2}\right).
    \end{equation}
    Therefore, we have:
    \begin{align}
        \left|\left.\frac{\dd}{\dd t}\right\vert_{t=0}K^{(kl)}_{\vec{g}+t\vec{h}}\right|&\le\left|\frac{(D^Qh)_{zw}D^Qg-(D^Qg)_{zw}D^Qh}{\pi(D^Qg)^2}\right\vert+\left\vert\frac{(D^Qh)_z(D^Qg)_w+(D^Qg)_z(D^Qh)_w}{\pi(D^Qg)^2}\right|\notag\\
        &+\left|\frac{2D^Qh(D^Qg)_z(D^Qg)_w}{\pi(D^Qg)^3}\right|\notag\\&\le C(\vec{g})\Vert \vec{h}\Vert_{C^3},
    \end{align}
   which implies \Cref{eq:differentialoriginal} immediately.
\end{proof}
\subsection{Quasiconformal variation of Grunsky determinant}
First fix a compact set $F_j\subset D_j$ for each $1\le j\le n$, and let $F=\cup_{j=1}^nF_j$. Choose $\vec{g}$, such that $g_j(0)\notin F_j$ for each $j$. Let $\mu_t\in L^\infty_1(F)$, such that $\mu_t=t\dot{\mu}+o(t)$, $\dot{\mu}=\oplus_{j=1}^n\dot{\mu}_j$, which means that $\dot{\mu}_j$ is supported on $F_j$. We consider quasiconformal deformation $w_t:\CC\to\CC$ as solution to Beltrami equation of $\mu_t$, such that $w_t(z)=z+O(z^{-1})$ at $\infty$. Let $D_{t,j}=w_t(D_j)$, $\Omega_{t}=w_t(\Omega)$, we construct a vector of conformal maps $\vec{g}_{t}=(g_{t,1}\dots,g_{t,n})$, where $g_{t,j}:\DD\to D_{t,j}$, such that $\vec{g}_{t}\in\mathcal{O}^{3,n}(\DD)$.\medskip\\
\textbf{Explicit construction of $\vec{g}_{t}$}: 
Consider for each $1\le j\le n$, $w_t\circ g_j:\DD\to D_{j,\mu_t}$. It has Beltrami differential as the pull-back of $\mu_{t,j}$ by $g_j$ and reflect it to $\DD^*$:
\begin{equation}\label{eq:nu}
    \nu_{t,j}(z)=\begin{cases}
        (g_j^*\mu_{t,j})(z)=(\mu_{t,j}\circ g_j)(z)\frac{\overline{g_j'(z)}}{g_j'(z)},\quad &z\in\DD,\\
        \overline{(g_j^*\mu_{t,j})(\frac{1}{\bar{z}})}\frac{z^2}{\bar{z}^2},\quad &z\in\DD^*.
    \end{cases}
\end{equation}

Then consider the solution to the Beltrami equation of $\nu_{t,j}$ with normalization:
\begin{equation}
    \partial_{\bar{z}}\phi_{t,j}=\nu_{t,j}\cdot\partial_z\phi_{t,j},\quad \phi_{t,j}:\CC\to\CC\text{ fixes }0,1,\infty.
\end{equation}
 Now take $g_{t,j}=w_t\circ g_j\circ \phi_{t,j}^{-1}$, then $g_{t,j}:\DD\to D_{t,j}$ conformally. By Kellogg's theorem (\Cref{prop:kellogg}), $g_j$ extends to be in $C^{3,\alpha}(\overline{\DD})$. Since both $\mu_t$ and $\nu_{j,t}$ are compactly supported, we have that $\phi_{t,j}$ and $w_t$ are analytic near $\mathbb{S}^1$ and $\gamma_j$ respectively. Thus $\vec{g}_{t}\in\mathcal{O}^{3,n}(\DD)$. We analyze the expansion of $\vec{g}_{t}$ as $t\to 0$ in the following.
\begin{lem}\label{lem:expansion}
    When $|t|$ is small, we have:
    \begin{equation}
        g_{t,j}(z)=g_j(z)+t\cdot h_{\mu,j}(z)+o(t),\quad z\in \DD,
    \end{equation}
    where:
    \begin{equation}
        h_{\mu,j}(z)=(V_{\mu,j}\circ g_j)(z)-X_{\mu,j}(z)g'_j(z),\quad V_{\mu,j}(z)=\frac{1}{\pi}\int_\CC\frac{\dot{\mu}_j(w)}{z-w}\dd^2w,
    \end{equation}
    and
    \begin{equation}\label{eq:Xmu}
        X_{\mu,j}(z)=(V_{\nu,j}(z)-V_{\nu,j}(0))-z(V_{\nu,j}(1)-V_{\nu,j}(0))=\frac{z(1-z)}{\pi}\int_\CC\frac{\dot{\nu}_j(w)}{w(z-w)(1-w)}\dd^2w.
    \end{equation}
\end{lem}
\begin{proof}
    By \Cref{prop2.1}, we have:
    \begin{equation}
        \phi_{t,j}(z)=z+t\cdot X_{\mu,j}(z)+o(t).
    \end{equation}
    Suppose $g_j(z)=\sum_{k=0}^\infty a_kz^k$, then:
    \begin{align}
        (g_j\circ \phi_{t,j}^{-1})(z)&=\sum_{k=0}^\infty a_k(z-t\cdot X_{\mu,j}+o(t))^k\notag\\
        &=\sum_{j=0}^\infty a_k(z^k-k\cdot t\cdot X_{\mu,j} z^{k-1})+o(t)\notag\\
        &=g_j(z)-tX_{\mu,j}\cdot g_j'(z)+o(t).
    \end{align}
    Therefore, we have:
    \begin{equation}
        g_{t,j}(z)=w_t((g_j\circ w^{-1}_{t,j})(z))=g_j(z)+t\left((V_{\mu,j}\circ g_j)(z)-X_{\mu,j}g'_j(z)\right)+o(t).
    \end{equation}
\end{proof}
Now we consider the functional: \begin{equation}
    \Theta(\mu)=\left.\frac{\dd}{\dd t}\right|_{t=0}\log\det(\mathbbm{1}-B_{\vec{g}_t}B_{\vec{g}_t}^*),
\end{equation}
where $\vec{g}_t$ is the family of vector of conformal maps constructed above by $\mu_t$. Using \Cref{prop:differential}, we get the following:
\begin{prop}\label{prop:bounded}
    There exists $C\ge 0$ independent of $\mu$, such that:
    \begin{equation}
         |\Theta(\mu)|\le C_F\Vert \dot{\mu}\Vert_{L^\infty(F)}.
    \end{equation}
    Moreover, $\Theta$ is a bounded linear functional on $L^\infty(F)$.
\end{prop}
\begin{proof}
    By \Cref{prop:differential}, $|\Theta(\mu)|\le C\Vert\vec{h}_\mu\Vert_{C^3}$. Since $\vec{h}_\mu$ is holomorphic, we have $\Vert \vec{h}_\mu\Vert_{C^3(\overline{\DD})}\le C\Vert \vec{h}_\mu\Vert_{C^3(\mathbb{S}^1)}$. To show this, take small collars $\gamma_j\subset U_j$, such that $\mathrm{dist}(U_j,F_j)=d_j>0$ for each $j$. By direct computation, we have $\Vert V_{\mu,j}\Vert_{C^3(U_j)}\le C(d_j)\Vert\dot{\mu}_j\Vert_{L^\infty(F_j)}$, which further implies that $\Vert V_{\mu,j}\circ g_j\Vert_{C^3(\mathbb{S}^1}\le C(d_j,g_j)\Vert\dot{\mu}_j\Vert_{L^\infty(F_j)}$. For $X_{\mu,j}$ similarly, by assumption that $g_j(0)\notin F_j$, there exists $0<r_j<R_j<1$, such that $\dot{\nu}_j$ is supported on $\A(r_j,R_j)\cup \A(R_j^{-1},r_j^{-1})\triangleq A_j\cup A_j^*$, then by \Cref{eq:Xmu}, we have:
    \begin{equation}
        \Vert X_{\mu,j}\Vert_{C^3(\mathbb{S}^1)}\le 2|V_{\nu,j}(0)|+|V_{\nu,j}(1)|+\Vert V_{\nu,j}\Vert_{C^3(\mathbb{S}^1)}\le C(r_j,R_j)\Vert\dot{\nu}_j\Vert_{L^1(\CC)}.
    \end{equation}
    For the $L^1$ term, by \Cref{eq:nu} we have:
    \begin{align}
        \Vert\dot{\nu}_j\Vert_{L^1(\CC)}&=\int_\DD |g_j^*\dot{\mu}_j|(z)\dd^2z+\int_{\DD^*}|g_j^*\dot{\mu}_j|\left(\frac{1}{\bar{z}}\right)\frac{z^2}{\bar{z}^2}\dd^2z\notag\\
        &\le(1+r_j^{-4})\int_\DD|g_j^*\dot{\mu}_j|(z)\dd^2z\notag\\
        &=(1+r_j^{-4})\int_{D_j}\dot{\mu}_j(w)|(g_j^{-1})'|^2\dd^2w\le C(g_j,r_j)\Vert \dot{\mu}_j\Vert_{L^\infty(F_j)}.
    \end{align}
    Combine the estimates above together, with \Cref{lem:expansion}, we have:
    \begin{equation}
        \Vert h_{\mu,j}\Vert_{C^3(\overline{\DD})}\le C(\Vert V_{\mu,j}\circ g_j\Vert_{C^3(\mathbb{S}^1)}+\Vert g'X_{\mu,j}\Vert_{C^3(\mathbb{S}^1)})\le C\Vert\dot{\mu}_j\Vert_{L^\infty(F_j)}.
    \end{equation}
\end{proof}
To proceed, we mention the notion of \emph{Bochner integration}. Let $(E,\mathcal{E},\mu)$ be a measurable space, $Y$ is a Banach space. A $Y$-valued step function on $E$ is given by $s(x)=\sum_{j=1}^ny_j\mathbbm{1}_{A_j}(x)$, where $\{A_j\}_{j=1}^n\subset \mathcal{E}$ are disjoint, and $\{y_j\}_{j=1}^n\subset Y$ are distinct points. Then we have:
\begin{equation}
    \int_Es(x)\mu(dx):=\sum_{j=1}^n\mu(A_j)y_j.
\end{equation}
For a $Y$-valued measurable function $f:E\to Y$, its Bochner integral is defined via approximation of step functions. There is one remarkable property of the Bochner integration. If $Y_1,Y_2$ are two Banach spaces, and $T:Y_1\to Y_2$ is a bounded linear operator, then for a measurable function $f:E\to Y_1$, one has:
\begin{equation}
    \int_ETf\dd\mu=T\int_Ef\dd\mu.
\end{equation}
\begin{thm}(Quasiconformal variation of the Grunsky determinant, compact support case). With the same notation, 
\begin{equation}
    \left.\frac{\dd}{\dd t}\right|_{t=0}\log D_F(w_t(\partial \Omega))=\frac{1}{3\pi}\mathrm{Re}\sum_{j=1}^n\int_{D_i}\dot{\mu}_j(z) S[g_j^{-1}](z)\dd^2z.\label{eq:qcfd}
\end{equation}
    
\end{thm}
\begin{proof}
    For each $1\le j\le n$, pich an open set $F_j\subset U_j$, and $\overline{U}_j\subset D_j$. Let $\delta_j=\mathrm{dist}(\partial U_j, F_j)>0$. For any $\omega\in F_j$, take small disk $B_r(\omega)$ with $r<\delta_j$. Then consider an extension (similar to \cite[Lemma 1]{Gar75}) of Schiffer's point deformation:
    \begin{equation}
        \varphi_{\mu,t}(z)=\begin{cases}
            z+\frac{t}{z-\omega}\dot{\mu}_j(\omega),\quad& z\in (B_r(\omega))^*\notag\\
            z+t\frac{\overline{z-\omega}}{r^2}\dot{\mu}_j(\omega),&z\in B_r(\omega).
        \end{cases}
    \end{equation}
    Then $\varphi_t:\CC\to\CC$ is quasiconformal when $|t|$ small, and its Beltrami differential $\beta_{\mu,r,t}(z)=\frac{t}{r}\dot{\mu}_j(\omega)\mathbbm{1}_{B_r(\omega)}(z)$ is supported on $\overline{U}_j$. By \Cref{thm:schvar}, we have:
    \begin{align}
        \Theta(\dot{\beta}_{\mu,r})(\omega)=\left.\frac{\dd}{\dd t}\right|_{t=0}\log D_F(\varphi_{\mu,t}(\partial \Omega))=\frac{1}{3\pi}\mathrm{Re}\{\dot{\mu}_j(\omega)S[g_j^{-1}](\omega)\}.
    \end{align}
    Now by \Cref{prop:bounded}, we can consider the Bochner integral, as $r\to 0$:
    \begin{equation}
        \int_{D_j}\Theta(\dot{\beta}_{\mu,r})(\omega)\dd^2\omega=\Theta\left(\int_{D_j}\dot{\beta}_{\mu,r}(\omega)\dd^2\omega\right)=\Theta\left(\frac{1}{r}\dot{\mu}_j*\mathbbm{1}_{B_r(0)}\right)\to \Theta(\dot{\mu}_j).
    \end{equation}
    By the same method as \Cref{thm:wpvarle}, this also works when $\dot{\mu}$ is supported on every $D_i$.
\end{proof}
\begin{thm}
    Let $\Omega$ be $n$-connected with $C^{3,\alpha}$ boundary, then:
    \begin{equation}\label{eq:fdle}
        \log D_F(\Omega)-\log D_F(\Omega_c)=-\frac{1}{12}I^L(\Omega).
    \end{equation}
\end{thm}
\begin{proof}
    Let $f_c:\Omega_c\to\Omega$ be conformal. By construction in \Cref{lem:path},
    there exists a quasiconformal map $w:\CC\to\CC$, $w(z)=z+O(z^{-1})$ near $\infty$, and $w|_{\Omega_c}=f_c$. The Beltrami differential $\mu=\oplus_{j=1}^n\mu_j$ of $w$ is supported on $B=\cup_{j=1}^nB_j$. 
     For small $\epsilon>0$, take a mollifier $\eta_{\epsilon,j}\in C^\infty_c(B_j)$, and let $\mu_{\epsilon,j}=\eta_{\epsilon,j}\mu_j$. Then take solution $w_{t,\epsilon}:\CC\to \CC$ to the Beltrami equation of $\mu_{t,\epsilon}=t\oplus_{j=1}^n\mu_{\epsilon, j}+o(t)$, such that $w_{t,\epsilon}(z)=z+O(z^{-1})$ at $\infty$. Along the path $\{\Omega_{t,\epsilon}\}_{0\le t\le 1}$, we integrate \Cref{eq:qcfd} and get that for any $\epsilon$,
     \begin{equation}
         \log D_F(\partial\Omega_{1,\epsilon})-\log D_F(\Omega_c)=-\frac{1}{12}I^L(\Omega_{1,\epsilon}).
     \end{equation}
     By \Cref{lem:path}, for each disk complement $D_{j,\epsilon}$ of $\Omega_{1,\epsilon}$, there exists conformal maps $g_{j,\epsilon}\to g_j$ in $C^{5,\alpha}(\overline{\DD})$, which implies that by \Cref{lem:kernel}, the kernels $K^{(kl)}_{\epsilon}\to K^{(kl)}$ in $C^{4,\alpha}$. Therefore we have $\log D_F(\partial\Omega_{1,\epsilon})\to\log D_F(\partial\Omega)$. The convergence of the Loewner energy $I^L(\Omega_{1,\epsilon})\to I^L(\Omega)$ follows from an argument similar to \Cref{prop:leapprox}, hence we get \Cref{eq:fdle} proved.
\end{proof}
    \section{Renormalized volume of multiply-connected domains}\label{sec4}
    \subsection{The Epstein map}
    In this subsection, we present the construction and collect properties we need of the Epstein surface following \cite{Ep2024}.
    
    \textbf{Definition and explicit formulation of Epstein map:}\medskip\\
    Let $\HH^{3}$ be the upper-half space model of the $3$-dimensional hyperbolic space form, with visual boundary at infinity as $\partial_\infty\HH^{3}=\widehat{\CC}$. The isometry group $\mathrm{Isom}^+(\HH^3)\cong\mathrm{PSL
    }_2(\CC)$ induces conformal diffeomorphisms on $\partial_\infty\HH^3$ with the spherical metric. Thus for a point $x\in\HH^3$, consider any $A_x\in\mathrm{Isom}^+(\HH^3)$ which maps $x$ to $(0,0,1)$, $\nu_x=A_x^*g_{\widehat{\CC}}$ is  conformal to the spherical metric on $\widehat{\CC}$, and being invariant under the action of any isometry fixing $x$. We call $\nu_x$ the \emph{visual metric} on $\partial_\infty \HH^3$ from $x$.
    \begin{lem}\label{lem:4.1}Let $g=e^\phi |\dd z|^2$ be a conformal metric on a domain $\Omega\subset\widehat{\CC}$. At any point $z\in\Omega$:
    \begin{itemize}
        \item[(a)] The set:
        \begin{equation}
            H_z(g)=\{x\in\HH^3:\nu_x(z)=g(z)\}
        \end{equation}
        is a horosphere based at $z$.
        \item[(b)] If $z\ne\infty$, $H_z(g)$ has Euclidean radius $e^{-\phi(z)/2}$. If $z=\infty$,  $H_z(g)$ is is a euclidean hyperplane parallel to $\CC\subset\HH^3$, with euclidean height $\lim_{w\to 0}e^{\phi(-1/w)/2}/(2\vert w\vert^2)$.
        \item[(c)] If $d_z$ is the hyperbolic distance from $(0,0,1)$ to $H_z(g)$, then:
        \begin{equation}
            d_z=\lim_{w\to z}\left|\frac{\phi(w)}{2}-\log\frac{2}{1+\vert w\vert^2}\right|.
        \end{equation}
    \end{itemize}
    \end{lem}
    \begin{proof}
        For $x=(u,t)\in \HH^3=\CC\times \R^+$, we may take $A_x(w,h)=((w-u)/t,h/t)$. On $\widehat{\CC}$, we have:
        \begin{equation}
            \nu_x(z)=(A_x^*g_{\widehat{\CC}})(z)=\frac{4t^2}{(|z-u|^2+t^2)^2}|\dd z|^2.
        \end{equation}
        If $z\in \CC$, solve the equation $\nu_x(z)=g(z)$ for $x\in\HH^3$, we have that:
        \begin{equation}
            H_z(g)=\left\{(u,t)\in\HH^3:|u-z|^2+(t-e^{-\phi(z)/2})^2=e^{-\phi(z)}\right\}.
        \end{equation}
        For $z=\infty$, we have: $H_\infty(g)=\{(u,t)\in\HH^3:t=\lim_{w\to 0}\frac{e^{\phi(-1/w)/2}}{2|w|^2}\}$. To compute the hyperbolic distance to $(0,0,1)$, we first consider the case $z=\infty$, then the distance is equal to $|\log(\lim_{w\to 0}\frac{e^{\phi(-1/w)/2}}{2|w|^2})|=\lim_{w\to 0}|\phi(-1/w)/2-\log (2|w|^2)|$. Then for $z\in\CC$, consider the hyperbolic isometry:
            \begin{equation}
                I_z(w,h)=\left(\frac{w-z}{|w-z|^2+h^2},\frac{h}{|w-z|^2+h^2}\right),\quad (w,h)\in \HH^3.
            \end{equation}
        On $\widehat{\CC}$, $I_z(z)=\infty$. By explicit computation, it will map $H_z(g)$ to the hyperplane with euclidean height $\frac{1}{2}e^{\phi(z)/2}$, and maps $(0,0,1)$ to $(-1/\bar{z},1/(1+|z|^2))$, then the hyperbolic distance is:
        \begin{equation}
            d_z=\left|\log\frac{\frac{1}{2}e^{\phi(z)/2}}{\frac{1}{1+|z|^2}}\right|=\left|\frac{\phi(z)}{2}-\log\frac{2}{1+|z|^2}\right|.
        \end{equation}
    \end{proof}
    The \emph{Epstein surface} $\Sigma_g$ over $\Omega$ associated with $g$ is the envelope of the family of horospheres $\{H_z(g)\}_{z\in\Omega}$. Equivalently, it can be regarded as the image of the \emph{Epstein map} $\mathrm{Ep}_g:\Omega\to \mathbb{H}^{3}$, where $\mathrm{Ep}_g(z)$ is the tangent point of $H_z(g)$ to the envelope. By explicit computation one has:
    \begin{prop}\label{prop:4.2}(\cite{Ep2024})
        The Epstein map of a smooth conformal metric $g=e^{\phi} dz^2$ on $\Omega\subset\widehat{\CC}$ has formula:
        \begin{equation}\label{(2.1)}
            \mathrm{Ep}_g(z)=\left(z+\frac{4\n\phi}{4e^\phi+\vert\n\phi\vert^2},\frac{8e^{\phi/2}}{4e^\phi+\vert\n\phi\vert^2}\right)\in\widehat{\CC}\times\R^+=\HH^{3}.
        \end{equation}
        Here the gradient $\n\phi$ and its norm are taken in the Euclidean metric.
    \end{prop}
    \begin{lem}\label{lem:4.3}
        (Naturality and Möbius invariance). Consider $A\in\mathrm{PSL}_2(\CC)$, then by definition of $H_z(g)$, one has $A^*\nu_{A\cdot x}=\nu_x$. Furthermore, this implies that:
    \begin{equation}
        A(H_z(A^*g))=H_{A\cdot z}(g),\text{ and }A\circ\mathrm{Ep}_{A^*g}=\mathrm{Ep}_g\circ A.
    \end{equation}
    \end{lem}
    
    \textbf{Fundamental forms at infinity and geodesic flow:}\medskip\\
    Consider an immersed surface $\Sigma$ given by $i:\Omega\to\HH^{3}$. If $\vec{n}$ is its unit outward normal vector field, then the shape operator $B:T\Sigma\to T\Sigma$ is given as $Bv=-\n_v\vec{n}$, which is self-adjoint with respect to the first fundamental form $\mathrm{I}=i^*g_{\HH^{3}}$. The second fundamental form is given by $\mathrm{I\!I}(u,v)=\mathrm{I}(Bu,v)$.
    In \cite{KS08}, Krasnov and Schlenker defined the \emph{fundamental forms at infinity} as follows.
    Consider the geodesic flow:
    \begin{align}
        \Phi_t:&T^1\HH^3\to T^1\HH^3,\notag\\
        &(x,v)\mapsto\left(\exp_x(tv),\frac{\dd}{\dd t}\exp_x(tv)\right).
    \end{align} On $\Sigma$, we may define $\Phi_t:\Sigma\to\HH^{3}$ by always choosing the direction given by the outer normal vector. Then let $\Sigma_t=\Phi_t(\Sigma)$, which is the equidistant surface of distance $t$ from $\Sigma$. Let $\mathrm{I}_t$ be the first fundamental form on $\Sigma_t$, then by \cite[Lemma 2.2]{KS08}, as $t\to\infty$, $\mathrm{I}_t$ is asymptotic to $(1/4)e^{2t}\mathrm{I}^*$, where $\mathrm{I}^*=\mathrm{I}(\mathbbm{1}+B,\mathbbm{1}+B)$. Similarly, the second fundamental form and shape operator at infinity are given by:
    \begin{equation}
        \mathrm{I\!I}^*=\mathrm{I}(\mathbbm{1}+B,\mathbbm{1}-B),\quad B^*=(\mathrm{I}^*)^{-1}\mathrm{I\!I}^*=(\mathbbm{1}+B)^{-1}(\mathbbm{1}-B).
    \end{equation}
    Consider the conformal metric $g=e^\phi\vert \dd z\vert^2$ on $\Omega$. If we let $g_t=e^{2t}g$ for $t\in\R$, then by definition, at $z\in\Omega$, $H_z(g_t)$ is exactly the corresponding horosphere for $\Phi_{-t}\circ\mathrm{Ep}_{g}$. Therefore:
    \begin{equation}
        \mathrm{Ep}_{g_t}=\Phi_{-t}\circ\mathrm{Ep}_g.
    \end{equation}
    Moreover:
    \begin{prop}\label{prop:4.4}\cite[Section 5,6,7]{Ep2024}
    As $t\to+\infty$, one has:
    \begin{itemize}
        \item[(a)]For any $z\in\Omega$, $\mathrm{Ep}_{g_t}(z)\to z$.
        \item[(b)]On $\Omega$, $\mathrm{I}^*=g$.
        \item[(c)]On $\Omega$, $\mathrm{I\!I}^*=\vartheta \dd z^2+\bar{\vartheta}\dd\bar{z}^2-K^*g$, where $\vartheta=\partial_{zz}\phi-\frac{1}{2}(\partial_z\phi)^2$, and $K^*$ is the sectional curvature of $\mathrm{I}^*$.
    \end{itemize}
    \end{prop}
    \begin{rmk}
        If $g$ is taken to be the hyperbolic metric on $\Omega$, then $\vartheta=\vartheta_\Omega=S[\pi^{-1}]$, where  $\pi:\DD\to\Omega$ is a universal covering map. 
    \end{rmk}
    \textbf{Curvatures at infinity:}\medskip\\
    Let $k_-\le k_+$ be the principal curvatures of $\mathrm{I\!I}$. Then mean curvature is given as $H=(k_-+k_+)/2$. Similarly, the mean curvature at infinity $H^*$ is defined  using $\mathrm{I\!I}^*$. The Gauss and mean curvatures at infinity can be related as follows:
    \begin{prop}\cite[Section 5,6,7]{Ep2024}\label{prop4.6} Let $\dd A$ and $\dd A^*$ be the area measures for $\mathrm{I\!}$ and $\mathrm{I\!}^*$ respectively, then:
        \begin{itemize}
            \item[(a)]The principal curvatures of $\mathrm{I\!I}^*$ are $k_{\pm}^*=-K^*\pm2\vert\vartheta\vert_g$. Moreover, $H^*=-K^*$.
            \item[(b)]The Epstein map $\mathrm{Ep}_g$ is an immersion when $k^*_+\ne-1$ and $k^*_-\ne-1$.
            \item[(c)]At one point where $\Sigma_g$ is immersive, $H\dd A=\left(\frac{1-(K^*)^2}{4}+|\vartheta|_g^2\right)\dd A^*$.
        \end{itemize}
    \end{prop}
    \subsection{The Epstein surface of a multiply connected domain}
    \begin{defn}Let $\Omega\subset\widehat{\CC}$ be a multiply connected domain with hyperbolic metric $g_\Omega$. The Epstein map of $\Omega$ is the Epstein map of $g=g_\Omega$, and we denote it as $\mathrm{Ep}_\Omega$. We also use the notation $\Sigma_\Omega$ for the Epstein surface of $\mathrm{Ep}_\Omega$.
    \end{defn}
    \begin{ex}
        Let $\Omega=\DD$ or $\DD^*$, then $g_\DD=\frac{4}{(1-\vert z\vert^2)^2} |\dd z|^2$, $g_{\DD^*}=\frac{4}{(\vert w\vert^2-1)^2} |\dd w|^2$ for $z\in\DD$, $w\in\DD^*$ respectively. So we have $\phi_\DD(z)=\log2-\log(1-\vert z\vert^2)$, $\phi_{\DD^*}(w)=\log2-\log(\vert w\vert^2-1)$. Using \Cref{prop:4.2}, one can compute that:
        \begin{equation}
            \mathrm{Ep}_\DD(z)=\left(\frac{2z}{1+\vert z\vert^2},\frac{1-\vert z\vert^2}{1+\vert z\vert^2}\right),\quad\mathrm{Ep}_{\DD^*}(w)=\left(\frac{2w}{1+\vert w\vert^2},\frac{\vert w\vert^2-1}{\vert w\vert^2+1}\right).
        \end{equation}
        So $\Sigma_\DD$, $\Sigma_{\DD^*}$ are both the totally geodesic plane in $\HH^3$ bounded by $\mathbb{S}^1$ (but with opposite orientations).
    \end{ex}
    \begin{ex}\label{ex:4.2}
    Let $\A(r,1)=\{z\in\CC:r<|z|<1\}$ be a round annulus, where $0<r<1$. The hyperbolic metric on $\A(r,1)$ is:
    \begin{equation}
        g_r(z)=\left(\frac{\pi}{\log r}\right)^2\frac{1}{|z|^2\sin^2\left(\pi\frac{\log |z|}{\log r}\right)}|\dd z|^2.
    \end{equation}
    By \Cref{prop:4.2}, $\Sigma_{\A(r,1)}$ is rotationally symmetric. Thus we may express it in polar coordinates  $\mathrm{Ep}_{\A(r,1)}(z)=(\rho(z)e^{i\theta(z)},h(z))\in \HH^3=\CC\times \R_+$. Then we have:
    \begin{align}
        \rho(z)&=\vert z\vert\cdot\frac{\pi^2+\sin^2\left(\frac{\pi\log|z|}{\log r}\right)\left(\pi^2\cot^2\left(\frac{\pi\log|z|}{\log r}\right)-(\log r)^2\right)}{\pi^2+\sin^2\left(\frac{\pi\log|z|}{\log r}\right)\left(\log r+\pi\cot\left(\frac{\pi\log|z|}{\log r}\right)\right)^2}.\\
        h(z)&=|z|\cdot\frac{-2\pi\log r\sin\left(\frac{\pi\log|z|}{\log r}\right)}{\pi^2+\sin^2\left(\frac{\pi\log|z|}{\log r}\right)\left(\log r+\pi\cot\left(\frac{\pi\log|z|}{\log r}\right)\right)^2}.
    \end{align}
    In particular, the equation $\rho(z)=0$ is equivalent to:
    \begin{equation}
        \tan^2\left(\frac{\pi\log|z|}{\log r}\right)=\frac{2\pi^2}{(\log r)^2-\pi^2}.
    \end{equation}
    According to the number of solutions to this equation, we may summarize as:
    \begin{itemize}
        \item If $e^{-\pi}<r<1$, then $\rho(z)=0$ has no zero. Thus $\Sigma_{\A(r,1)}$ has no self-intersection.
        \item If $r=e^{-\pi}$, then $\rho(z)=0$ has exactly one zero $|z|=e^{-\pi/2}$. Thus $\Sigma_{\A(r,1)}$ is tangent to itself at $(0,0,e^{-\pi/2})$.
        \item If $0<r<e^{-\pi}$, then $\rho(z)=0$ has two distinct zeros. Thus $\A(r,1)$ has two intersections.
    \end{itemize}
    \end{ex}
    \begin{figure}
        \centering
        \includegraphics[width=1\textwidth]{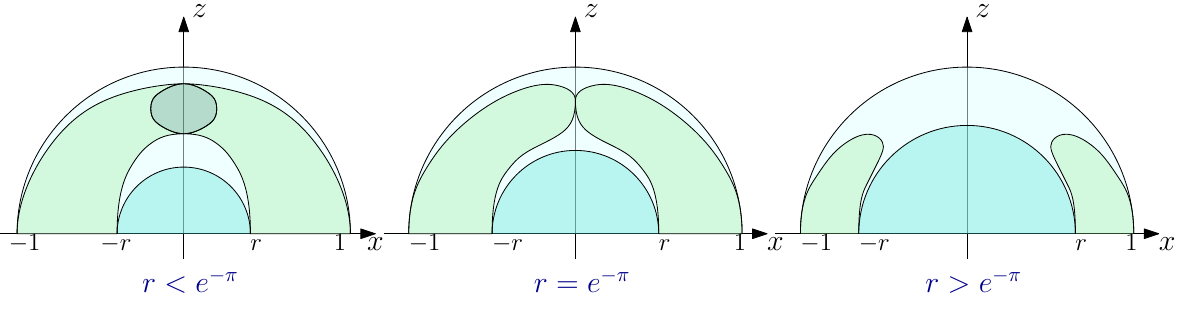}
        \caption{The cross sectional views of the Epstein surfaces $\{\Sigma_{\A(r,1)},\Sigma_{r\DD},\Sigma_{\DD^*}\}$ on the plane $\{y=0\}\subset\HH^3$ for three different values of $r$.}
        \label{fig:Multiply}
    \end{figure}
   For an $n$-connected domain $\Omega$, it has complement domains $D_1,\dots, D_n$. In \Cref{ex:4.2} we see that there can be self-intersection of $\Sigma_\Omega$. However, there will be no intersection between different surfaces in this ensemble.
    \begin{prop}\label{prop4.7}
        The surfaces $\Sigma_{D_1},\dots,\Sigma_{D_n},\Sigma_\Omega$ are pairwise disjoint.
    \end{prop}
    \begin{proof}
        Consider any two domains among $(D_1,\dots, D_n,\Omega)$, e.g., $D_j$ and $\Omega$. Then $\Omega\subset D_j^*=\widehat{\CC}\setminus \overline{D_j}$.  As the Epstein surface is the envelop of the horospheres, it's enough to show that for any $z\in D_j$, $w\in\Omega$, the horospheres $H_z^{D_j}$ and $H_w^\Omega$ are disjoint. To do so, let the hyperbolic metric on $D_i^*$ and $\Omega$ be $g_j^*=\rho_j^2 |\dd z|^2$ and $g_\Omega=\rho_\Omega^2|\dd z|^2$ respectively, then by Schwarz--Pick lemma, we have that $\rho_j(w)\le\rho_\Omega(w)$. Together with \Cref{lem:4.1}, this implies that the horosphre $H_w^{\Omega}$ lies in the interior of $H_w^{D_j^*}$. By \cite[Proposition 4.2]{BBVW25}, $H_z^{D_j}\cap H_w^{D_j^*}=\emptyset$, then it follows that $H_z^{D_j}\cap H_w^{\Omega}=\emptyset$ as well. 
    \end{proof}
    The boundary behaviour of Epstein map of multiply connected domains can be affected by the regularity of the boundary curves. We state the following result which generalizes \cite[Corollary 3.14]{BBVW25}.
    \begin{prop}\label{prop:4.8}
        If $\Omega$ is bounded by asymptotically conformal curves, then $\mathrm{Ep}_\Omega$ extends continuously to the identity on $\partial\Omega$. Moreover, it is an embedding near each boundary curve.
    \end{prop}
    \begin{proof} 
        When $z\in\Omega$ is close to one boundary curve $\gamma$, by \cite[Corollary 3.13]{BBVW25} we have that:
        \begin{equation}
            (\sqrt{5}-2)d(z,\gamma)\le d(\mathrm{Ep}_\Omega(z),\gamma)\le 5d(z,\gamma),
        \end{equation}
        where $d$ is the Euclidean distance in $\R^3$. This implies that $\mathrm{Ep}_\Omega$ extends continuously to the identity on $\gamma$. By \Cref{prop:quadraWP}, $\partial\Omega$ being asymptotically conformal implies that when $z$ is close to $\partial\Omega$, $|\vartheta|_\Omega(z)<1$, and thus by \Cref{prop4.6} $\mathrm{Ep}_\Omega$ is an immersion at $z$. It follows further from the proof of \cite[Corollary 3.14]{BBVW25} that $\mathrm{Ep}_\Omega$ is embedded near $\partial\Omega$.
    \end{proof}
    \subsection{Definition of renormalized volume}
        Following from \Cref{prop4.7}, the surfaces $\Sigma_{D_1},\dots,\Sigma_{D_n},\Sigma_\Omega$ bound a 3-dimensional subset inside $\HH^3$, with boundary at infinity as $\partial\Omega$. The \emph{signed volume} between Epstein surfaces is defined as follows. We assume that $\partial\Omega$ is asymptotically conformal. For $\epsilon>0$, let $\zeta_\epsilon\in C^\infty_0(\R^+)$, such that $\zeta_\epsilon\equiv1$ on $(2\epsilon,+\infty)$, and $\zeta_\epsilon\equiv0$ on $(0,\epsilon)$. Take the coordinate $(x,y,h)\in\HH^3$, then it has hyperbolic volume form $\mathrm{dVol}=\dd x\dd y\dd h/h^3$. Thus we may define the cut-off volume form as $\mathrm{dVol}_\epsilon=\zeta_\epsilon \mathrm{dVol}$. Consider any piecewise extension $\Phi:\overline{\HH}^3\to\overline{\HH}^3$, such that $\Phi\vert_\Omega=\mathrm{Ep}_\Omega$, and $\Phi\vert_{D_i}=\mathrm{Ep}_{D_i}$. By \Cref{prop:4.8} $\mathrm{Ep}_\Omega$ is embedding near $\partial\Omega$, so we may further assume that the Jacobian of $\Phi$ is nonnegative near $\partial\Omega$ inside $\HH^3$. Then let:
    \begin{equation}
        V(\Omega)=\lim_{\epsilon\to 0^+}\int_{\HH^3}\Phi^*\dd\mathrm{Vol}_\epsilon.
    \end{equation}
    \begin{prop}\cite[Subsection 4.2]{BBVW25}\label{prop4.9}
        The volume $V(\Omega)$ doesn't depend on the choice of extension $\Phi$, and the limit exists in $(-\infty,+\infty]$. Also, $V(\Omega)$ is M\"obius invariant.
    \end{prop}
    For $\Omega$ with $C^{5,\alpha}$ boundary, $z_0\in\partial\Omega$, we may always define the osculation circle $\mathcal{C}_\Omega(z_0)$ of $\partial\Omega$ at $z_0$ as the round circle with the same curvature as $\partial\Omega$. The osculation plane $\mathcal{P}_\Omega(z_0)$ is the geodesic plane bounded by $\mathcal{C}_\Omega(z_0)$ in $\HH^3$.
    \begin{lem}
        (Osculation of multiply connected boundary)\label{lem4.10}. Suppose that $\Omega$ has $C^{5,\alpha}$ boundary curves. Then at any $z_0\in\partial\Omega$, $\Sigma_\Omega$ and $\mathcal{P}_\Omega(z_0)$ are tangent at $(z_0,0)$ and agree up to second order.
    \end{lem}
    \begin{proof}
        By a M\"obius map we may assume that $z_0=0$, and $\mathcal{C}_\Omega(z_0)=\partial\HH^2\subset\widehat{\CC}$. Then $\mathcal{P}_\Omega(z_0)=\{y=0\}\subset\HH^3$. Let $\pi:\HH^2\to\Omega$ be a locally isometric universal covering map, such that by M\"obius map, assume that:
        \begin{equation}
            \pi(0)=0,\pi'(0)=1,\pi''(0)=0.
        \end{equation}
        Then in a small chart containing $0$, $(g_\Omega\circ\pi)\cdot|\pi'|^2=g_{\HH^2}$. In this chart, by \Cref{prop:4.2}, we may compute that:
        \begin{equation}
            (\mathrm{Ep}_\Omega\circ\pi)(z)=\left(\pi(z)+\frac{u(z)\overline{v(z)}}{1+|v(z)|^2},\frac{|u(z)|}{1+|v(z)|^2}\right),
        \end{equation}
        where for $z=x+iy$:
        \begin{equation}
            u(z)=-2iy\pi'(z)=-2iy+o(|z|^2),v(z)=1+iy\frac{\pi''(z)}{\pi'(z)}=1+o(|z|).
        \end{equation}
        Therefore, $\mathrm{Ep}_\Omega(z)=(x,0,y)+o(|z|^2)$, which coincides with $\mathcal{P}_\Omega(z_0)$ up to second order.
    \end{proof}
    By \Cref{prop:4.8} and \Cref{prop:kellogg}, we may find $\epsilon>0$, such that for every $U\in\{\Omega,D_1,\dots,D_n\}$ below the hyperplane of euclidean height $h=\epsilon$, $\Sigma_U$ is embedded. Then for every $1\le j\le n$, $\Sigma_\Omega,\Sigma_{D_j}$ and $\{h=\epsilon\}$ bounded a small oriented region in $\HH^3$. We may parameterize it using the following coordinate defined on a collar of $\gamma_j$ times $(0,\epsilon]$:
    \begin{equation}\label{eq:91}
        G_j(s,a,h)=\left(\gamma_j(s)+a\vec{n}_j(s),h\right)\in\widehat{\CC}\times\R_+=\HH^3,
    \end{equation}
    where $\vec{n}_j$ is the outward normal vector of $D_j$.
    Using this coordinate, we may prove the following finiteness result analogous to \cite[Proposition 4.6]{BBVW25}. 
    \begin{prop}
        Let $\Omega$ be a multiply connected domain  with $C^{5,\alpha}$ boundary. Then the volume $V(\Omega)$ is finite.
    \end{prop}
    \begin{proof}
        By \Cref{prop:kellogg}, we know that there is a conformal map $f\in C^{5,\alpha}(\overline{\Omega_c})$. Then by \Cref{prop:4.2} we have that $\mathrm{Ep}_U$ extends to be $C^{3,\alpha}$ on $\overline{U}$, where $U\in\{\Omega,D_1,\dots,D_n\}$. For every $1\le j\le n$ and $0<h<\epsilon$, in the coordinate $G_j$, there exists functions $a^\pm(s,h)$, such that $G_j(s,a^+(s,h),h)\in\Sigma_\Omega$ and $G_j(s,a^-(s,h),h)\in\Sigma_{D_j}$. By inverse function theorem, $a^\pm(s,h)$ are $C^{3,\alpha}$. Thus by \Cref{lem4.10}, we have $|a^+(s,h)-a^-(s,h)|=O(h^3).$ Using this bound, and the computation in \cite[Proposition 4.6]{BBVW25}, we have that $V(\Omega)<\infty$.
    \end{proof}
    \begin{defn}\label{def4.2}
         Let $\Omega\subset\widehat{\CC}$ be a multiply connected domain with $C^{5,\alpha}$ boundary, we define its renormalized volume as:
        \begin{equation}
            V_R(\Omega)=V(\Omega)-\frac{1}{2}\int_{(\cup_{i=1}^n\Sigma_i)\cup\Sigma_\Omega}H\dd A.
        \end{equation}
    \end{defn}
    \begin{prop}
        $V_R$ can also be expressed using the data of fundamental forms at infinity:
        \begin{equation}
            V_R(\Omega)=V(\Omega)-\frac{1}{2}\int_{\Omega}\frac{|\vartheta_\Omega|^2}{\rho^2_\Omega}\dd^2z-\frac{1}{2}\sum_{j=1}^n\int_{D_j}\frac{|S[g_j^{-1}]|^2}{\rho^2_{D_j}}\dd^2z.
        \end{equation}
        In particular, $V_R$ is M\"obius invariant.
    \end{prop}
    \begin{proof}
        By \Cref{prop:4.4}, $\mathrm{Ep}_\Omega$ is immersive almost everywhere on $\Omega$, and at those points, $H\dd A=(\frac{1-(K^*)^2}{4}+|\vartheta|^2_g)\dd A^*$. Since $\mathrm{\!I}^*=g_\Omega$ is hyperbolic, we have $K^*=-1$. Thus $H\dd A=|\vartheta|_g^2\dd A^*$ is M\"obius invariant. From \Cref{lem:4.3} we also have that $V(\Omega)$ is M\"obius invariant.
    \end{proof}
    \section{Quasiconformal variation of the renormalized volume}\label{sec5}
    Throughout this section, we assume that $\Omega$ is $n$-connected with $C^{5,\alpha}$ boundary, such that $\infty\in\Omega$. We use the variation formula to prove that under this regularity, the renormalized volume of $\Omega$ is equal to the Loewner energy up to two constants.
    \subsection{The first variation formula}
    We consider the following Beltrami differentials:
    \begin{equation}\label{eq99}
        \mathrm{Bel}^{4,\alpha}(\Omega):=\{\mu\in L^\infty_1(\CC):\mu\text{ extends to be }C^{4,\alpha} \text{ on }\overline{D_j} \text{ and } \overline{\Omega}, \forall 1\le j\le n\}.
    \end{equation}
    We emphasize that we do not require the boundary value of the extension of $\mu|_{D_j}$ and $\mu|_{\Omega}$ to be the same on $\gamma_j$ for each $j$. Now consider a path $t\mapsto \mu_t$ which lies in $C^1(I,\mathrm{Bel}^{4,\alpha}(\Omega))$, where $I$ is an interval containing $0$. For this path,  take the normalized solution $w_t:\CC\to \CC$ to $\partial_{\bar{z}}w_t=\mu_t\cdot\partial_z w_t$ such that $w_t$ fixes three points $0,1,\infty$. We prove regularity results of $(t,z)\mapsto w_t(z)$. To do this, we need the elliptic regularity of Beltrami equations.
    \begin{prop}\label{lem5.1}
        For every $0<\beta<\alpha$, the map $t\mapsto w_t$ is in $C^1(I,C^{5,\beta}(\overline{U}))$ for every $U\in\{D_1,\dots, D_n,\Omega\}$. Moreover, $\gamma_{t,j}=w_t(\gamma_j)$ is of $C^{5,\beta}$ regularity. 
    \end{prop}
    We postpone the proof of this proposition to \Cref{appendixB}.
    An explicit example of such a construction can be given as follows:
    \begin{lem}\label{lem:path}
        There is a $C^{5,\alpha}$ family of domains $\{\Omega_t\}_{t\in[0,1]}$, such that each $\Omega_t$ is conformally equivalent to $\Omega_c$, $\Omega_0=\Omega_c$, and $\Omega_1=\Omega$. The family of domains is parameterized by a family of quasiconformal maps $w_t$, i.e., $\partial\Omega_t=w_t(\partial\Omega_c)$, such that the corresponding Beltrami differentials $t\mapsto\mu_t\in \mathrm{Bel}^{4,\alpha}(\Omega)$.
    \end{lem}
\begin{proof}
    Take a conformal map $f_c:\Omega_c\to\Omega$ with $f_c(\infty)=\infty$. By \Cref{prop:kellogg}, for each $1\le j\le n$, $f_c$ extends to a $C^{5,\alpha}$ diffeomorphism $\varphi_j=f|_{C_j}:C_j\to\gamma_j$. We construct a quasiconformal map from $\CC$ to $\CC$ which equals to $f_c$ when restricted to $\Omega_c$. Indeed, take conformal maps $g_j:B_j\to D_j$, and consider the welding map $\phi_j:g_j^{-1}\circ f_c|_{C_j}:C_j\to C_j$ for each $j$. Then $\phi_j(z_j+r_je^{i\theta})=z_j+r_je^{i\beta_j(\theta)}$ is an orientation-preserving  $C^{5,\alpha}$ circle diffeomorphism, which implies that $\beta_j'>0$, and it is bounded on the unit circle. Take a smooth cut-off function $\chi:[0,1]\to [0,1]$, such that $\chi=0$ near $0$ and $\chi=1$ near $1$. Let $\hat{w}_j:B_j\to B_j$ defined as:
    \begin{equation}
        \hat{w}_j(z_j+\lambda r_je^{i\theta})=z_j+\lambda r_j\exp(i(\theta+\chi(\lambda)(\beta_j(\theta)-\theta))),\quad 0\le \lambda\le 1.
    \end{equation} Then $\hat{w}_j$ has complex dilatation  $\mu_j(z_j+\lambda r_je^{i\theta})=e^{2i\theta}\frac{\chi(\lambda)(1-\beta_j'(\theta))}{2-\chi(\lambda)(1-\beta_j'(\theta))}$, which implies that $\hat{w}_j$ is quasiconformal. In this way, the maps $\{g_j\circ \hat{w}_j\}_{1\le j\le n}$ and $f_c$ glues into a quasiconformal map $w_1:\CC\to\CC$, with Beltrami differential $\mu=\sum_{j=1}^n\mu_j\cdot \mathbbm{1}_{B_j}$. For $0\le t\le 1$, consider the unique solution $w_t:\CC\to\CC$ of $\partial_{\bar{z}}w_t=t\mu\cdot\partial_zw_t$ such that $w_t$ fixes three points $0,1,\infty$. 
    By the explicit expression of $\mu$, we have that $\mu$ is $C^{4,\alpha}$ on $\cup_j\overline{B_j}$. 
\end{proof}
    To proceed with the variation of renormalized volume, we construct a cut-off construction similar to  \cite[Section 5.2]{BBVW25}. Take $t\mapsto \mu_t\in\mathrm{Bel}^{4,\alpha}(\Omega)$, and $t\mapsto w_t$ is the corresponding family of normalized quasiconformal maps. Let $\Omega_t=w_t(\Omega)$, $D_{t,j}=w_t(D_j)$, and $\gamma_{t,j}=w_t(\gamma_j)$. 
    Take $\epsilon>0$ small, such that for each $1\le j\le n$, under the level set $\{h=\epsilon\}$, the Epstein surfaces $\Sigma_{D_{t,j}}$ and $\Sigma_{\Omega_t}$ are embedded. Then near each $\gamma_{t,j}$, there is a collar neighbourhood, on which there is a coordinate defined as \Cref{eq:91}:
    \begin{equation}
        G_{t,j}(s,a,h)=(\gamma_{t,j}(s)+a\vec{n}_{t,j}(s),h),\quad 0<h\le \epsilon,
    \end{equation}
    where $\vec{n}_{t,j}$ is the outward unit normal vector of $\gamma_{t,j}$.
    Using this coordinate, we may parameterize the space bounded by $\Sigma_{\Omega_t}$, $\Sigma_{D_{t,j}}$, and the horizontal piece $\{h=\epsilon\}$. Then there is a small $r>0$, such that for every $1\le j\le n$ and $t\in I$, $G_{t,j}$ is defined on $\mathbb{S}^1_s\times(-r,r)_a\times(0,\epsilon]_h$. Let $\Pi_{s,h}$ be the projection of $G$-coordinate onto $(s,h)$ components. For each $j$, define: \begin{equation}
        P_{t,j}^\pm=(\Pi_{s,h}\circ G_{t,j}\circ \mathrm{Ep}_{U_{t,j}^\pm})^{-1}:\mathbb{S}^1_s\times(0,\epsilon]_h\to U_{t,j}^\pm,\quad U^+_{t,j}=\Omega_t,U_{t,j}^-=D_{t,j}.
    \end{equation}
    This means that for given $(s,h)$, we find the point on the Epstein surface $\Sigma_{U_{t,  j}^\pm}$ with corresponding $(s,h)$ in $G$-coordinate, then map it back to $U_{t,j}^\pm$ by the inverse of Epstein map.
    Let $\Omega(\epsilon)=\Omega\setminus P_{0,j}^+(\mathbb{S}^1_s\times(0,\epsilon])$, $D_j(\epsilon)=D_j\setminus P_{0,j}^-(\mathbb{S}^1_s\times(0,\epsilon]_h)$.
    Then we may define the following piecewisely defined homeomorphism which maps from $\widehat{\CC}\to\HH^3$:
    \begin{equation}
        E_{t,\epsilon}(z)=\begin{cases}
            
            (\mathrm{Ep}_{\Omega_t}\circ w_t)(z),\quad&z\in \Omega(\epsilon),\\
            (\mathrm{Ep}_{D_{t,j}}\circ w_t)(z),&z\in D_j(\epsilon),\\
            \mathrm{Ep}_{D_{t,j}}\circ P_{t,j}^-\circ (P_{0,j}^-)^{-1}(z), &z\in D_{j}\setminus D_j(\epsilon),\\
            \mathrm{Ep}_{\Omega_t}\circ P_{t,j}^+\circ (P_{0,j}^+)^{-1}(z), &z\in \Omega\setminus \Omega(\epsilon).\\
        \end{cases}
    \end{equation}
    \begin{prop}\label{prop5.3}
        For every $0<\beta<\alpha$ and every fixed $\epsilon>0$ such that the Epstein surfaces under the height $\{h=\epsilon\}$ are embedded, the piecewise maps $t\mapsto E_{t,\epsilon}\in C^1(I,C^{3,\beta})$ on its own domain respectively. 
    \end{prop}
    \begin{proof}
        We first consider the pieces on $\Omega(\epsilon)$ and $D_j(\epsilon)$. Consider $U\in \{D_{1},\dots, D_{n},\Omega\}$ with universal covering map $\pi:\DD\to U$. Let $\nu_t=\pi^*\mu_t\in L^\infty_1(\DD)$, and take the reflected solution $w_{\nu_t}:\DD\to\DD$ of Beltrami equation (see notations in \Cref{def2.1}). By \Cref{lem5.1} again, $t\mapsto w_{\nu_t}\in C^1(I,C^{5,\beta}(\overline{\DD}))$. Now we have that $\pi_t=w_t\circ \pi\circ w^{-1}_{\nu_t}:\DD\to U_t$ is a universal covering map, where $U_t=w_t(U)$. Let $g_t=\rho^2_tdz^2$ be the hyperbolic metric of $U_t$, then we may compute explicitly that:
        \begin{equation}
            \rho_t\circ w_t\circ \pi=\rho_t\circ \pi_t\circ w_{\nu_t}=\frac{2}{1-|w_{\nu_t}|^2}\frac{1}{|\pi'|}\sqrt{\frac{J_{w_{\nu_t}}\circ \pi}{ J_{w_t}}}\in C^1(I,C^{4,\beta}(\DD)),
        \end{equation}
        where $J_\cdot$ is the Jacobian. Therefore, $(\n\log\rho_t)\circ w_t=\n(\log\rho_t\circ w_t)(dw_t)^{-1}$ is $C^1(I,C^{3,\beta})$ on a chart of $\pi$. For $D_j(\epsilon)$, it can be covered by finitely many charts, so the regularity can be upgraded to $C^1(I,C^{3,\beta}(\overline{D_j(\epsilon)}))$.
        Thus by \Cref{prop:4.2}, we have:
        \begin{equation}
            E_{t,\epsilon}=\left(w_t+\frac{(\n\log\rho_t)\circ w_t}{2\rho_t^2\circ w_t+2|(\n\log\rho_t)\circ w_t|^2},\frac{2\rho_t\circ w_t}{2\rho_t^2\circ w_t+2|(\n\log\rho_t)\circ w_t|^2}\right),
        \end{equation}
        which implies that $t\mapsto E_{\epsilon,t}$ is in $C^1(I,C^{3,\beta}))$ on $D_j(\epsilon)$. For $\Omega(\epsilon)$, there is an unbounded subdomain containing $\infty$. However, we may take an inversion $z\mapsto -1/z$ to map it into a bounded one. By the naturality of Epstein map (\Cref{lem:4.3}), the finite covering argument still works. Finally, we consider $E_{t,\epsilon}$ when defined on $\Omega\setminus D_j(\epsilon)$ or $\Omega\setminus\Omega(\epsilon)$. Consider the $G_{t,j}$-coordinate near $\gamma_{t,j}$, for given $(s,h)$, by embeddedness, there are unique $a^\pm_t(s,h)$, such that $G_{t,j}(s,a^\pm_t(s,h),h)\in \Sigma_{U_{t,j}^\pm}$, where $U^+_{t,j}=\Omega_t$, $U^-_{t,j}=D_{t,j}$. Thus by inverse function theorem, $t\mapsto a^\pm_t$ is in $C^1(I,C^{3,\beta})$. Note that the image of $E_{t,\epsilon}$ of those parts can be parameterized by the $G_{t,j}$-coordinate as:
        \begin{equation}
            (\gamma_{t,j}(s)+(\lambda a^+_t(s,\epsilon)+(1-\lambda)a^-_t(s,\epsilon))\vec{n}_{t,j},\epsilon),\quad s\in \mathbb{S}^1,0\le \lambda\le 1.
        \end{equation}
        Therefore, by the regularity of $t\mapsto \gamma_{t,j}$ and $t\mapsto a_t^\pm$, we have that $t\mapsto E_{t,\epsilon}$ is also $C^1(I,C^{3,\beta})$ on $\Omega\setminus\Omega(\epsilon)$ or $D_j\setminus D_j(\epsilon)$.
    \end{proof}
    
\begin{thm}\label{thm:varre}
        Let $t\mapsto \mu_t$ be a path in $C^1(I,\mathrm{Bel}^{4,\alpha}(\Omega))$ for  $\alpha>0$, and $w_t$ is the corresponding normalized solution to Beltrami equation. If $\mu_t=t(\oplus_{j=1}^n\dot{\mu}_j+\dot{\mu}_\Omega)+o(t)$, which means that $\dot{\mu}_j,\dot{\mu}_\Omega$ are supported on $D_j,\Omega$ respectively. Then:
        \begin{equation}\label{eq5.4}
            \delta V_R(\Omega)=\left.\frac{\dd}{\dd t}\right|_{t=0}V_R(\Omega_t)=-\mathrm{Re}\left\{\sum_{j=1}^n\int_{D_j}\dot{\mu}_j\cdot S[g_j^{-1}] \dd ^2z+\int_{\Omega}\dot{\mu}_\Omega\cdot\vartheta_{\Omega}\dd^2 z\right\},
        \end{equation}
        where we recall that  $g_j:B_j\to D_j$ is  conformal, and $\vartheta_\Omega=S[\pi^{-1}]$, with $\pi:\HH\to\Omega$ being the universal covering map.
    \end{thm}
    \begin{proof}
    By \Cref{prop5.3}, \cite[Theorem 4]{Sou04}, and \cite[Theorem 5.7]{BBVW25}, we have that:
    \begin{align}
        \delta V(\Omega)=\sum_{j=1}^n\int_{D_j}\mathrm{Ep}_{D_j}^*\left(\delta H+\frac{1}{4}\langle\delta\mathrm{I},\mathrm{I\!I}\rangle \dd A\right)+\int_{\Omega}\mathrm{Ep}_{\Omega}^*\left(\delta H+\frac{1}{4}\langle\delta\mathrm{I},\mathrm{I\!I}\rangle \dd A\right),
    \end{align}
    where $\langle P,Q\rangle=\mathrm{tr}(\mathrm{I}^{-1}P\mathrm{I}^{-1}Q)$. Then by \Cref{def4.2}, \Cref{prop:4.2}, and \cite[Corollary 6.2]{KS08}, we have:
    \begin{align}
        \delta V_R(\Omega)&=\delta V(\Omega)-\frac{1}{2}\delta\left(\int_{\Sigma_\Omega\cup(\cup_{j=1}^n \Sigma_{D_j})}H\dd A\right)\notag\\
        &=-\frac{1}{8}\int_{\Omega\cup(\cup_{j=1}^nD_j)}\langle \delta\mathrm{I}^*,\mathrm{I\!I}^*_0\rangle \dd A^*,
    \end{align}
    where $\mathrm{I\!I}^*_0$ is the traceless part of $\mathrm{I\!I}^*$. For a symmetric $(0,2)$-tensor $P$, we may express it in matrix form by writing it as $P=p_{11}\dd z^2+2p_{12}\dd z\dd\bar{z}+p_{22}\dd\bar{z}^2$.
    Then we have \Cref{prop:4.2}:
    \begin{align}
        \mathrm{I^*}^{-1}=\begin{pmatrix}
            0 & 2e^{-2\rho}\\
            2e^{-2\rho} & 0
        \end{pmatrix},\quad \mathrm{I\!I}^*_0=\begin{pmatrix}
            \vartheta & 0\\
            0 & \overline{\vartheta}
        \end{pmatrix}.
    \end{align}
    For $\delta\mathrm{I}^*=\left.\frac{d}{dt}\right|_{t=0}w_t^*\mathrm{I}^*$, by explicit computation we have:
    \begin{align}
        w_t^*\mathrm{I}^*&=e^{2\rho\circ w_t}\dd w_t\dd\overline{w}_t\notag\\
        &=e^{2\rho}(1+2t\rho'\delta w)(\dd z+t(\partial_z\delta w \dd z+\dot{\mu}
        \dd\bar{z}))(\dd\bar{z}+t(\overline{\partial_z\delta w}\dd z+\overline{\dot{\mu}}\dd z))+o(t)\notag\\
        &=\mathrm{I}^*+te^{2\rho}\begin{pmatrix}
            \overline{\dot{\mu}} & \mathrm{Re}(\partial_z\delta w)+\rho'\delta w\\
            \mathrm{Re}(\partial_z\delta w)+\rho'\delta w &\dot{\mu}
        \end{pmatrix}+o(t).
    \end{align}
    Therefore,
    \begin{align}
        \delta V_R(\Omega)&=-\frac{1}{8}\int_{\Omega\cup(\cup_{j=1}^n D_j)}\mathrm{tr}(\mathrm{I^{-1}}^*\delta\mathrm{I}^*\mathrm{I^*}^{-1}\mathrm{I\!I}^*_0)e^{2\rho}
        \dd^2z\notag\\
        &=-\mathrm{Re}\int_{\Omega\cup(\cup_{j=1}^n D_j)}\dot{\mu}\cdot\vartheta \dd^2z.
    \end{align}
    \Cref{eq5.4} follows from the fact that on $D_j$, $\vartheta_{D_j}=S[g_j^{-1}]$.
    \end{proof}
    \begin{cor}
        Let $\Omega$ be a multiply connected domain in $\widehat{\CC}$ with $C^{5,\alpha}$ boundary, such that it is a critical point of the variation of the renormalized volume, then it is a circle domain.
    \end{cor}
    \begin{proof}
        By \Cref{thm:varre}, we have that for each $1\le j\le n$, the Schwarzian derivative $S[g_j^{-1}]=0$. This implies that every $g_j^{-1}$ is a Möbius map. Therefore, the boundary curves $\gamma_j=g^{-1}_j(\mathbb{S}^1)$ is a round circle in $\widehat{\CC}$, which means that $\Omega$ is a circle domain. 
    \end{proof}
    \begin{cor}\label{cor5.7} When $\Omega$ has $C^{5,\alpha}$ boundary, we have:
    \begin{align}
        V_R(\Omega)-V_R(\Omega_c)&=\frac{\pi}{4}I^L(\Omega).
    \end{align}
\end{cor}
\begin{proof}
   Consider the family $\{\Omega_t\}_{0\le t\le 1}$ defined in \Cref{lem:path}, with $\Omega_0=\Omega_c$ and $\Omega_1=\Omega$. By  \Cref{thm:wpvarle}, \Cref{thm:varre}, $V_R(\Omega)-V_R(\Omega_c)-\frac{\pi}{4}I^L(\Omega)$ has variation $0$, which implies that $t\mapsto V_R(\Omega_t)-V_R(\Omega_c)-\frac{\pi}{4}I^L(\Omega_t)$ has derivative $0$ on $(0,1)$. This implies that:
   \begin{equation}
       V_R(\Omega_t)-V_R(\Omega_c)-\frac{\pi}{4}I^L(\Omega_t)\equiv V_R(\Omega_c)-V_R(\Omega_c)-\frac{\pi}{4}I^L(\Omega_c)=0.
   \end{equation}
\end{proof}
\begin{prop}(The modulus-preserving second variation).
    Let $w_t$ a quasiconformal variation driven by harmonic Beltrami differentials $\mu_t=t\dot{\mu}=t\oplus_{j=1}^n\dot{\mu}_j$ supported on $D$.
    Then we have:
    \begin{align}\label{eq5.16}
        \delta^2 V_R(\Omega)&=\left.\frac{\dd^2}{\dd t^2}\right|_{t=0}V_R(\Omega_t)\notag\\
        &=-\frac{6}{\pi}\sum_{i,j=1}^n\mathrm{Re}\left\{\int_{D_i\times D_j}\frac{\dot{\mu}_i(z)\dot{\mu}_j(w)}{(w-z)^4}\dd ^2w\dd^2z\right\}+\frac{1}{2}\sum_{j=1}^n\Vert\dot{\mu}_j\Vert^2_{\mathrm{WP}}.
    \end{align}
\end{prop}
\begin{proof}
    We take a complex family of Beltrami differentials, namely, $\mu_\epsilon=\epsilon\cdot\dot{\mu}$, $\epsilon$ lies in a small ball containing $0$ in $\CC$. Let $h_j=g_j^{-1}:D_j\to \DD$ conformally. By \Cref{eq5.4}, the complex first variation formula is:
    \begin{equation}
        \left.\frac{\partial}{\partial\epsilon}\right|_{\epsilon=0}V_R(w_\epsilon(\Omega))=-\frac{1}{2}\sum_{j=1}^n\int_{D_j}\dot{\mu}_j\cdot S[h_j]\dd^2z.
    \end{equation}
    Let $D_{\epsilon,j}=w_\epsilon(D_j)$, and $h_{\epsilon,j}:D_{\epsilon,j}\to \DD$ is conformal. By the computation in \cite[Theorem 5.5]{Ah06}, 
    we have:
    \begin{align}
        \frac{\partial}{\partial\epsilon}V_R(w_\epsilon(\Omega))
        &=-\frac{1}{2}\sum_{j=1}^n\int_{D_j}\dot{\mu}_j\cdot(S[h_{\epsilon,j}]\circ w_\epsilon)\cdot(\partial_zw_{\epsilon})^2 \dd^2z.
    \end{align}
    Let $v_{\epsilon,j}=h_{\epsilon,j}\circ w_\epsilon\circ g_j:\DD\to\DD$, it is a quasiconformal map with Beltrami differential $\epsilon g_j^*\dot{\mu}_j$.
    Formal computation of Schwarzian gives:
    \begin{equation}
        (S[h_{\epsilon,j}]\circ w_\epsilon)\cdot(\partial_zw_\epsilon)^2=(S[v_{\epsilon,j}]\circ h_j)\cdot(h_j')^2+S[h_j]-S[w_\epsilon].
    \end{equation}
    Since $w_\epsilon$ depends holomorphically on $\epsilon$, by \cite[Theorem 6.12.5]{Hub06}, we have:
    \begin{align}
        \left.\frac{\partial}{\partial\epsilon}\right|_{\epsilon=0}S[w_{\epsilon}](z)=-\frac{6}{\pi}\sum_{j=1}^n\int_{D_j}\frac{\dot{\mu_j}(w)}{(w-z)^4}\dd^2w,\quad \left.\frac{\partial}{\partial\bar\epsilon}\right|_{\epsilon=0}S[w_{\epsilon}](z)=0.
    \end{align}
    To compute the derivative of $S[v_{\epsilon,j}]$, take reflection $\widehat{\epsilon(g_j^*\dot{\mu}_j)}(z)=\overline{\epsilon(g_j^*\dot{\mu}_j)(1/\bar{z})}(z/\bar{z})^2$ onto $\DD^*$, and extend $v_{\epsilon,j}$ to be on $\CC$. Then $v_\epsilon$ depends holomorphically on $\epsilon$ over $\DD$, and anti-holomorphically over $\DD^*$. By \cite[Lemma 1.1.12]{TT06} and computation as above, we have:
    \begin{align}
    \left.\frac{\partial}{\partial\epsilon}\right|_{\epsilon=0}S[v_{\epsilon,j}](z)=0,\quad \left.\frac{\partial}{\partial\bar\epsilon}\right|_{\epsilon=0}S[v_{\epsilon,j}](z)=-\frac{6}{\pi}\int_{\DD}\frac{\overline{(g_j^*\dot{\mu}_j)(w)}}{(1-z\overline{w})^4}\dd^2w.
    \end{align}
    By \cite[Lemma 5.10]{BBVW25} and  \cite[Proposition 5.4.9]{Hub06}, we have:
    \begin{align}
        \left.\frac{\partial^2}{\partial\epsilon^2}\right|_{\epsilon=0}V_R(w_{\epsilon}(\Omega))&=-\frac{3}{\pi}\sum_{i,j=1}^n\int_{D_i}\dot{\mu}_i(z)\left(\int_{D_j}\frac{\dot{\mu}_j(w)}{(w-z)^4}\dd^2w\right)\dd^2z,\notag\\
        \left.\frac{\partial^2}{\partial\bar{\epsilon}\partial\epsilon}\right|_{\epsilon=0}V_R(w_{\epsilon}(\Omega))&=\frac{1}{4}\sum_{j=1}^n\int_{D_j}\dot{\mu}_j\cdot(\overline{(g_j^*\dot{\mu}_j)}\rho^2_{\DD})\circ h_j(h_j')^2\dd^2z\notag\\
        &=\frac{1}{4}\sum_{j=1}^n\int_{D_j}|\dot{\mu}_j|^2(z)\rho^2_{D_j}(z)\dd^2z=\frac{1}{4}\sum_{j=1}^n\Vert\dot{\mu}_j\Vert^2_{\mathrm{WP}}.
    \end{align}
    The formula we want follows from the fact that:
    \begin{equation}
        \left.\frac{\dd^2}{\dd t^2}\right|_{t=0}V_R(w_t(\Omega))=2\mathrm{Re}\left(\left.\frac{\partial^2}{\partial\epsilon^2}\right|_{\epsilon=0}V_R+\left.\frac{\partial^2}{\partial\bar{\epsilon}\partial\epsilon}\right|_{\epsilon=0}V_R\right).
    \end{equation}
\end{proof}
\begin{cor}
        At the critical point, namely, $\Omega=\Omega_c$ is a circle domain, the second variation formula is positively definite. Moreover, there is a constant $0<b<1$, such that:
    \begin{equation}
        \frac{1-b}{2}\sum_{j=1}^n\Vert\dot{\mu}_j\Vert^2_{\mathrm{WP}}\le\delta^2V_R(\Omega)\le\frac{1+b}{2}\sum_{j=1}^n\Vert\dot{\mu}\Vert^2_{\mathrm{WP}}.
    \end{equation}
    \end{cor}
    \begin{proof}
        At $\Omega=\Omega_c$, every $D_j$ is a round disk, then by \cite[Lemma 1.1.12]{TT06}, the diagonal integrals in \Cref{eq5.16} vanish. Thus:
        \begin{equation}\label{eq6.25}
            \delta^2V_R(\Omega)=\frac{1}{2}\sum_{j=1}^n\Vert \dot{\mu}_j\Vert^2_{\mathrm{WP}}-\frac{6}{\pi}\sum_{1\le i\ne j\le n}\mathrm{Re}\left\{\int_{D_i\times D_j}\frac{\dot{\mu}_i(z)\dot{\mu}_j(w)}{(w-z)^4}\dd ^2w\dd^2z\right\}.
        \end{equation}
        We use the following $A_2$ version block Grunsky operator similar to what we have in \Cref{eq:4.5}. Let $H'_D=\oplus_{j=1}^nA_2(D_j)$, then:
        \begin{align}
            B_{A_2}:\overline{H'_D}&\to H_D',\notag\\
                    \overline{u}&\mapsto\sum_{j=1}^n\int_{D_j}\frac{12}{\pi}\frac{\overline{u_j(w)}}{(z-w)^4}\rho_{D_j}^{-2}(w)\dd^2w.
        \end{align}
        By classical theory of Bergman kernels \cite{Aro50} and \Cref{rmk4.3} , we have that $\Vert B_{A_2}\Vert_{\mathrm{op}}\le\Vert B_{A_2^1}\Vert^2_{\mathrm{op}}<1$, where we use $B_{A_2^1}$ to denote the block Grunsky operator of $\Omega$ we defined in \Cref{sec6.1}. By assumption that $\dot{\mu}$ is harmonic, we may take $\dot{\mu}_j=\overline{\varphi_j}\rho_{D_j}^{-2}$, where $\varphi_j\in A_2(D_j)$, and $\Vert \dot{\mu}_j\Vert_{\mathrm{WP}}=\Vert\varphi_j\Vert_{A_2}$. By \Cref{eq6.25} we have:
        \begin{align}
            |\delta^2V_R(\Omega)|&=\frac{1}{2}\left|\sum_{j=1}^n\Vert \dot{\mu}_j\Vert^2_{\mathrm{WP}}-\mathrm{Re}\langle B_{A_2}\overline{\varphi},\varphi\rangle_{H_D'}\right|\notag\\
            &\le \frac{1}{2}(1+\Vert B_{A_2}\Vert_{\mathrm{op}})\sum_{j=1}^n\Vert\dot{\mu}_j\Vert^2_{\mathrm{WP}},
        \end{align}
        where we use the notation $\varphi=(\varphi_1,\dots,\varphi_n)\in H_D'$.
    \end{proof}
\subsection{Comparison between domain functionals on round annuli}
\textbf{Renormalized volume of annuli:}
\begin{prop}\label{prop5.9}(Variation of annuli). For the round annulus $\A(r,1)=\{z\in\CC:r<|z|<1\}$, where $0<r<1$, we have:
    \begin{equation}
         \frac{\dd}{\dd r}V_R(\A(r,1))=\frac{\pi}{2r}\left(1+\frac{\pi^2}{(\log r)^2}\right).
    \end{equation} 
    Thus by integrating back the formula above, we have:
    \begin{equation}
        V_R(\A(r,1))=\frac{\pi}{2}\left(\log r-\frac{\pi^2}{\log r}\right).
    \end{equation}
    \end{prop}
    \begin{proof}
        Consider $0<r<1$, and a small $\vert t\vert<\min(r,1-r)$. On $\A(r,1)$, in polar coordinate $z=\rho e^{i\theta}$, let $F_t(z)=\exp(\frac{\log(r+t)}{\log r}\log \rho+i\theta)$. Then $F_t$ maps $\A(r,1)$ to $\A(r+t,1)$ quasiconformally, with Beltrami differential $\mu_t=t\dot{\nu}+o(1)$, where $\dot{\nu}=\frac{1}{2r\log r}\frac{z}{\overline{z}}$. We may extend $F_t$ to $\widehat{\CC}$ by letting $F_t(z)=\log((r+t)z)$ on $r\DD$, and $F_t(z)=z$ on $\DD^*$. By standard gluing argument we have that $F_t$ is quasiconformal on $\widehat{\CC}$, with Beltrami differential $\mu_t\cdot\mathbbm{1}_{\A(r,1)}$. Using the quasiconformal variation given by $F_t$, by \Cref{thm:varre} we have:
        \begin{align}
             \frac{\dd}{\dd r}V_R(\A(r,1))=\left.\frac{\dd}{\dd t}\right|_{t=0}V_R(\A(r+t,1))&=-\mathrm{Re}\int_{\A_r}\dot{\nu}\cdot\vartheta_r \dd^2z\notag\\
             &=-\int_{\A_r}\frac{1}{4r\log r}\frac{1}{\vert z\vert^2}\left(1+\frac{\pi^2}{(\log r)^2}\right) \dd^2z\notag\\
             &=\frac{\pi}{2r}\left(1+\frac{\pi^2}{(\log r)^2}\right).
        \end{align}
        Integrate this derivative formula gives the formula of $V_R(\A(r,1))$ with a constant $C$. By explicit computation at $r=e^{-\pi}$ (see \Cref{appendixA} for details), we get that $V_R(\A(e^{-\pi,1}))=0$. Then we get that $C=0$.
        \end{proof}
        \textbf{Loewner potential of round annuli:}\medskip\\
        Since $\A(r,1)$ is a circle domain, we consider the Loewner potential $\mathcal{H}_2(\A(r,1))$ alternatively. In \cite[Proposition 4.2]{LM25}, we computed the following formula:
        \begin{equation}
            \mathcal{H}_2(\A(r,1))=-\log\log\frac{1}{r}-2\sum_{k=1}^\infty\log(1-r^{2k})+C.
        \end{equation}
        Where $C$ is a constant independent of $r$.\medskip\\
        \textbf{Grunsky determinant of round annuli}:\medskip\\
        For $\A(r,1)$, the complement domains are $D_1=r\DD$, $D_2=\DD^*$. So the conformal maps are $g_1(z)=rz$, $g_2(z)=1/z$. Since they are M\"obius, by \Cref{lem:kernel}, the diagonal blocks of the Grunsky operator vanish. For the off-diagonal blocks, we have:
        \begin{align}
            B^{(12)}:\overline{A_2(\DD)}&\to A_2(\DD),\notag\\
            u&\mapsto\left(z\mapsto-\frac{r}{\pi}\int_{\DD}\frac{\overline{u(w)}}{(1-rzw)^2}\dd^2w=-\frac{1}{\pi}\int_{\DD}\sum_{k=1}^\infty kr^{k}z^{k-1}w^{k-1}\overline{u(w)}\dd^2w\right).
        \end{align}
        Under the orthonormal basis of $A_2^1(\DD)$ as $\{e_n(z)=\sqrt{\frac{n}{\pi}}z^{n-1}\}_{n\in\N^*}$, we have that $B^{(12)}$ and $B^{(21)}$ are diagonal with entries $\{-r^k\}_{k\in\N^*}$. Therefore, $B_{\vec{g}}B^*_{\vec{g}}$ is diagonal, with eigenvalues $\{r^{2k}\}_{k\in\N^*}$, each with multiplicity $2$. Therefore,
        \begin{equation}
            \log D_F(\partial\A(r,1))=2\sum_{k=1}^\infty\log(1-r^{2k}).
        \end{equation}


\section{The uniformization of planar domains}\label{sec7}
In this section, we discuss the relations to uniformization theorem of our domain functionals. We first assume that we \textbf{do not know } Koebe's uniformization theorem of finitely connected planar domains and prove it using the Loewner potential. Then we discuss possible approach towards infinitely connected case via the renormalized volume.
\subsection{The Koebe's uniformization theorem via Loewner potential}\label{sec7.1}
\begin{lem}\label{lem3.9}
    For any $n$-connected domain $\Omega\subset\widehat{\CC}$, there is a domain $\tilde{\Omega}$ with analytic boundary, such that it is  conformally equivalent to $\Omega$.
\end{lem} 
\begin{proof}
    Apply the Koebe's construction  described in \Cref{sec3.2} to $\Omega$, at the first step, $\gamma_1$ is mapped to a round circle $\gamma_1^{(1)}$. Then at the second step, since $\varphi^{(2)}$ is conformal on $\gamma_1^{(1)}$, $\gamma_1^{(2)}=\varphi^{(2)}(\gamma_1^{(1)})$ keeps to be analytic. So after $n$ steps, every boundary curve is deformed to be analytic. 
\end{proof}
Now we may assume that $\Omega\subset\widehat{\CC}$ is an $n$-connected domain with analytic boundary, such that $\infty\in\Omega$. Let $\mathcal{F}(\Omega)$ to be the space of conformal maps on $\Omega$, such that $f(z)=z+O(z^{-1})$ near $\infty$. For $f\in\mathcal{F}(\Omega)$, let $f(\Omega)=\Omega_f$, $\widehat{\CC}\setminus\overline{\Omega_f}=D_f=\cup_{j=1}^nD_{f,j}$, and $\partial D_{f,j}=f(\gamma_j)=\gamma_{f,j}$
\begin{lem}\label{lem3.10}
    There exists $R=R(\Omega)>0$, such that for any $f\in \mathcal{F}(\Omega)$, $\overline{D_f}\subset B_R(0)$.
\end{lem}
\begin{proof}
    By assumption, $D=\widehat{\CC}\setminus\overline{\Omega}$ is bounded, thus there is $M=M(\Omega)>0$, such that $\overline{B_M(0)}\subset\Omega$. For $f\in \mathcal{F}(\Omega)$, let $f_M(w)=M^{-1}f(Mw)$, then $f_M$ is a univalent function on $\DD^*$. Suppose that $f_M(w)=w+\sum_{k=1}^\infty a_kw^{-k}$, then by the area theorem (see e.g. \cite[Theorem 2.1.2]{Lehto1987}), we have $\sum_{k=1}^\infty k|a_k|^2\le 1$. Then for $|w|=2$, by Cauchy--Schwarz inequality:
    \begin{align}
        |f_M(w)-w|&\le \left|\sum_{k=1}^\infty a_kw^{-k}\right|\notag\\
        &\le\left(\sum_{k=1}^\infty k|a_k|^2\right)^{1/2}\left(\sum_{k=1}^\infty\frac{4^{-k}}{k}\right)^{1/2}\le \sqrt{\log(4/3)}.
    \end{align}
    Therefore, for $|z|=2M$, $|f(z)|\le (2+\sqrt{\log(4/3)})M< 3M$. We may take $R=3M$, then $\overline{D_f}\subset B_R(0)$.
\end{proof}
To proceed, we recall some potential theory which can be found in \cite{GM08}. Let $G_\Omega$ be the Green's function of $\Omega$, where we use the normalization $-\D G_\Omega(\cdot,z)=\delta_z$. The \emph{harmonic measure} of $\Omega$ from $\infty$ is a probability measure $\omega_\Omega(\cdot,\infty)$on $\partial\Omega$, which can be given explicitly as:
\begin{equation}
    \omega_\Omega(\dd z,\infty)=-\frac{1}{2\pi}\partial_{\vec{n}}G_\Omega(z,\infty)|\dd z|,
\end{equation}
where $\vec{n}$ is the outer normal vector at $z\in\partial\Omega$. Let $\mathcal{P}(\partial\Omega)$ be the space of probability measures on $\partial\Omega$. Then the logarithmic energy for $\mu\in\mathcal{P}(\partial\Omega)$ is given by:
\begin{equation}
    E_{\partial\Omega}(\mu)=\int_{\partial\Omega}\int_{\partial\Omega}\log\frac{1}{|z-w|}\mu(\dd z)\mu(\dd w).
\end{equation}
The harmonic measure of $\Omega$ from $\infty$ is also called the equilibrium measure of $\partial\Omega$, as it is the unique minimizer of $E_{\partial\Omega}$ over $\mathcal{P}(\partial\Omega)$. We use the notation $\mu_{\partial\Omega}=\omega_\Omega(\cdot,\infty)$, and define the \emph{logarithmic capacity} of $\partial\Omega$ as $\mathrm{Cap}(\partial\Omega)=\exp(-E_{\partial\Omega}(\mu_{\partial\Omega}))$. For every boundary curve $\gamma_j$, let $\alpha_j=\mu_{\partial\Omega}(\gamma_j)>0$. The conformal invariance of Green's function, and our assumption that $f(\infty)=\infty$ implies that:
\begin{equation}
    \alpha_j=\omega_{\Omega}(\gamma_j,\infty)=\omega_{\Omega_f}(\gamma_{f,j},\infty)=\mu_{\partial\Omega_f}(\gamma_{f,j}).
\end{equation}
\begin{lem}\label{lem7.3}
    There exists $\beta=\beta(\Omega)>0$, such that for every $1\le j\le n$ and $f\in \mathcal{F}(\Omega)$, $\mathrm{Cap}(\gamma_{f,j})\ge\beta$.
\end{lem}
\begin{proof}
    Let $\mu_{f,j}=\alpha_j^{-1}\mu_{\partial\Omega_f}|_{f(\gamma_j)}$, then it is a probability measure on $\gamma_{f,j}$. Consider \Cref{lem3.10}, we have that for any $z,w\in \partial\Omega_f$, $|z-w|\le 2R$. Thus we have:
    \begin{align}
        \log\frac{2R}{\mathrm{Cap}(\partial\Omega_f)}&=\int_{\partial\Omega_f}\int_{\partial\Omega}\log\frac{2R}{|z-w|}\mu_{f}(\dd z)\mu_f(\dd w)\notag\\
        &\ge\alpha_j^2\int_{\gamma_{f,j}}\int_{\gamma_{f,j}}\log\frac{2R}{|z-w|}\mu_{f,j}(\dd z)\mu_{f,j}(\dd w)\notag\\
        &\ge\alpha_j^2\log\frac{2R}{\mathrm{Cap}(\gamma_{f,j})}.
    \end{align}
    Therefore, we have $\mathrm{Cap}(\gamma_{f,j})\ge\max_{1\le i\le n}2R(C/2R)^{\alpha_j^{-2}}>0$.
\end{proof}
We consider the conformal metric $g=e^\phi |\dd z|^2$ on $\widehat{\CC}$, such that $\phi\equiv 0$ on $B_R(0)$, and $\phi=4/(1+|z|^2)^2$ near $\infty$. Then in particular we know that $g$ is flat on $D_f$ for every $f\in\mathcal{F}(\Omega)$. 
\begin{prop}\label{prop7.4}
    Let $D_{f,j}^*=\widehat{\CC}\setminus\overline{D_{f,j}}$, and let $0<\lambda_{1,f,j}\le\lambda_{2,f,j}\le\dots$ be the spectrum of Laplace--Beltrami operator $-\D_{f,j}$ on $D_{f,j}$. Then there is a universal spectral gap which holds for every $f\in\mathcal{F}(\Omega)$, namely, there is $\lambda>0$ independent of $f$, such that:
    \begin{equation}
        \lambda_{1,f,j}\ge \lambda,\quad \forall 1\le j\le n, f\in\mathcal{F}(\Omega).
    \end{equation}
\end{prop}
\begin{proof}
    By the variational principle for first eigenvalue (see e.g. \cite[Theorem 6.5.2]{Eva98}), and the density of $C_0^\infty(D_{f,j}^*)$ in $H_0^1(D_{f,j}^*)$, we only need to show that there is $\lambda>0$ independent of $f$, such that for every $u\in C^\infty_0(D_{f,j}^*)$, 
    \begin{equation}
        \lambda \int_{D_{f,j}^*}|u|^2\dd^2z\le \int_{D_{f,j}^*}|\n u|^2\dd^2z.
    \end{equation}
    To do this, consider the Green's function $G_g$ of $\widehat{\CC}$ with respect to $g$, and we normalize it as $-\D G_g(\cdot,z)=\delta_z-1/A_g$, where $A_g=\mathrm{Area}_g(\widehat{\CC}).$  Since $g$ is flat on $B_R(0)$, $G_g(z,w)=-\log|z-w|+U_g(z,w)$ when $z,w\in B_R(0)$, and $U_g$ is bounded on $B_R(0)$. Let $\rho_u=(\D u)\dd A_g$, since $u$ is compactly supported, we have $\rho_u(\widehat{\CC})=0$. Let $\theta=\mu_{f,j}-A_g^{-1}\dd A_g$, then $\theta(\widehat{\CC})=0$ as well. Thus for any $t\in\R$, $(\theta+t\rho_u)(\widehat{\CC})=0$. By \cite{GM08},
    \begin{align}\label{eq163}
       0\le E^G_g(\theta+t\rho_u)=\int_{\widehat{\CC}}\int_{\widehat{\CC}}G_g(z,w)(\theta+t\rho_u)(\dd z)(\theta+t\rho_u)(\dd w). 
    \end{align}
    We compute the right hand side of the equation above in the following.
    By \Cref{lem7.3}, we have:
    \begin{align}
        \int_{\widehat{\CC}}\int_{\widehat{\CC}}G_g(z,w)\theta(\dd z)\theta(\dd w)&=\int_{\widehat{\CC}}\int_{\widehat{\CC}}G_g(z,w)\mu_{f,j}(\dd z)\mu_{f,j}(\dd w)\notag\\
        &=\int_{\gamma_{f,j}}\int_{\gamma_{f,j}}\left(\log\frac{1}{|z-w|}+U_g(z,w)\right)\mu_{f,j}(\dd z)\mu_{f,j}(\dd w)\notag\\
        &\le-\log\beta+\sup_{z,w\in B_R(0)}|U_g(z,w)|\le Q,
    \end{align}
    where $Q=Q(g,R,\beta)$ is a constant independent of $f$. For the $\rho_u$ term, we have:
    \begin{equation}
        \int_{\widehat{\CC}}\int_{\widehat{\CC}}G_g(z,w)\rho_u(\dd z)\rho_u(\dd w)=\int_{\widehat{\CC}}\left(\frac{1}{A_g}-u(w)\right)\rho_u(\dd w)=\int_{D_{f,j}^*}|\n u|^2\dd^2z.
    \end{equation}
     
    For the mixed term, if we let $\overline{u}=A^{-1}_g\int_{\widehat{\CC}}u\dd A_g$, then:
    \begin{align}
        \int_{\widehat{\CC}}\int_{\widehat{\CC}}G_g(z,w)\rho_u(\dd z)\theta(\dd w)&=\int_{\widehat{\CC}}\left(\frac{1}{A_g}-u(w)\right)\theta(\dd w)=\overline{u}.
    \end{align}
    Then \Cref{eq163} implies that for any $t\in\R$,
    \begin{equation}
        0\le t^2\int_{D_{f,j}^*}|\n u|^2+2t\overline{u}+Q.
    \end{equation}
    Therefore, the discriminant of the quadratic equation above is non-positive, i.e.,
    \begin{equation}
        \overline{u}^2\le Q\int_{D_{f,j}^*}|\n u|^2.
    \end{equation}
    Finally, by Poincaré inequality on $B_R(0)$, we have:
    \begin{align}
        \int_{D_{f,j}^*}|u|^2&\le 2\int_{D_{f,j}^*}(|\overline{u}|^2+|u-\overline{u}|^2)\notag\\
        &\le 2A_g|\overline{u}|^2+ C\int_{D_{f,j}^*}|\n u|^2\notag\\
        &\le (2A_gQ+C)\int_{D_{f,j}^*}|\n u|^2,
    \end{align}
    where $C$ is a constant independent of $f$.
\end{proof}
We consider the definition of Loewner potentials via zeta function \Cref{eq34}. For simplification of notations, we use $\Omega'=\Omega_f$, $D_j'=D_{f,j}$ temporarily. Then we have:
\begin{align}\label{eq170}
    \mathcal{H}_{n,g}(\Omega')-\sum_{j=1}^n\mathcal{H}_{g}(\gamma_j')&=\log\frac{\prod_{j=1}^n\det_\zeta(-\D_{g,D_j'^*})}{\det_\zeta(-\D_{g,\Omega'})\det^{n-1}_\zeta(-\D_{g,\widehat{\CC}})\mathrm{Area}^{n-1}_g(\widehat{\CC})}\notag\\
    &=\left.\frac{\dd}{\dd s}\right|_{s=0}\left((n-1)\zeta_{-\D_{g,\widehat{\CC}}}(s)+\zeta_{-\D_{g,\Omega'}}(s)-\sum_{j=1}^n\zeta_{-\D_{g,D_j'^*}}(s)\right)+\mathrm{Const}.\notag\\
    &=\left.\frac{\dd}{\dd s}\right|_{s=0}\frac{1}{\Gamma(s)}\int_0^\infty t^{s-1}(T_g(t)-(n-1))\dd t+\mathrm{Const}.
\end{align}
Here
\begin{equation}
    T_g(t)=(n-1)\mathrm{tr}(e^{t\D_{g,\widehat{\CC}}})+\mathrm{tr}(e^{t\D_{g,\Omega'}})-\sum_{j=1}^n\mathrm{tr}(e^{t\D_{g,D_j'^*}}).
\end{equation}
\begin{lem}\label{lem7.5}
    With the notations above, for every $t\ge 0$, $T_g(t)\ge 0$.
\end{lem}
\begin{proof}
    Note that the heat trace $\mathrm{tr}(e^{t\D_{g,M}})$ is the total mass of all the rooted Brownian loops on $M$ which return to its starting point at time $2t$ (see e.g. \cite[Section 3]{LW04}). Now consider one such loop in $\widehat{\CC}$ which intersects $k$ of the sets $\{\overline{D'}_j\}_{1\le j\le n}$, $0\le k\le n$. Then this loop is completely contained in $n-k$ of the sets $\{D_j'^*\}_{1\le j\le n}$, and when $k=0$, it is also completely contained in $\Omega'$. Thus for this loop, in the expression of $T_g(t)$ above, it is counted $(n-1)+\mathbbm{1}_{\{k=0\}}-(n-k)=k-1+\mathbbm{1}_{\{k=0\}}$ times, which is nonnegative. 
\end{proof}
\begin{prop}\label{prop3.11}
    There is a constant $C=C(\Omega)$, such that
    \begin{equation}
        \inf_{f\in\mathcal{F}(\Omega)}\mathcal{H}_{n,g}(\Omega_f)\ge C.
    \end{equation}
\end{prop}
\begin{proof}
    We show that there is a constant $C_1=C_1(\Omega)$, such that $\mathcal{H}_{n,g}(\Omega')-\sum_{j=1}^n\mathcal{H}_g(\gamma_j')\ge C_1$. Together with this, since $I^L(\gamma_j')=12(\mathcal{H}_g(\gamma_j')-\mathcal{H}_g(\mathbb{S}^1))\ge 0$, we will have that:
    \begin{equation}
        \mathcal{H}_{n,g}(\Omega')\ge C_1+n\mathcal{H}_{g}(\mathbb{S}^1)=C.
    \end{equation}
    To do this, consider the zeta function expression in \Cref{eq170}. Since $1/\Gamma(s)=s+\upgamma s^2+o(s^2)$ near $0$, where $\upgamma$ is the Euler constant, we have:
    \begin{align}
        &\frac{1}{\Gamma(s)}\int_0^\infty t^{s-1}(T_g(t)-(n-1))\dd t\notag\\
        =&\frac{1}{\Gamma(s)}\left(\int_0^1t^{s}\frac{{}T_g(t)}{t}\dd t+\int_0^1t^{s-1}(1-n)\dd t\right)+\frac{1}{\Gamma(s)}\int_1^\infty t^{s}\frac{T_g(t)-(n-1)}{t}\dd t\notag\\
        =&s\left(\int_0^1\frac{T_g(t)}{t}\dd t+(1-n)\upgamma+\int_1^\infty\frac{T_g(t)-(n-1)}{t}\dd t\right)+o(s).
    \end{align}
    Therefore,
    \begin{equation}
        \mathcal{H}_{n,g}(\Omega')-\sum_{j=1}^n\mathcal{H}_g(\gamma_j')=\int_0^1\frac{T_g(t)}{t}\dd t+\int_1^\infty\frac{T_g(t)-(n-1)}{t}\dd t+\mathrm{Const}.
    \end{equation}
    By \Cref{prop7.4}, we have that:
    \begin{align}
        \mathrm{tr}(e^{t\D_{g,D_j'^*}})=\sum_{k=1}^\infty e^{-t\lambda_{k}'}=\sum_{k=1}^\infty e^{-(t-1)\lambda_k'}e^{-\lambda_k'}\le e^{-(t-1)\lambda}\mathrm{tr}(e^{\D_{g,D_j'^*}}).
    \end{align}
    Then we have:
    \begin{align}
        T_g(t)&=(n-1)\mathrm{tr}(e^{t\D_{g,\widehat{\CC}}})+\mathrm{tr}(e^{t\D_{g,\Omega'}})-\sum_{j=1}^n\mathrm{tr}(e^{t\D_{g,D_j'^*}})\notag\\
        &\ge (n-1)-e^{-(t-1)\lambda}\sum_{j=1}^n\mathrm{tr}(e^{\D_{g,D_j'^*}})\notag\\
        &\ge(n-1)-ne^{-(t-1)\lambda}\mathrm{tr}(e^{\D_{g,\widehat{\CC}}}).
    \end{align}
    By \Cref{lem7.5}, we have that:
    \begin{align}
        \mathcal{H}_{n,g}(\Omega')-\sum_{j=1}^n\mathcal{H}_g(\gamma_j')&=\int_0^1\frac{T_g(t)}{t}\dd t+\int_1^\infty\frac{T_g(t)-(n-1)}{t}\dd t+\mathrm{Const}.\notag\\
        &\ge-\int_1^\infty\frac{ne^{-(t-1)\lambda}\mathrm{tr}(e^{\D_{g,\widehat{\CC}}})}{t}\dd t+\mathrm{Const}.\notag\\
        &=-\frac{n\cdot\mathrm{tr}(e^{\D_{g,\widehat{\CC}}})}{\lambda}+\mathrm{Const}.
    \end{align}
    Therefore, the final result is a constant independent of $f$.
\end{proof}
\begin{thm}\label{thm7.7}
    There is a conformal map $f_0\in\mathcal{F}(\Omega)$, such that
    $\Omega_{f_0}$ is a circle domain. Thus any finitely connected domain can always be conformally uniformized into a circle domain. 
\end{thm}
\begin{proof}
    By \Cref{lem3.9}, we may assume that $\Omega$ has analytic boundary. We apply Koebe's construction (see \Cref{sec3.2}) to $\Omega$ but change a little bit. Let $\Omega=\Omega^{(0)}$. Assume that we've already constructed $\Omega^{(k)}=f^{(k)}(\Omega)$, where $f^{(k)}\in\mathcal{F}(\Omega)$. For $\Omega^{(k)}$, let $1\le m_k\le n$ be the index such that $I^L(\gamma_{m_k}^{(k)})=\max_{1\le j\le n}I^L(\gamma_j^{(k)})$. Then define $\Omega^{(k+1)}=\varphi^{(k)}(\Omega^{(k)})$, where $\varphi^{(k)}:\widehat{\CC}\setminus D_{m_k}^{(k)}\to \DD^*$ conformally such that $\varphi^{(k)}(z)=z+O(z^{-1})$ at $\infty$, and let $f^{(k+1)}=\varphi^{(k)}\circ f^{(k)}$. Note that since the sequence $\{f^{(k)}\}_{k\in\N}\subset\mathcal{F}(\Omega)$, they form a normal family. By passing to a subsequence we may assume that $f^{(k)}\to f\in\mathcal{F}(\Omega)$ locally uniformly in $\CC$. By \Cref{lem3.10}, for each $1\le j\le n$, $\{\gamma_j^{(k)}\}_{k\in\N}$ converges uniformly to a curve $\gamma_j^{(\infty)}$. By the proof of \Cref{thm:lpnonnegative} and \Cref{prop3.11}, we have that:
    \begin{equation}
        \sum_{k=0}^\infty\max_{1\le j\le n}I^L(\gamma_j^{(k)})=\lim_{N\to+\infty}\sum_{k=0}^N12(\mathcal{H}_{n,g}(\Omega^{(k)})-\mathcal{H}_{n,g}(\Omega^{(k+1)}))\le 12(\mathcal{H}_{n,g}(\Omega)-C)<\infty.
    \end{equation}
    Therefore, for each $1\le j\le n$, the sequence $\{I^L(\gamma_j^{(k)})\}_{k\in\N}$ converges to $0$ as $k\to\infty$. By \cite[Lemma 2.12]{RW21}, we have that $I^L(\gamma_j^{(\infty)})=0$. This implies that $\Omega^{(\infty)}=f(\Omega)$ is a circle domain.
\end{proof}
\subsection{Towards infinitely connected uniformization}\label{sec7.2}
Let $\Omega$ be countably connected, that is to say, $D=\widehat{\CC}\setminus\overline{\Omega}$ has countably many components. We always assume that  $\infty\in\Omega$. By He--Schramm's uniformization theorem \cite{HeSchramm1993}, there is a circle domain $\Omega_c$ with a conformal map $f_c:\Omega_c\to\Omega$, such that $f_c(z)=z+O(z^{-1})$ near $\infty$. We label all the components of $D$ as $\{D_j\}_{j\in\N^*}$. Let $J\subset\N^*$ be a finite set, we define $\Omega_J=\widehat{\CC}\setminus\cup_{j\in J}\overline{D_j}$, and a conformal map $f_{c,J}:\Omega_{c,J}\to\Omega_J$ such that $f_{c,J}(z)=z+O(z^{-1})$. 
\begin{defn}
    Assume that every boundary component of $\Omega$ has $C^{5,\alpha}$ regularity. We define the countably connected renormalized volume of $\Omega$ as:
    \begin{equation}
        V_\infty(\Omega)=\sup_{J}(V_R(\Omega_J)-V_R(\Omega_{c,J})),
    \end{equation}
    where the supremum is taken over all $J\subset\N^*,\#J<\infty$.
\end{defn}
\begin{lem}\label{lem7.1}
    If $J,K\subset \N^*$, such that $J\subset K$, then we have:
    \begin{equation}
        V_R(\Omega_J)-V_R(\Omega_{c,J})\le V_R(\Omega_K)-V_R(\Omega_{c,K}).
    \end{equation}
    Thus for any increasing exhaustion $\{J_l\}_{l\in \N^*}$ of $\N^*$, we have $V_\infty(\Omega)=\lim_{l\to\infty}(V_R(\Omega_{J_l})-V_R(\Omega_{c,J_l}))$. 
\end{lem}
\begin{figure}[h]
        \centering
        \includegraphics[width=0.75\textwidth]{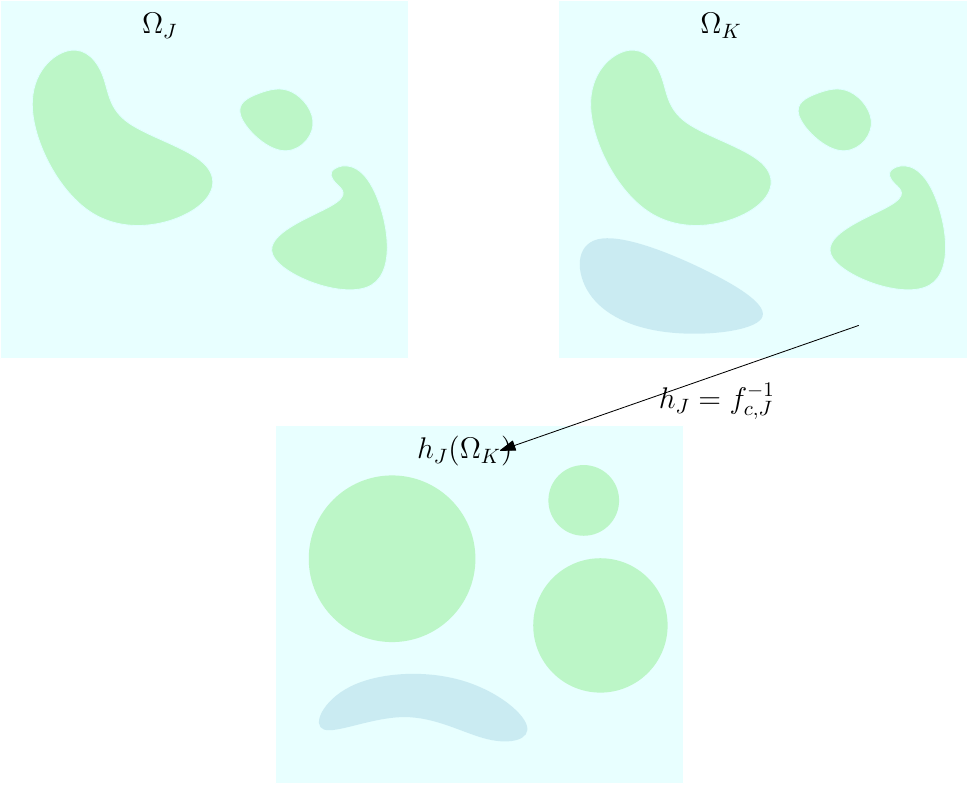}
        \caption{Illustration of the domains and conformal map used in \Cref{lem7.1}.}
        \label{fig:Multiply}
    \end{figure}
\begin{proof}
    Assume that $\#J=n_J,\#K=n_k$. Consider the conformal map $h_J=f_{c,J}^{-1}:\Omega_{J}\to\Omega_{c,J}$, we have for any metric $g$ conformal to the spherical metric:
    \begin{align}
       \mathcal{H}_{n_K,g}& (h_J(\Omega_K))-\mathcal{H}_{n_K,g}(\Omega_K) \notag\\
        &=\log\frac{\det_\zeta(-\D_{g,\Omega_K})\prod_{k\in K}\det_\zeta(-\D_{g,D_k})}{\det_{\zeta}(-\D_{g,h_J(\Omega_K)})\prod_{i\in K\setminus J}\det_\zeta(-\D_{g,h_J(D_i)})\prod_{j\in J}\det_\zeta(-\D_{g,B_j})}\notag\\
        &=\mathcal{H}_{n_J,g}(\Omega_{c,J})-\mathcal{H}_{n_J,g}(\Omega_J)+\log\frac{\det_\zeta(-\D_{g,\Omega_K})}{\det_\zeta(-\D_{g,h_J(\Omega_K)})}\notag\\
        &\quad +\log\frac{\det_\zeta(-\D_{g,\Omega_{c,J}})\prod_{k\in K}\det_\zeta(-\D_{g,D_k})}{\det_\zeta(-\D_{g,\Omega_J})\prod_{i\in K\setminus J}\det_\zeta(-\D_{g,h_J(D_i)})\prod_{j\in J}\det_\zeta(-\D_{g,D_j})}.
    \end{align}
    Here we assume that $B_1,\dots,B_{n_J}$ are the round disk complements of $\Omega_{J}$.
    We prove that the zeta-regularized determinants parts in the equation above vanish. Indeed, since $h$ is conformal on $\Omega_K,\Omega_J$, and $D_i$ for $i\in K\setminus J$, then:
    \begin{equation}
        {\det}_\zeta(-\D_{g,h_J(D_i)})={\det}_\zeta(-\D_{h_J^*g,D_i}),\quad i\in K\setminus J,
    \end{equation}
    and the same formula also holds for $h(\Omega_K)$ and $h(\Omega_J)=\Omega_{c,J}$. Then by Polyakov--Alvarez formula (\Cref{eq36}), all the corresponding boundary terms will cancel, and we thus have:
    \begin{equation}\label{eq167}
        \mathcal{H}_{n_K,g}(h(\Omega_K))-\mathcal{H}_{n_K,g}(\Omega_K)=\mathcal{H}_{n_J,g}(\Omega_{c,J})-\mathcal{H}_{n_J,g}(\Omega_J).
    \end{equation}
    Then by \Cref{thm:lpnonnegative}, \Cref{cor5.7}  we have:
    \begin{align}
        V_R(\Omega_J)-V_R(\Omega_{c,J})&=3\pi(\mathcal{H}_{n_J,g}(\Omega_J)-\mathcal{H}_{n_J,g}(\Omega_{c,J}))\notag\\
        &=3\pi(\mathcal{H}_{n_K,g}(\Omega_K)-\mathcal{H}_{n_K,g}(\Omega_{c,K})\notag\\
        &-(\mathcal{H}_{n_K,g}(h(\Omega_K))-\mathcal{H}_{n_K,g}(\Omega_{c,K})))\notag\\
        &\le V_R(\Omega_K)-V_R(\Omega_{c,K}).
    \end{align}
\end{proof}
\begin{prop}\label{prop7.9}
    Let $w_t:\widehat{\CC}\to \widehat{\CC}$ be a family of quasiconformal maps, such that $w_t(z)=z+O(z^{-1})$ near $\infty$, and the Beltrami differential $\mu_t=t\dot{\mu}+o(t)$ is supported on  $D_J=\{D_j\}_{j\in J}$, where $J$ is finite. Then we have:
    \begin{equation}
        \delta V_\infty(\Omega)=\left.\frac{\dd}{\dd t}\right|_{t=0}V_\infty(w_t(\Omega))=-\mathrm{Re}\sum_{j\in J}\int_{D_j}\dot{\mu}_j\cdot S[g_j^{-1}]\dd^2z.
    \end{equation}
    Moreover, if $\Omega$ is a critical point for every such variation, then $\Omega$ is a circle domain.
\end{prop}
\begin{proof}
    Consider any $K\subset\N^*$ such that $J\subset K$. Since 
      $w_t$ is conformal on $\Omega_J$ and $\Omega_K\subset\Omega_J$, by the proof of \Cref{lem7.1}, we have:
      \begin{equation}
          V_R(\Omega_K)-V_R(w_t(\Omega_K))=V_R(\Omega_{J})-V_R(w_t(\Omega_J)),
      \end{equation}
      which is similar to \Cref{eq167}. By \Cref{lem7.1}, we have that:
      \begin{align}
          V_\infty(w_t(\Omega))&=\sup_K(V_R(w_t(\Omega_K))-V_R(\Omega_{c,K}))\notag\\
          &=\sup_K(V_R(w_t(\Omega_K))-V_R(\Omega_K)+V_R(\Omega_K)-V_R(\Omega_{K,c}))\notag\\
          &=V_R(w_t(\Omega_J))-V_R(\Omega_J)+V_\infty(\Omega)
      \end{align}
      Thus apply \Cref{thm:varre} to the equation above, we can get the variation formula. If $\Omega$ is a critical point for this variation, then $D_J$ is a union of disks. 
\end{proof}
\begin{rmk}
    Note that no matter $V_\infty(\Omega)$ is finite or not, we have that $V_\infty(w_t(\Omega))-V_\infty(\Omega)=V_R(w_t(\Omega_J))-V_R(\Omega_J)$ is finite. Thus the variation always makes sense.
\end{rmk}
    \appendix
    \section{Explicit computation of Epstein surface}\label{appendixA}
    In $\HH^3=\CC\times\R^+$, consider the polar coordinate $x=(\rho e^{i\theta},h)$. So we have the variables \begin{equation}
        \mathrm{Ep}_{\A(r,1)}(z)=\mathrm{Ep}_{\A(r,1)}(se^{i\tau})=(\rho(s,\tau)e^{i\theta(s,\tau)},h(s,\tau)).
    \end{equation}
    For simplicity of notation, we introduce $v=\pi\log s/\log r\in(0,\pi)$, $u=\cot v\in \R$, and $\alpha=-(\log r)/\pi$. On $\A(r,1)$, the complete hyperbolic metric is rotationally symmetric:
    \begin{equation}
        g_r(z)=\left(\frac{\pi}{\log r}\right)^2\frac{1}{|z|^2\sin^2\left(\pi\frac{\log |z|}{\log r}\right)}|dz|^2=e^{2\alpha v}\frac{1}{\alpha^2\sin^2 v}|dz|^2.    
    \end{equation}
    \begin{lem}\label{lem:A.1}
        For simplicity of notation, let $v=\pi\log s/\log r\in(0,\pi)$, $\alpha=-\log r/\pi$, then:
        \begin{align}
            \rho&=\rho(v)=e^{-\alpha v}\cdot\frac{1+\cos^2v-\alpha^2\sin^2v}{1+\left(\cos v-\alpha\sin v\right)^2},\notag\\
            h&=h(v)=e^{-\alpha v}\cdot\frac{2\alpha\sin v}{1+\left(\cos v-\alpha\sin v\right)^2}.\notag
        \end{align}
    \end{lem}
    \begin{prop}
        For $\A(r,1)$,
        \begin{equation}
            V(\A(r,1))=-\frac{1}{8\pi}(\log r)^3+\frac{\pi}{4}\log r-\frac{5\pi^3}{8\log r}.
        \end{equation}
    \end{prop}
    \begin{proof}
        Consider a 2-form in $\HH^3$:
        \begin{equation}
            \omega=\frac{x\dd y-y\dd x}{2h^3}\wedge \dd h,
        \end{equation}
        then we can verify that $\dd\omega=\dd\mathrm{Vol}_{\HH^3}$. To make our computation rigorous, let $\epsilon>0$ small, and take a truncation $h>\epsilon$. Let $\Sigma_{r\DD,\epsilon},\Sigma_{\DD^*,\epsilon},\Sigma_{\A_r,\epsilon}$ be the truncated Epstein surfaces. Thus by Stokes' formula, we have:
        \begin{equation}
            V(\A(r,1))=\lim_{\epsilon\to 0}\int_{\Sigma_{\A_r,\epsilon}+\Sigma_{r\DD,\epsilon}-\Sigma_{\DD^*,\epsilon}}\omega=\lim_{\epsilon\to0}\pi\int_{\Sigma_{\A_r,\epsilon}+\Sigma_{r\DD,\epsilon}-\Sigma_{\DD^*,\epsilon}}\frac{\rho^2(s)h'(s)}{h^3(s)}\dd s.
        \end{equation}
        The integral on two hemispheres can be easily computed:
        \begin{align}
            I_{\Sigma_{\DD^*}}(\epsilon)&=\pi\int_0^{\sqrt{1-\epsilon^2}}\frac{-s^3}{(1-s^2)^2}\dd s=\pi\left(-\log\epsilon-\frac{1}{2\epsilon^2}+\frac{1}{2}\right),\\
            I_{\Sigma_{r\DD}}(\epsilon)&=\pi\int_0^{\sqrt{r^2-\epsilon^2}}\frac{-s^3}{(r^2-s^2)^2}\dd s=\pi\left(\log\frac{r}{\epsilon}-\frac{r^2}{2\epsilon^2}+\frac{1}{2}\right).
        \end{align}
        For the annulus part, we do the change of variable for $u=\cot v=\cot(\pi\log s/\log r)$. Let $u_\epsilon^+>u_\epsilon^-$ be the two solutions of $h(u)=\epsilon$. Using the previous computation in \Cref{lem:A.1}, we have:
        \begin{align}
            \frac{I_{\A(r,1)}(\epsilon)}{\pi}&=\int_{u_\epsilon^-}^{u_\epsilon^+}\left(\frac{\rho(u)}{h(u)}\right)^2\frac{\dd\log h}{\dd u}\dd u\notag\\
            &=\int_{u_\epsilon^-}^{u_\epsilon^+}\frac{(2u^2+1-\alpha^2)^2}{4\alpha^2(1+u^2)}\left(\frac{u+\alpha}{u^2+1}-\frac{4u-2\alpha}{2u^2-2\alpha u+\alpha^2+1}\right)\dd u\notag\\
            &=\log(2u^2-2\alpha u+\alpha^2+1)-\frac{\log(u^2+1)}{2}+\frac{\alpha^5-2\alpha^3+5\alpha}{8\alpha^2}\arctan u\notag\\
            &-\left.\frac{u^2}{2\alpha^2}+\frac{(\alpha^5-2\alpha^3+\alpha)u-\alpha^4-2\alpha^2-1}{8\alpha^2(u^2+1)}\right|^{u_\epsilon^+}_{u_\epsilon^-}.
        \end{align}
        To get the final result, we analyze the asymptotic behaviour of $u_\epsilon^\pm$.
        \begin{itemize}
            \item When $v\to 0$, $h(v)=\alpha v+o(v^2)$. Thus $v\sim\epsilon/\alpha$. By $u(v)=1/v-v/3+o(v^2)$ we have $u_\epsilon^+=\alpha/\epsilon+O(1)$.
            \item When $v\to\pi$, let $w=\pi-v$, then $h(w)=e^{-\alpha\pi}\alpha w+o(w^2)$. Thus $w\sim\epsilon e^{\alpha\pi}/\alpha$. By $u(w)=\cot(\pi-w)=-\cot w=-1/w+w/3+o(w^2)$ we have: $u_\epsilon^-=-\alpha e^{-\alpha\pi}/\epsilon+O(1)$.
        \end{itemize}
        Using the asymptotics of $u_\epsilon^\pm$, we have:
        \begin{align}
            \frac{I_{\A(r,1)}(\epsilon)}{\pi}=\alpha\pi+\frac{\alpha^5-2\alpha^3+5\alpha}{8\alpha^2}\pi-\frac{1-e^{-2\alpha\pi}}{2\epsilon^2}+o(1).
        \end{align}
        To sum up, we have:
        \begin{align}
            \frac{V(\A(r,1))}{\pi}&=\lim_{\epsilon\to 0}\left(\log\frac{r}{\epsilon}-\frac{r^2}{2\epsilon^2}+\frac{1}{2}\right)-\left(-\log\epsilon-\frac{1}{2\epsilon^2}+\frac{1}{2}\right)\notag\\
            &+\alpha\pi+\frac{\alpha^5-2\alpha^3+5\alpha}{8\alpha^2}\pi-\frac{1-e^{-2\alpha\pi}}{2\epsilon^2}\notag\\
            &=-\frac{(\log r)^3}{8\pi^2}+\frac{\log r}{4}-\frac{5\pi^2}{8\log r}.
        \end{align}
        Note that in the last step we've used the substitution $\alpha=-\log r/\pi$.
    \end{proof}
    \begin{lem}
        For $\A(r,1)$, 
        \begin{equation}
            \int_{\Sigma_{r\DD}\cup\Sigma_{\A(r,1)}\cup\Sigma_{\DD^*}}H\dd A=\int_{\A(r,1)}\frac{|\vartheta_r|^2}{\rho_{\A(r,1)}^2}\dd^2z=-\frac{1}{4\pi}(\log r)^3-\frac{\pi}{2}\log r-\frac{\pi^3}{4\log r}.
        \end{equation}
    \end{lem}
    \begin{proof}
        Since $\Sigma_{r\DD}$ and $\Sigma_{\DD^*}$ are totally geodesic, the mean curvature vanishes on them. On $\A(r,1)$, we have:
        \begin{equation}
            \vartheta_r=\frac{1}{2z^2}\left(1+\frac{\pi^2}{(\log r)^2}\right).
        \end{equation}
        Thus:
        \begin{align}
            \int_{\Sigma_{\A(r,1)}}H\dd A&=\frac{(\log r)^2}{4\pi^2}\left(1+\frac{\pi^2}{(\log r)^2}\right)^2\int_{\A(r,1)}\sin^2\left(\pi\frac{\log|z|}{\log r}\right)\frac{\dd^2z}{|z|^2}\notag\\
            &=-\frac{(\log r)^3}{4\pi}-\frac{\pi\log r}{2}-\frac{\pi^3}{4\log r}.
        \end{align}
    \end{proof}
    \begin{cor}
        \begin{equation}
            V_R(\A(r,1))=V(\A(r,1))-\frac{1}{2}\int_{\Sigma_{r\DD}\cup\Sigma_{\A(r,1)}\cup\Sigma_{\DD^*}}H\dd A=\frac{\pi}{2}\left(\log r-\frac{\pi^2}{\log r}\right).
        \end{equation}
    \end{cor}
    This verifies the result in \Cref{prop5.9} where we get from the variational formula.
    \section{Elliptic regularity of Beltrami equation}\label{appendixB}
    The aim of this appendix is to prove \Cref{lem5.1}.
    Let $\mu\in L^\infty_1(\CC)$, $w=u+iv$ is a solution to the Beltrami equation of $\mu$, then the Beltrami equation can be written as the following elliptic system (see e.g. \cite[Chapter 10]{AIM09}):
    \begin{equation}\label{eqb1}
        \begin{cases}
            \mathrm{div}(A_\mu\n u)=0,\\
            \mathrm{div}(A_\mu\n v)=0,
        \end{cases}\quad
        A_\mu=\frac{1}{1-|\mu|^2}\begin{pmatrix}
            |1-\mu|^2& -2\mathrm{Im}\mu\\
            -2\mathrm{Im}\mu&|1+\mu|^2
        \end{pmatrix}.
    \end{equation}
    We will use the following regularity result of this type of equation:
    \begin{lem}(\cite[Theorem 1]{Zhuge21})\label{lemb1}
        Let $D$ be a bounded Lipschitz domain in $\R^d$, $D^\pm$ are two disjoint subdomains, such that $D^+=D\setminus\overline{D^-}$. Let $S=\overline{D^+}\cap\overline{D^-}$ be the interface between $D^\pm$, such that $S$ is a $C^{k,\alpha}$ curve. Let $A^\pm\in C^{k-1,\alpha}(\overline{D^\pm},\R^{d\times d})$, $h^\pm\in C^{k-1,\alpha}(\overline{D^\pm},\R^d)$. If $u$ is a weak solution of the following equation:
        \begin{equation}
            \begin{cases}
                \mathrm{div}(A\n u)=\mathrm{div}(h),\quad &\text{in }D\setminus S=D^+\cup D^-,\\
                u^+=u^-, &\text{on }S,
            \end{cases}
        \end{equation}
        then $u\in C^{k,\alpha}(\overline{D^\pm})$. Moreover,
        \begin{equation}
            \sum_{\pm}\Vert u\Vert_{C^{k,\alpha}(\overline{D^\pm})}\le C\left(\sum_{\pm}\Vert h^\pm\Vert_{C^{k-1,\alpha}(\overline{D^\pm})}+\Vert u\Vert_{L^2(D)}\right),
        \end{equation}
        where $C$ only depends on $A,d,k,\alpha,S$. 
    \end{lem} 
    \begin{lem}\label{lemb2}
        Using the notations in \Cref{lem5.1}, then for any $0<\beta<\alpha$, we have that $t\mapsto w_t\in C^0(I,C^{5,\beta}(\overline{U}))$, where $U\in \{D_1,\dots, D_n,\Omega\}$.
    \end{lem}
    \begin{proof}
        We prove that for any $t_0\in I$, with $t_n\to t_0$, then $\Vert w_{t_n}-w_{t_0}\Vert_{C^{5,\beta}}\to 0$ as $n\to \infty$. To do this, consider a compact interval $I_1$ containing the sequence $\{t_n\}_{n\in \N}$, then by \cite[Theorem 7]{AB60} and \Cref{lemb1}, we have that:
        \begin{equation}
            \sup_{t\in I_1}\Vert w_t\Vert_{C^{5,\alpha}}\le C\sup_{t\in I_1}\Vert w_t\Vert_{L^2}<+\infty,
        \end{equation}
        on every $\overline{U}$. Since $C^{5,\alpha}\subset C^{5,\beta}$ is a compact embedding of Banach spaces, we have that for every subsequence of $\{w_{t_n}\}$, there is a further subsequence $\{w_{t_{n_k}}\}$, such that $w_{t_{n_k}}\to w$ in $C^{5,\beta}$. Since the inclusion of $C^{5,\beta}\subset L^2(\overline{U})$ is continuous, this convergence also happens in $L^2$. By \cite[Theorem 8]{AB60}, we have that $w_{t_n}\to w_{t_0}$ in $L^2(\overline{U})$ as well, so $w=w_{t_0}$ a.e. This implies that $w_{t_n}\to w_{t_0}$ in $C^{5,\beta}(\overline{U})$.
    \end{proof}
    \begin{lem}
        Suppose that $w_t=u_t+iv_t$, then there exists $t\mapsto\dot{u}_t\in C^0(I,C^{5,\beta}(\overline{U}))$ for any $0<\beta<\alpha$, such that:
        \begin{equation}
            \lim_{s\to 0}\left\Vert\frac{u_{t+s}-u_t}{s}-\dot{u}_t\right\Vert_{C^{5,\beta}(\overline{U})}= 0.
        \end{equation}
        The same conclusion holds for $v_t$ as well. 
    \end{lem}
    \begin{proof}
        Since $t\mapsto \mu_t\in C^1(I,C^{4,\alpha}(\overline{U}))\subset C^1(I,L^\infty(\overline{U}))$, by \cite[Theorem 10]{AB60}, we have that $t\mapsto u_t\in C^1(I,L^2(\overline{U}))$. Thus we may take the differential at $t$ as $\dot{u}_t$. We need to show that it is also the differential of $t\mapsto u_t$ in the $C^{5,\beta}$ topology. To do this, we first show that $t\mapsto u_t$ is differentiable. Let $A_t$ be the matrix of $\mu_t$ as in \Cref{eqb1}. Then since $t\mapsto\mu_t\in C^1(I,C^{4,\alpha}(\overline{U}))$, we have that $t\mapsto A_t\in C^1(I,C^{4,\alpha}(\overline{U},\R^{2\times 2}))$ as well. Then take $\dot{A}_t$ be the differential of $t\mapsto A_t$ in $C^{4,\alpha}$ topology. For $s$ small, by direct computation using \Cref{eqb1} we have:
        \begin{align}
            \mathrm{div}\left(\frac{A_{t+s}-A_t}{s}\n u_{t+s}\right)=-\mathrm{div}\left(\frac{A_t\n u_{t+s}}{s}\right)=-\mathrm{div}\left(A_t\frac{\n u_{t+s}-\n u_t}{s}\right).
        \end{align}
        Let $s\to 0$, we have that:
        \begin{equation}\label{eqb7}
            \mathrm{div}(\dot{A}_t\n u_t)=-\mathrm{div}(A_t\n\dot{u}_t).
        \end{equation} 
        Let $D^Qu_{s,t}=(u_{t+s}-u_t)/s$ be the difference quotient, then we have:
        \begin{align}
            \mathrm{div}(A_{t+s}\n(D^Qu_{s,t}-\dot{u}_t))&=-\frac{1}{s}\mathrm{div}(A_{t+s}\n u_t)-\mathrm{div}(A_{t+s}\n\dot{u}_t)\notag\\
            &=-\mathrm{div}\left(\left(\frac{A_{t+s}-A_t}{s}-\dot{A}_t\right)\n u_t+(A_{t+s}-A_t)\n \dot{u}_t\right).
        \end{align}
        Then apply the inhomogeneous version of \Cref{lemb1} to the equation above, we have:
        \begin{align}
            \Vert D^Qu_{s,t}-\dot{u}_t\Vert_{C^{5,\alpha}}\le C\left(\left\Vert \frac{A_{t+s}-A_t}{s}-\dot{A}_t\right\Vert_{C^{4,\alpha}}+\Vert A_{t+s}-A_t\Vert_{C^{4,\alpha}}+\Vert D^Qu_{s,t}-\dot{u}_t\Vert_{L^2}\right),
        \end{align}
        which converges to $0$ as $s\to 0$. This implies that $t\mapsto u_t$ is differentiable in $C^{5,\alpha}$ topology. The remaining task for us is to show that $t\mapsto \dot{u}_t$ is continuous. To this end, we apply the same idea to \Cref{eqb7} as follows:
        \begin{align}\label{eqb10}
            \mathrm{div}(A_t\n (\dot{u}_{t+s}-\dot{u}_t))&=\mathrm{div}((A_{t+s}-A_t)\n \dot{u}_{t+s}))+\mathrm{div}(\dot{A}_{t+s}\n u_{t+s})+\mathrm{div}(\dot{A}_t\n u_t)\notag\\
            &=\mathrm{div}((A_{t+s}-A_t)\n\dot{u}_{t+s}+(\dot{A}_{t+s}-\dot{A}_t)\n u_{t+s}+\dot{A}_t(\n u_t-\n u_{t+s})).
        \end{align}
        Thus by \Cref{lemb1} and \Cref{lemb2}, we have:
        \begin{align}
            \Vert\dot{u}_{t+s}-\dot{u}_t\Vert_{C^{5,\beta}}\le C(\Vert A_{t+s}-A_t\Vert_{C^{4,\beta}}+\Vert\dot{A}_{t+s}-\dot{A}_t\Vert_{C^{4,\beta}}+\Vert u_t-u_{t+s}\Vert_{C^{5,\beta}}+\Vert \dot{u}_{t+s}-\dot{u}_t\Vert_{L^2}),
        \end{align}
        which converges to $0$ as $s\to 0$.
    \end{proof}
    Apply this lemma to both $u$ and $v$, we get the proof of \Cref{lem5.1}. We remark that the reason we have to lower $C^{5,\alpha}$ to $C^{5,\beta}$ is that the last term $\n u_t-\n u_{t+s}$ in \Cref{eqb10} only lies in $C^{5,\beta}$.
    \bibliographystyle{alpha}
    \bibliography{bibliography-2-4}
\end{document}